%% file: arxiv_main.tex
\documentclass[english,ruled]{article}
\usepackage[T1]{fontenc}
\usepackage[latin9]{inputenc}
\usepackage{verbatim}
\usepackage{subcaption}
 \usepackage[linesnumbered,ruled]{algorithm2e}
\usepackage{amsmath}
\usepackage{amsthm}
\usepackage{etoolbox}
\usepackage{amssymb}
\usepackage{graphicx}
\usepackage{xcolor}
\usepackage{mathtools}
\usepackage{multicol}
\usepackage{multirow}
\usepackage{geometry}
\usepackage{booktabs}
\usepackage{enumitem}
\usepackage{setspace}
\makeatletter
\usepackage[toc,page,header]{appendix}
\usepackage{minitoc}
\usepackage{ifthen}
\newboolean{doublecolumn}
\setboolean{doublecolumn}{false}
\newboolean{arxiv}
\setboolean{arxiv}{true}
\usepackage{pgfplots}
\usepackage{pdflscape}
\usepackage[pagebackref=true]{hyperref}
\usepackage{cleveref}
\usepackage{authblk}
\usepackage{todonotes}
\pgfplotsset{compat=1.15}

\newboolean{proofinpaper}
\setboolean{proofinpaper}{false}

\newcommand{\inputinpaper}[1]{%
  \ifbool{proofinpaper}{\begin{proof}
						\input{#1}{}%
						\end{proof}}}

\newcommand{\inputinappendix}[1]{%
\ifbool{proofinpaper}{}{\input{#1}}%
}

\makeatletter
\renewcommand\paragraph{%
  \@startsection{paragraph}{4}{\z@}%
    {-1.5ex\@plus -0.5ex \@minus -0.2ex}
    {-.5em \@plus -.1em}%
    {\normalfont\normalsize\bfseries}%
}
\makeatother
\newcommand{\assign}{\coloneqq}
\newcommand{\nin}{\not\in}

\newcommand{\tmop}[1]{\ensuremath{\operatorname{#1}}}

\newcommand{\tmem}[1]{\textit{#1}}
\newcommand{\mathd}{\mathrm{d}}

\definecolor{deepskyblue}{rgb}{0.0, 0.75, 1.0}

\input{macros.tex}

\newcommand{\Comment}[1]{\hspace*{\fill}\textsf{$\triangleright$ ~#1}}

\renewcommand*\backref[1]{\ifx#1\relax \else (cited on #1) \fi}
\theoremstyle{plain}
\newtheorem{lem}{\protect\lemmaname}[section]
\theoremstyle{remark}
\newtheorem{rem}{\protect\remarkname}
\theoremstyle{plain}
\newtheorem{thm}{\protect\theoremname}[section]
\theoremstyle{plain}

\providecommand{\corollaryname}{Corollary}
\theoremstyle{plain}

\theoremstyle{plain}

\theoremstyle{plain}

\theoremstyle{plain}
\newtheorem{ques}{\protect\questionname}

\providecommand{\lemmaname}{Lemma}
\providecommand{\remarkname}{Remark}
\providecommand{\theoremname}{Theorem}
\providecommand{\examplename}{Example}
\providecommand{\propositionname}{Proposition}

\providecommand{\questionname}{Question}

\newcommand{\adam}{{\texttt{Adam}}}
\newcommand{\bfgs}{{\texttt{BFGS}}}
\newcommand{\lbfgs}{{\texttt{L-BFGS}}}
\newcommand{\ob}{\texttt{OSGM-Best}}
\newcommand{\ov}{\texttt{Vanilla OSGM-H}}
\newcommand{\oh}{\texttt{OSGM-H}}
\newcommand{\osr}{\texttt{OSGM-R}}
\newcommand{\os}{\texttt{OSGM}}
\newcommand{\ohb}{\texttt{OSGM-HB}}

\newcommand{\ohblook}{\texttt{Monotone Lookahead OSGM-HB}}

\newcommand{\ogd}{{\texttt{OGD}}}
\newcommand{\hdm}{{\texttt{HDM}}}
\newcommand{\hdmclassic}{{\texttt{Classic-HDM}}}

\newcommand{\agd}{{\texttt{AGD}}}
\newcommand{\osgm}{{\texttt{OSGM}}}
\newcommand{\osgmrx}{\texttt{OSGM-R}}

\newcommand{\osgmhandsonhx}{\texttt{Monotone Lookahead OSGM-H}}
\newcommand{\osgmhandsoffhx}{\texttt{Monotone OSGM-H}}

\newcommand{\Mlook}{\mathcal{M}_{\text{look}}}

\newcommand{\Fhdm}{\mathcal{F}_{{\hdm}}}
\newcommand{\Fosgm}{\mathcal{F}_{{\osgm}}}

\crefdefaultlabelformat{#2\textbf{#1}#3} 
\crefname{section}{\textbf{section}}{\textbf{sections}}
\Crefname{section}{\textbf{Section}}{\textbf{Sections}}
\crefname{thm}{\textbf{Theorem}}{\textbf{theorems}}
\Crefname{thm}{\textbf{Theorem}}{\textbf{Theorems}}
\crefname{lem}{\textbf{Lemma}}{\textbf{lemmas}}
\Crefname{lem}{\textbf{Lemma}}{\textbf{Lemmas}}
\crefname{prop}{\textbf{proposition}}{\textbf{propositions}}
\Crefname{prop}{\textbf{Proposition}}{\textbf{Propositions}}
\crefname{algorithm}{\textbf{algorithm}}{\textbf{algorithms}}
\Crefname{algorithm}{\textbf{Algorithm}}{\textbf{Algorithms}}
\crefname{coro}{\textbf{Corollary}}{\textbf{corollaries}}
\Crefname{coro}{\textbf{Corollary}}{\textbf{corollaries}}
\crefname{definition}{\textbf{Definition}}{\textbf{definitions}}
\Crefname{definition}{\textbf{Definition}}{\textbf{definitions}}
\crefname{table}{\textbf{Table}}{\textbf{tables}}
\Crefname{table}{\textbf{Table}}{\textbf{tables}}
\crefname{figure}{\textbf{Figure}}{\textbf{figures}}
\Crefname{figure}{\textbf{Figure}}{\textbf{figures}}
\crefname{exple}{\textbf{Example}}{\textbf{examples}}
\Crefname{exple}{\textbf{Example}}{\textbf{examples}}

\begin{document}

\title{Gradient Methods with Online Scaling \\ \vspace{10pt}
\Large Part II. Practical Aspects \\
}

\author[1]{Ya-Chi Chu\thanks{ycchu97@stanford.edu, equal contribution}}
\author[2]{Wenzhi Gao\thanks{gwz@stanford.edu, equal contribution}}
\author[2,3]{Yinyu Ye\thanks{yyye@stanford.edu}}
\author[2,3]{Madeleine Udell\thanks{udell@stanford.edu}}
\affil[1]{Department of Mathematics, Stanford University}
\affil[2]{ICME, Stanford University}
\affil[3]{Department of Management Science and Engineering, Stanford University}

\maketitle

\begin{abstract}
Part I of this work \cite{gao2025gradient} establishes online scaled gradient methods ({\os})\footnote{This paper extends two previous works \cite{gao2024gradient} and \cite{chu2025provable}.}, a framework that utilizes online convex optimization to adapt stepsizes in gradient methods. This paper focuses on the practical aspects of {\os}. We leverage the {\os} framework to design new adaptive first-order methods and provide insights into their empirical behavior. The resulting method, {\ob}, matches the performance of quasi-Newton variants while requiring less memory and cheaper iterations. We also extend {\os} to nonconvex optimization and outline directions that connect {\os} to existing branches of optimization theory and practice.
\end{abstract}

\input{sec_intro.tex}
\input{sec_osgmprac.tex}
\input{sec_nonconvex.tex}
\input{sec_exp.tex}

\input{sec_directions.tex}

\input{conclusions.tex}

\renewcommand \thepart{}
\renewcommand \partname{}

\bibliography{ref.bib}
\bibliographystyle{plain}

\doparttoc
\faketableofcontents
\part{}

\newpage
\appendix

\addcontentsline{toc}{section}{Appendix}
\part{Appendix}
\parttoc

\newpage
\input{app_prac.tex}

\input{app_nonconvex.tex}
\input{app_extensions.tex}

\end{document}

%% file: macros.tex
\global\long\def\vertiii#1{\left\vert \kern-0.25ex  \left\vert \kern-0.25ex  \left\vert #1\right\vert \kern-0.25ex  \right\vert \kern-0.25ex  \right\vert }%

\global\long\def\diam{\mathrm{diam}}%

\global\long\def\argmin{\operatornamewithlimits{arg\,min}}%

\global\long\def\diag{\mathrm{Diag}}%

\global\long\def\and{\mathrm{and}}%

\global\long\def\Rbb{\mathbb{R}}%

\global\long\def\Acal{\mathcal{A}}%

\global\long\def\Bcal{\mathcal{B}}%

\global\long\def\Fcal{\mathcal{F}}%

\global\long\def\Mcal{\mathcal{M}}%

\global\long\def\Ocal{\mathcal{O}}%

\global\long\def\Pcal{\mathcal{P}}%

\global\long\def\Scal{\mathcal{S}}%

\global\long\def\Xcal{\mathcal{X}}%

%% file: sec_intro.tex
\section{Introduction} \label{sec:intro}

Consider gradient descent applied to a smooth convex problem $f^\star \coloneqq \min_{x \in
\mathbb{R}^n} f (x)$:
\begin{equation} \label{eqn:osgm-update}
	x^{k + 1} = x^k - P_k \nabla f (x^k),
\end{equation}
where $P_k \in \mathbb{R}^{n \times n}$ is a matrix stepsize. The first part of this work {\cite{gao2025gradient}} introduces {\os}, a framework that
uses online learning to adjust stepsize $\{ P_k \}$ on the fly
and achieves problem-dependent acceleration.\\

We briefly summarize the idea of {\os} introduced in Part I \cite{gao2025gradient}. The analyses of algorithms for smooth optimization often show that each iteration yields some improvement. For instance, in linearly convergent methods, a common improvement is a contraction of the form $\tfrac{f (x^{k + 1}) - f^{\star}}{f (x^k) - f^{\star}} \leq 1 - \tfrac{1}{\kappa}$. Here $\kappa \geq 1$ often denotes some condition number. If this contraction ratio bound holds uniformly for all $k$, chaining the per iteration improvements over $K$ steps yields the rate
\[ \tfrac{f (x^{K + 1}) - f^{\star}}{f (x^1) - f^{\star}} = \textstyle \prod_{k = 1}^K
   \tfrac{f (x^{k + 1}) - f^{\star}}{f (x^k) - f^{\star}} \leq ( 1 -
   \tfrac{1}{\kappa} )^K . \]

Many proofs first establish improvement at every iteration and then chain these improvements to prove convergence. 
However, this approach requires monotonic improvement. {\os} instead chains the contraction first:
\begin{equation}
	\tfrac{f (x^{K + 1}) - f^{\star}}{f (x^1) - f^{\star}} = \textstyle \prod_{k = 1}^K
   \tfrac{f (x^{k + 1}) - f^{\star}}{f (x^k) - f^{\star}} \leq (
   \tfrac{1}{K} \textstyle \sum_{k = 1}^K \tfrac{f (x^{k + 1}) - f^{\star}}{f (x^k) -
   f^{\star}} )^K . \label{eqn:ratio-bound}
\end{equation}
To guarantee convergence, it suffices to have the {\tmem{average}} contraction
ratio $\tfrac{1}{K} \textstyle \sum_{k = 1}^K \tfrac{f (x^{k + 1}) - f^{\star}}{f (x^k) -
f^{\star}}$ be sufficiently small. In other words, convergence is possible even if some iterations worsen the objective. When the iterates $\{x^k\}$ are generated via the update \eqref{eqn:osgm-update}, the contraction ratio attributed to one iteration, starting from $x$, can be written as $r_x (P) = \tfrac{f(x - P \nabla f (x)) - f^{\star}}{f (x) - f^{\star}}$, which depends on the stepsize matrix. In view of \eqref{eqn:ratio-bound}, it suffices to select $\{ P_k \}$ sequentially to minimize $\tfrac{1}{K} \textstyle \sum_{k = 1}^K r_{x^k} (P_k)$, the average contraction ratio. This perspective formulates stepsize selection as a sequential decision-making problem, a central topic in online learning \cite{orabona2019modern,hazan2016introduction}. 
If $\{ P_k \}$ are chosen by a no-regret online learning algorithm, we obtain the guarantee:
$\tfrac{1}{K} \textstyle \sum_{k = 1}^K r_{x^k} (P_k) \leq \tfrac{1}{K} \textstyle \sum_{k = 1}^K
r_{x^k} (\hat{P}) + o (1)$ for any benchmark stepsize $\hat{P}$, including the optimal preconditioner \cite{gao2023scalable,qu2024optimal}. The relation further yields the convergence guarantee
\[ \tfrac{f (x^{K + 1}) - f^{\star}}{f (x^1) - f^{\star}} \leq (
   \tfrac{1}{K} \textstyle \sum_{k = 1}^K r_{x^k} (\hat{P}) + o (1) )^K \text{~~~for any $\hat{P}$.} \]
More generally, {\os} models stepsize selection in a gradient-based method as a sequential decision-making problem between a stepsize scheduling agent and a landscape agent. At each iteration:
\begin{enumerate}[leftmargin=14pt]
  \item Scheduler makes decision $P_k$ and proposes an update $x^{k + 1 / 2} =
  x^k - P_k \nabla f (x^k)$.
  
  \item Landscape chooses the next iterate $x^{k + 1} =\mathcal{M} (x^k, x^{k
  + 1/2})$ and provides feedback $\ell_{x^k}
  (P_k)$ to the scheduler.
  
  \item Scheduler updates the stepsize using an online learning algorithm
  $\mathcal{A} (P_k, \{ \ell_{x^j} \}_{j \leq k})$.
\end{enumerate}
Different instantiations of $(\ell_x, \mathcal{M}, \mathcal{A})$ yield different algorithm variants, each with its own convergence properties. Part I of this work develops the theoretical foundations of {\os}, providing specific instantiations that achieve both global and local convergence guarantees.

\subsection{Contributions}

Building on the theoretical foundations established in Part I, the second part focuses on the practical aspects of {\os} and outlines several future directions.

\paragraph{Insights into algorithm behavior.}
While Part I focuses on complexity analysis, this work investigates the empirical behavior of {\os}, with a particular focus on {\oh} \cite{gao2025gradient}. We develop a more detailed analysis that explains the empirical behavior of {\oh} and other hypergradient-type methods \cite{gunes2018online}. These insights lead to concrete algorithmic recommendations. 
	
\paragraph{A practical variant of {\os}.}
Part I developed {\os} variants based on gradient descent. 
This paper shows how to use the same ideas to choose stepsizes in other gradient-based methods, such as heavy-ball momentum.
The strategies above, when combined with heavy-ball momentum, give rise to {\ob}: a robust variant of {\os} that performs competitively with {\lbfgs} \cite{liu1989limited} while requiring less memory and cheaper iterations.

\paragraph{Nonconvex optimization.}
Although {\os} is grounded in online convex optimization, we demonstrate that the analysis can be extended to smooth nonconvex problems by introducing regularization in the stepsize space. Under the Polyak-{\L}ojasiewicz (PL) condition, the resulting {\os} variants provably adapt to locally convex regions of the landscape. Experiments on standard nonconvex benchmarks support these theoretical claims.

\paragraph{Connection to other optimization techniques.}
We highlight connections between {\os} and a range of ideas in optimization, including the Barzilai-Borwein step \cite{barzilai1988two}, accelerated gradient descent \cite{d2021acceleration}, proximal gradient methods \cite{beck2017first}, and the performance estimation problem (PEP) \cite{upadhyaya2025automated}. These connections highlight the potential of {\os} as a valuable addition to the modern optimization toolbox.

\subsection{Notations}

We use $\| \cdot \|$ to denote the Euclidean norm. Letters $A, a$ denote matrices and scalars. 
$\| A \|_F \assign \sqrt{\sum_{i j} a_{i j}^2}$ denotes the matrix Frobenius norm. We use $\langle \cdot, \cdot
\rangle$ for Euclidean or Frobenius inner product.
$\Pi_{\mathcal{P}} [\cdot]$ denotes the orthogonal projection onto a closed convex set $\mathcal{P}$. 
We use $\mathcal{X}^{\star} = \{ x : f (x) = f^{\star} \}$ to denote the optimal
set of $f$ and $\Scal^\star$ to denote the set of stationary points of $f$; $\tmop{dist} (P, \mathcal{P}) \assign \| P - \Pi_{\mathcal{P}}
[P] \|_F$ denotes the distance between a point $P$ and a closed convex set
$\mathcal{P}$; $\tmop{diam} (\mathcal{P}) = \max_{X, Y \in \mathcal{P}}  \| X
- Y \|_F$ denotes the diameter of set $\Pcal$ in Frobenius norm. 
We use superscript $x^k$ to index algorithm iterates and subscript $P_k$ to index
the stepsize sequence. A function $f$ satisfies $\mu$-PL condition if $f(x) - f^\star \leq \frac{1}{2\mu} \|\nabla f(x) \|^2$.
To emphasize that the feedback functions have stepsize $P$ as their argument, we will use $h_x(P)$ instead of $h_x$ in most contexts.

\paragraph{Structure of the paper.}This paper is organized as follows. \Cref{sec:prac} provides insights into the convergence behavior of {\oh} observed in the literature. These insights, combined with a heavy-ball momentum enhancement,
constitute {\ob}, a variant that yields competitive practical performance.
\Cref{sec:nonconvex} generalizes {\os} to smooth nonconvex optimization. \Cref{sec:exp}
presents numerical experiments on both convex and nonconvex optimization
problems. \Cref{sec:extensions} outlines several future directions and shows the connection between {\os} and existing optimization theory.

%% file: sec_osgmprac.tex
\section{Algorithmic behavior and ingredients of practical {\osgm}} \label{sec:prac}

In this section, we motivate the key ingredients of {\ob} that improve the practical convergence of {\os}.  
We provide insights into each ingredient through either rigorous analysis or illustrative examples. Our analysis focuses on {\oh} (\Cref{alg:osgmhandsonhx}), {\os} using hypergradient feedback, although most of the insights
extend naturally to other variants of {\os}. The hypergradient feedback \cite{gao2025gradient} takes the form

\[ h_x(P) = \tfrac{f(x - P\nabla f(x)) - f(x)}{\|\nabla f(x)\|^2}.\]

\begin{figure}[!h]
\begin{minipage}[t]{0.56\textwidth}
\begin{algorithm}[H]
{\textbf{input:}  $x^1$, $P_1 \in \Rbb^{n\times n}$, online gradient stepsize $\eta > 0$}\\
\For{$k = 1, 2, \dots$}{
      $x^{k + 1/2} = x^{k} - P_k \nabla f(x^k)$\\
      Let $x^{k+1}$ satisfy $f(x^{k+1}) \leq \min \{ f(x^{k + 1 / 2}), f(x^k) \}$ \label{alg:oh-monotone} \\
      $P_{k+1} = P_k - \eta \nabla h_{x^k}(P_k)$ }
\caption{{\osgmhandsoffhx}\label{alg:osgmhandsonhx}}
\end{algorithm}
\end{minipage}
\quad
\begin{minipage}[t]{0.4\textwidth}
\begin{algorithm}[H]
{\textbf{input:}  $x^1$, $P_1 \in \Pcal = \Rbb^{n\times n}$,  $\eta > 0$}\\
\For{$k = 1, 2, \dots$}{
      $P_{k+1} = P_k - \eta \nabla h_{x^k}(P_k)$\\
      $x^{k + 1} = x^{k} - P_{k+1} \nabla f(x^k)$ }
\caption{{\hdmclassic}\label{alg:classic}}
\end{algorithm}
\end{minipage}
	
\end{figure}

The feedback $h_x(P)$ is motivated by the classical descent lemma $f(x - \frac{1}{L} \nabla f(x) ) \leq - \frac{1}{2L} \|\nabla f(x)\|^2$ and is closely related to the classic hypergradient descent heuristic (\Cref{alg:classic}) studied in the literature \cite{gunes2018online,rubio2017convergence}. However, there are two key differences: the monotonicity condition imposed on \Cref{alg:oh-monotone} and the order in which $P_k$ and $x^{k}$ are updated. These seemingly subtle distinctions turn out to be critical for an effective variant of {\oh} in both theory and in practice. In this section, we detail our design and address the following aspects:

\begin{itemize}[leftmargin=10pt]
  \item \emph{Monotone landscape action for {\oh}.}
  Is it necessary to enforce monotonicity in {\oh}?
  \item \emph{Limitations of hypergradient descent in the literature.}  
  How does the order of the updates in {\oh} matter?
  
  \item \emph{Heavy-ball momentum and modified feedback.}
  Can momentum accelerate the convergence of {\oh}?
\end{itemize}
After discussing each ingredient in turn, we combine them into one recipe: {\ob}.

\input{sec_practical_monotone.tex}
\input{sec_practical_precient.tex}
\input{sec_practical_hbfeedback.tex}
\input{sec_osgm_best.tex}

%% file: sec_practical_monotone.tex
\subsection{Monotone landscape action for hypergradient feedback}

\begin{figure}[h]
\centering
\includegraphics[scale=0.282]{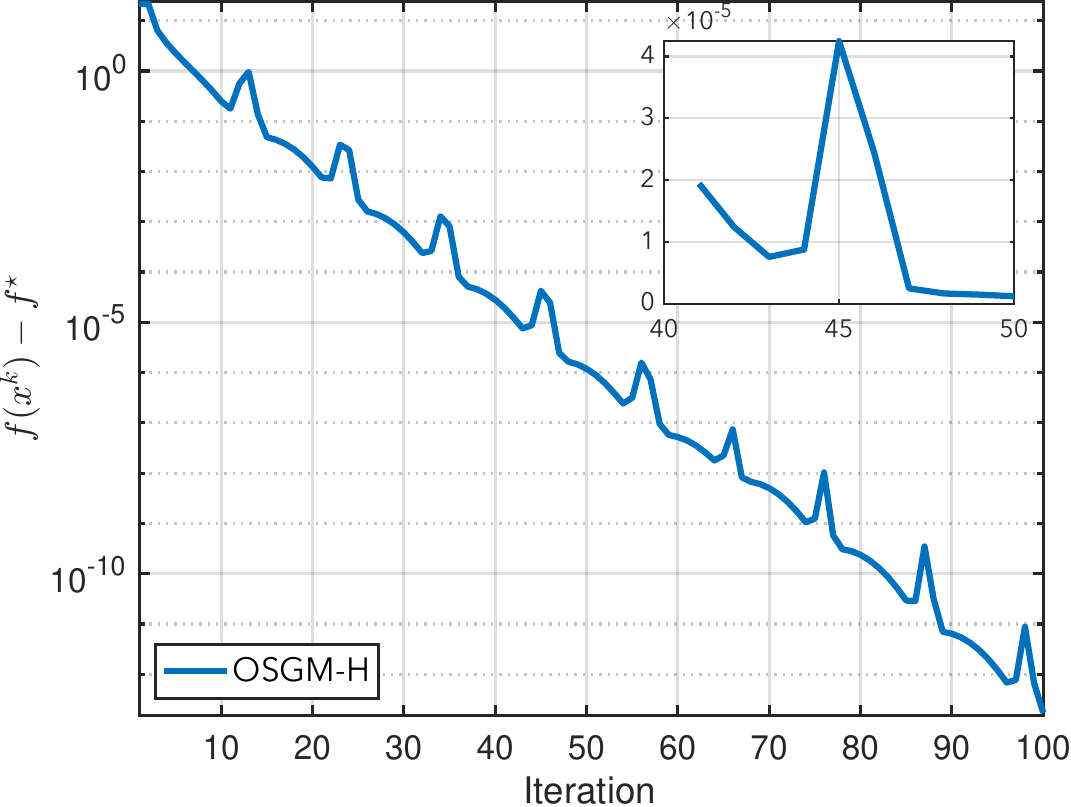}\qquad
\includegraphics[scale=0.16]{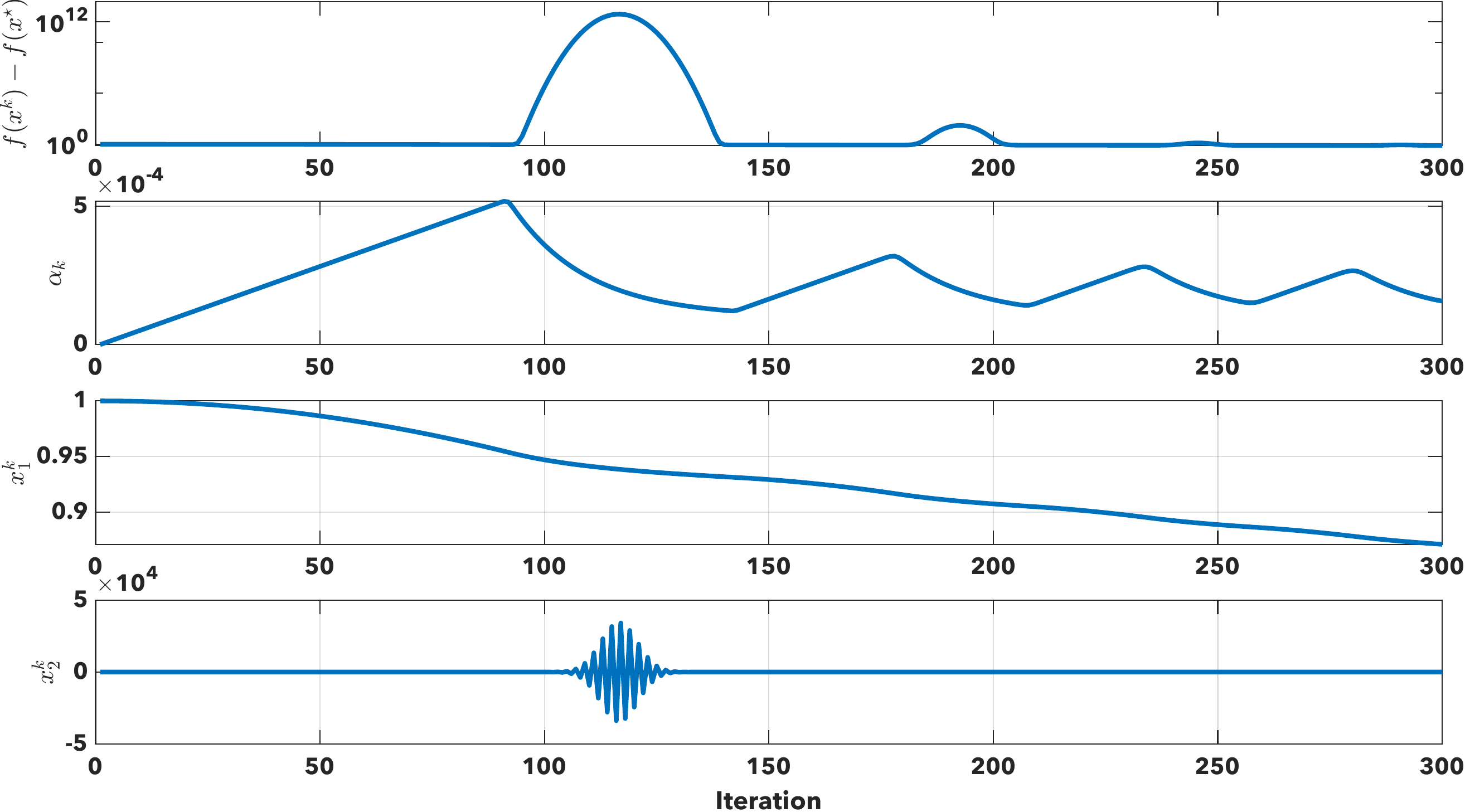}
 \caption{Spiky behavior of {\ov} on a quadratic function. \label{fig:spiky}}
\end{figure}

The convergence analysis for {\oh} in Part I \cite{gao2025gradient} always requires a monotone landscape action to ensure at least $f (x^{k + 1}) \leq f (x^k)$.
In practice, {\oh} is still observed to converge without enforcing monotonicity, but it exhibits periodic ``spikes'' in the objective (\Cref{fig:spiky} (left)). 
Similar spiky patterns have also been observed in other methods, such as heavy-ball momentum \cite{danilova2020non} and the Barzilai-Borwein method {\cite{fletcher2005barzilai}}.

In this section, we show that this requirement is not an artifact of analysis. Actually, one can construct an example for which the spike of {\oh} can be arbitrarily high and can happen arbitrarily late. This precludes the possibility of establishing a standard convergence rate in terms of $K$ without the monotone landscape action. 

\paragraph{Example.}
Consider a simple 2D quadratic function parametrized by $\kappa \gg 1$:
\[ \min_{x = (x_1, x_2) \in \mathbb{R}^2} f (x) = \tfrac{1}{2} x_1^2 +
   \tfrac{\kappa}{2} x^2_2 . \]
The objective $f$ is $\kappa$-smooth and $1$-strongly convex with condition
number $\kappa$ and  $x^{\star} = (0, 0)$. The minimax optimal stepsize for this problem is $\alpha^\star = \frac{2}{1 + \kappa} \approx \frac{2}{\kappa}$, and stepsizes larger than $\alpha^\star$ can lead to divergence. Scalar version of {\oh} (without enforcing monotonicity) reads
\begin{align}
  \left(\begin{smallmatrix}
    x^{k + 1}_1\\
    x^{k + 1}_2
  \end{smallmatrix}\right) & =  \left(\begin{smallmatrix}
    (1 - \alpha_k) x^k_1\\
    (1 - \kappa \alpha_k) x^k_2
  \end{smallmatrix} \right), \nonumber\\
  \alpha_{k + 1} & = \alpha_k + \eta \Big[ \underbrace{\tfrac{(x^k_1)^2}{(x^k_1)^2 +
  \kappa^2 (x^k_2)^2} (1 - \alpha_k) + \tfrac{\kappa^2 (x^k_2)^2}{(x^k_1)^2 +
  \kappa^2 (x^k_2)^2} (1 - \kappa \alpha_k)}_{\Delta_k} \Big] . \label{eqn:mono-act-alpha-update}
\end{align}

The sign of stepsize change $\Delta_k$ is determined by both $1 - \alpha_k$ and $1 - \kappa \alpha_k$. When $\alpha_k > \frac{1}{\kappa}$, the value $1 - \kappa \alpha_k < 0$ from the second coordinate will try reducing the stepsize, while $1 - \alpha_k$ from the first coordinate tends to increase the stepsize until $\alpha_k \geq 1 \gg \alpha^\star$. However, these two terms are weighted by quantities that depend on the location of $x^k$ on the landscape. In particular, if $x^k_1 \gg \kappa x^k_2$, then $\tfrac{\kappa^2 (x^k_2)^2}{(x^k_1)^2 +
  \kappa^2 (x^k_2)^2}$ is close to $0$, and the value $1 - \kappa \alpha_k$ from the second coordinate will be disregarded.  It can therefore lead to overly large stepsize and divergence.  Suppose {\oh} is initialized such that $x^k_1 \gg \kappa x^k_2$. Then each spike in \Cref{fig:spiky} can be decomposed into three phases. \\

\emph{Myopic stepsize increase.} When $x^k_1 \gg \kappa x^k_2$, $\Delta_k \approx 1 - \alpha_k$, and  {\oh} drives the stepsize to $1$ to accelerate the convergence in the direction of $x_1$. However, the magnitude of $x_2$ explodes since $| 1 - \kappa \alpha_k | \gg 1$. \\

\emph{Divergence and recovery.} Though the divergence of $x_2$ starts from a small initial value, it finally becomes large enough that the second term $\tfrac{\kappa^2 (x^k_2)^2}{(x^k_1)^2 +
  \kappa^2 (x^k_2)^2}$ in \eqref{eqn:mono-act-alpha-update} dominates, decreasing $\alpha_k$. 
However, {\oh} requires many iterations to decrease $\alpha_k$ to $\tfrac{2}{\kappa}$, and $x^k_2$ continues to increase during this phase, resulting in a spike in the objective. \\

\emph{Periodic spiky behavior.}
After {\oh} decreases the stepsize below $\tfrac{2}{\kappa}$, the second coordinate $x^k_2$ decreases sharply. However, $x^k_2$ decreases much faster than $x^k_1$ does since $1 - \kappa \alpha_k \ll 1 - \alpha_k$ if $\alpha_k = \mathcal{O}(\tfrac{1}{\kappa})$. 
Then {\oh} returns to the state $x_1^k \gg \kappa x_2^k$, and repeats the same divergence-recovery loop.
Notice that the magnitude $|x_1|$ consistently decreases (\Cref{fig:spiky} (right)), since $\alpha_k$ remains less than $1$ and thus $|1 - \alpha_k| < 1$. Hence, the algorithm ultimately converges after several loops.
But this analysis also precludes a standard convergence proof, as we can pick $\delta \rightarrow 0$ to construct a bad starting point $(1, \delta)$ for any (sufficiently large) iteration $K$ so that a spike happens around $K$ and $f(x^K) > f(x^1)$. \\

To quantify the spiky behavior, 
suppose {\oh} starts from $(x^1_1, x^1_2) = (1, \delta)$ for some
$\delta > 0$ to be determined later, and the stepsize is initialized by
$\alpha_1 = 0$. \Cref{lem:osgmhx-divergence} shows that {\oh}
increases stepsize by at least a constant factor of $\eta \leq
\tfrac{1}{\kappa}$ when $\alpha_k$ is small enough and $x^k_1 \gg \kappa
x^k_2$.

\begin{lem} \label{lem:osgmhx-divergence}
  Let $\kappa \geq 2$. If $\alpha_k \leq \tfrac{1}{2}$ and $| x^k_1 | \geq \sqrt{2} \kappa^{3 / 2} |
  x^k_2 |$, then $\alpha_{k + 1} \geq \alpha_k + \tfrac{\eta}{4}$.
\end{lem}

\input{mpproofs/proof-lem-osgmhx-divergence.tex}

Assume $\kappa$ is large enough to have a constant $a$ such that
$\tfrac{2}{\kappa} \ll \tfrac{a}{2} < \tfrac{1}{2}$, say $a = \tfrac{1}{3}$.
Note that if the condition in \Cref{lem:osgmhx-divergence} is always satisfied,
{\oh} takes at most $\lceil \tfrac{4 a}{\eta} \rceil$ iterations to increase the
stepsize from $\alpha_1 = 0$ to at least $a$. 
During the stepsize increasing period, the magnitude $| x^k_2 |$ is at most
expanded by a factor $| 1 - \kappa \alpha_k | < \kappa - 1$ since $\alpha_k <
1$. Therefore, we can choose $\delta > 0$ small enough such that the iterates
$\{ (x^k_1, x^k_2) \}_{k \leq K}$ always satisfy the condition in
\Cref{lem:osgmhx-divergence} but $(x^{K + 1}_1, x^{K + 1}_2)$ violates the same
condition. Then {\oh} starts to decrease the stepsize after iteration $K + 1$, but the stepsize $\alpha_{K + 1} > a \gg \tfrac{2}{\kappa}$
and thus $x_2$ continues to expand, causing a spike in the objective.
{\oh} takes $\Theta ( \tfrac{1}{\eta \kappa})$ iterations to decrease the stepsize from $a$ to $\tfrac{a}{2}$ since
for $\alpha_k \leq a$, every iteration can decrease the stepsize by at most $\eta
(\kappa a - 1)$:
\[ \alpha_{k + 1} - \alpha_k = \eta \tfrac{(1 - \alpha_k) (x^k_1)^2 + \kappa^2
   (1 - \kappa \alpha_k) (x^k_2)^2}{(x^k_1)^2 + \kappa^2 (x^k_2)^2} \geq -
   \eta (\kappa a - 1) . \]
Every iteration during this phase expands $x_2$ by at least a factor $| 1 -
\alpha_k \kappa | \geq \tfrac{a \kappa}{2} - 1 = \tfrac{\kappa}{6} - 1$ and the spike grows as $\Omega ( \kappa^{\frac{1}{\eta \kappa}})$. Hence, given any $K$, we can pick $\delta, \eta$, and $\kappa$ so that a spike appears around $K$ (\Cref{fig:spiky} (right)).

%% file: proofs/proof-lem-osgmhx-divergence.tex
The assumption $| x_1^k | \geq \sqrt{2} \kappa^{3 / 2} | x_2^k |$ and that $\kappa \geq 2 $ implies
\begin{equation} \label{eqn:pf-divergence-1}
  (x_1^k)^2 \geq 2 \kappa^3 (x_2^k)^2 \geq 2 \kappa^2 (\kappa - 2) (x_2^k)^2.
\end{equation}
Then we deduce that
\begin{align}
  \tfrac{\alpha_{k + 1} - \alpha_k}{\eta} ={} & \tfrac{(x_1^k)^2 (1 - \alpha_k)
  + \kappa^2 (1 - \kappa \alpha_k) (x_2^k)^2}{(x_1^k)^2 + \kappa^2 (x_2^k)^2}
  \tag{by \eqref{eqn:mono-act-alpha-update}}\\
  \geq{} & \tfrac{1}{4} \tfrac{2 (x_1^k)^2 - 2 \kappa^2 (\kappa - 2) (x_2^k)^2}{(x_1^k)^2 +
  \kappa^2 (x_2^k)^2} \tag{by $\alpha_k \leq \tfrac{1}{2}$}\\
  \geq{} & \tfrac{1}{4} \tfrac{(x_1^k)^2}{(x_1^k)^2 + \kappa^2 (x_2^k)^2}
  \geq \tfrac{1}{4}.
  \tag{by \eqref{eqn:pf-divergence-1}}
\end{align}
This completes the proof.

%% file: sec_practical_precient.tex
\subsection{Limitations of \hdmclassic} \label{sec:prescient}

The classic version of hypergradient descent
(\hdmclassic, \Cref{alg:classic}) \cite{gunes2018online,rubio2017convergence} uses a different order for the primal update and stepsize update compared to {\oh}:
\[ \begin{array}{ll}
     {\hdmclassic} & \oh\\
     P_{k + 1} = P_k - \eta \nabla h_{x^k} (P_k) & x^{k + 1} = x^k - P_k
     \nabla f (x^k)\\
     x^{k + 1} = x^k - P_{k + 1} \nabla f (x^k) \qquad & P_{k + 1} = P_k -
     \eta \nabla h_{x^k} (P_k)
   \end{array} \]

\begin{figure}[h]
  \centering
  \includegraphics[width=0.5\linewidth]{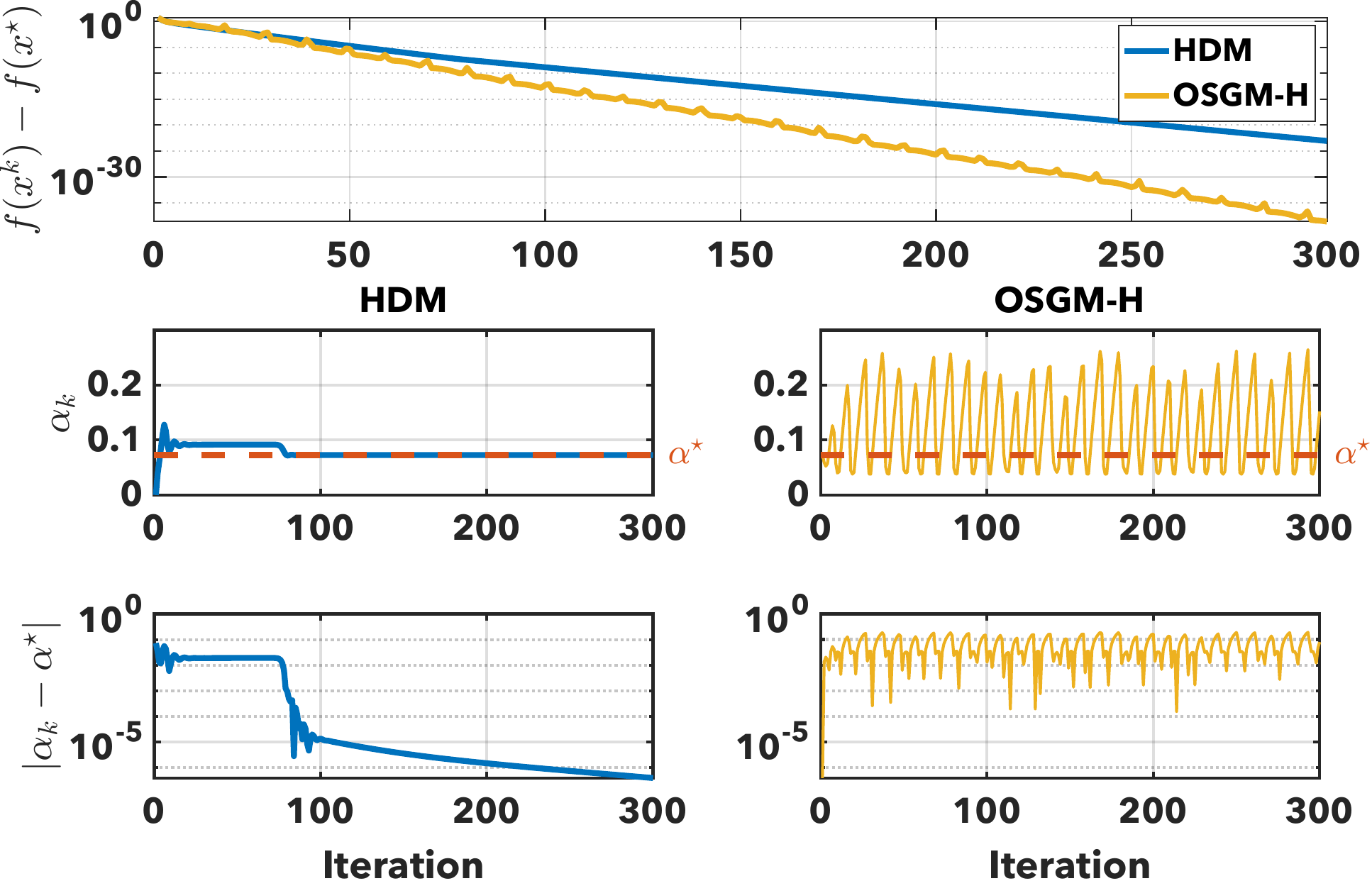}
  \includegraphics[width=0.45\linewidth]{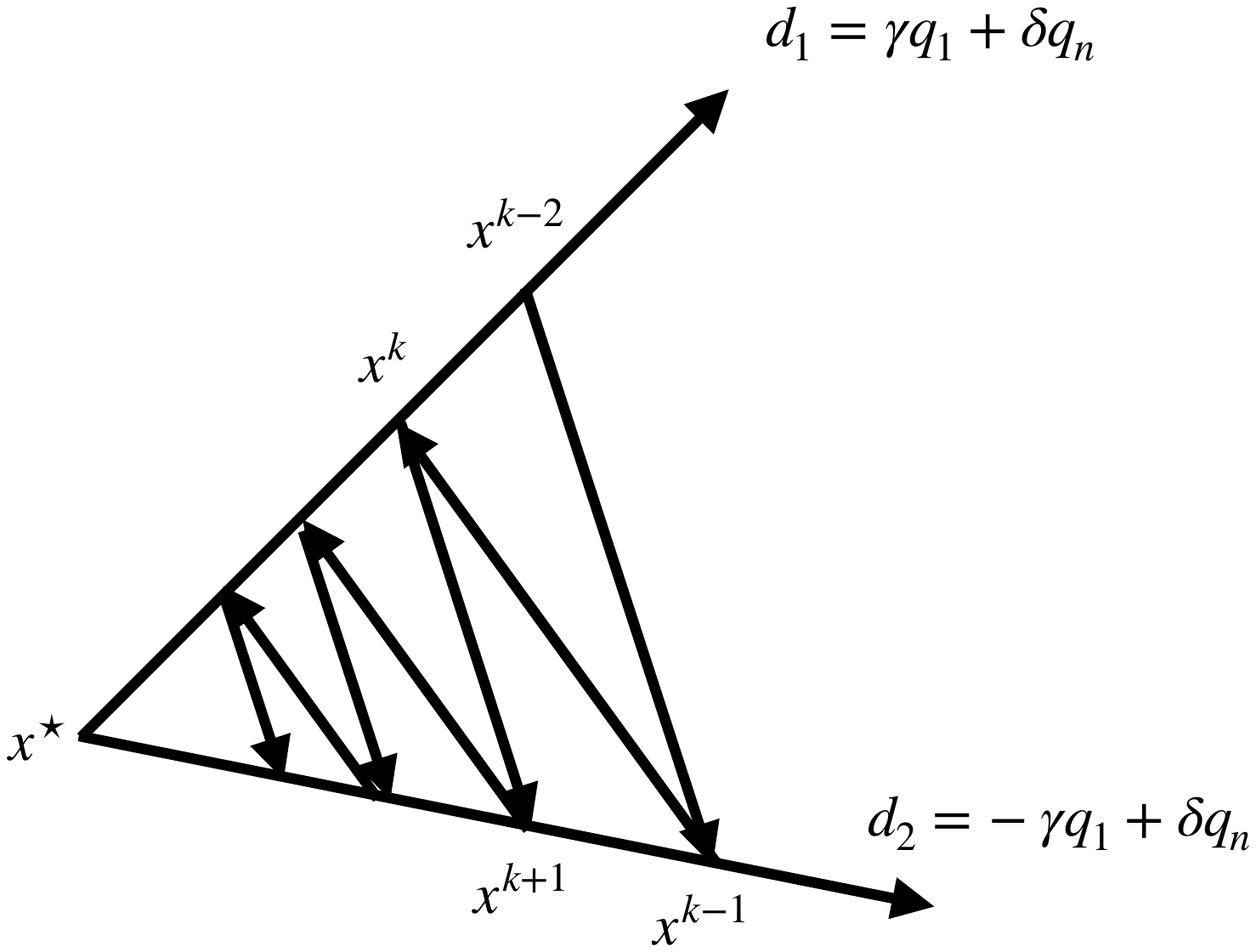}
  \caption{
  Left: stepsizes $\{ \alpha_k \}$ converge to $\alpha^\star = \tfrac{2}{L + \mu}$ in {\hdmclassic} and oscillate around $\alpha^\star = \tfrac{2}{L + \mu}$ in {\os}. Right: as $\alpha_k$ converges, {\hdmclassic} will converge to $x^\star$ along two directions.}
  \label{fig:dynamics}
\end{figure}

To facilitate our analysis, for this section only, we drop the monotone landscape action for {\oh}. This simplification does not compromise the intuitions we want to convey.  {\hdmclassic} updates the stepsize using the feedback function at $x^k$ first
and then applies the newly obtained stepsize $P_{k + 1}$ to the gradient step
on $x^k$; while {\oh} directly applies the stepsize $P_k$. From an online learning perspective, the decisions in {\hdmclassic} utilize information from the future. Intuitively, it should yield better performance. However, the performance of {\hdmclassic} is inferior to {\oh}.  \Cref{fig:dynamics} (left) illustrates the behavior of {\hdmclassic} applied to a quadratic function. There are two observations:
\begin{enumerate}[leftmargin=15pt]
\item {\oh} converges faster than {\hdmclassic} under the same parameter settings.
\item The stepsize $\{\alpha_k\}$ in {\hdmclassic} converges, while $\{\alpha_k\}$ keeps oscillating in {\oh}.
\end{enumerate}
The convergence of $\{\alpha_k\}$ was observed by \cite{gunes2018online} and was partially explained in \cite{chu2025provable} in the case where $P$ is a full-dimensional matrix. This section further characterizes this phenomenon on a toy quadratic problem:
\[ \min_{x\in\Rbb^n}  \tfrac{1}{2} \langle x, A x \rangle, \]

where $A \succ 0$ has distinct eigenvalues.
We will demonstrate that $\{\alpha_k\}$ in {\hdmclassic} stabilizes around the minimax optimal stepsize $\alpha^\star = \tfrac{2}{L + \mu}$. Besides, the iterates $\{x^k\}$ generated by {\hdmclassic} alternates between two extremal eigen directions of $A$ and lead to slow convergence. On the other hand, in {\oh}, $\{\alpha_k\}$ does not stabilize and keeps adapting.\\

\textit{Dynamical system formulation and orbits.}
We can further simplify the analysis of {\hdmclassic} and
{\oh} by working on the basis that consists of the eigenvectors of $A$. Let $A = Q
\Lambda Q^{\top}$ be the eigendecomposition. Transform the variables to $z
\assign Q x$. Then {\hdmclassic} and {\oh} with scalar stepsize
$P_k = \alpha_k I$ read as
\[ \begin{array}{ll}
     {\hdmclassic} & \oh\\
     \alpha_{k + 1} = \alpha_k + \eta - \alpha_k \eta \tfrac{\langle z^k,
     \Lambda^3 z^k \rangle}{\langle z^k, \Lambda^2 z^k \rangle} \qquad & z^{k
     + 1} = (I - \alpha_k \Lambda) z^k\\
     z^{k + 1} = (I - \alpha_{k + 1} \Lambda) z^k \qquad & \alpha_{k + 1} =
     \alpha_k + \eta - \alpha_k \eta \tfrac{\langle z^k, \Lambda^3 z^k
     \rangle}{\langle z^k, \Lambda^2 z^k \rangle}
   \end{array} \]
Since the problem is homogeneous, we consider the dynamical system whose state is given by the stepsize $\alpha$ and the
direction $\hat{z} \assign \tfrac{z}{\| z \|}$:
\begin{equation}\label{eqn:dynamical-system}
	{\Fhdm} (\alpha, \hat{z}) \assign \big( \alpha^+
   (\hat{z}), \tfrac{(I - \alpha^+ (\hat{z}) \Lambda) \hat{z}}{\| (I -
   \alpha^+ (\hat{z}) \Lambda) \hat{z} \|} \big) \quad \text{and} \quad
   {\Fosgm} (\alpha, \hat{z}) \assign \big(
   \alpha^+ (\hat{z}), \tfrac{(I - \alpha \Lambda) \hat{z}}{\| (I - \alpha
   \Lambda) \hat{z} \|} \big)
\end{equation}
where the updated stepsize $\alpha^+ (\hat{z}) \assign \alpha + \eta - \alpha \eta \tfrac{\langle
\hat{z}, \Lambda^3 \hat{z} \rangle}{\langle \hat{z}, \Lambda^2 \hat{z}
	\rangle}$ starting at state $(\alpha, \hat{z})$ is the same in both algorithms.
Assume $A$ has distinct eigenvalues for simplicity. We
show that both ${\Fhdm}$ and
${\Fosgm}$ have a period-two orbit.

\begin{lem}[Orbit of the dynamical system] \label{lem:orbit}
  Suppose $A$ has distinct eigenvalues between $[\mu, L]$: $L = \lambda_1
  > \lambda_2 > \cdots > \lambda_{n - 1} > \lambda_n = \mu$. Let $q_i$ denote
  the unit eigenvectors of $A$ associated with eigenvalues $\lambda_i$. Then
  the dynamical systems $\mathcal{F}_{{\hdm}}$ and
  $\mathcal{F}_{\os}$ have a common period-$2$ orbit:
  \[ \{ (\alpha^{\star}, \hat{z}^{\star}_1), (\alpha^{\star},
     \hat{z}^{\star}_2) \} \assign \big\{ ( \tfrac{2}{L + \mu}, \gamma
     e_1 + \delta e_n ), ( \tfrac{2}{L + \mu}, - \gamma e_1 +
     \delta e_n ) \big\}, \]
  where $e_1 = (1, 0 \ldots, 0)$, $e_n = (0, \ldots, 0, 1)$, $\gamma =
  \tfrac{\mu}{\sqrt{L^2 + \mu^2}}$, and $\delta = \pm \tfrac{L}{\sqrt{L^2 +
  \mu^2}}$.
\end{lem}

\input{mpproofs/proof-lem-orbit.tex}

\begin{rem}
Both {\hdmclassic} and {\oh} have additional period-two orbits when the iterates fall in a subspace orthogonal to certain eigenvectors of $A$. For example, when $\langle q_1, z^1 \rangle = 0$, then the orbit can be derived based on $L = \lambda_2$ (i.e. the maximum eigenvalue restricted on subspace). Without loss of generality, we focus on the orbit generated by the iterations that satisfy $|\langle q_1, z^k \rangle| > 0, |\langle q_n, z^k \rangle| > 0$ for all $k$. In this case, the algorithm iterations always have non-vanishing components along both top and bottom eigen directions.
\end{rem}

Although both {\hdmclassic} and {\oh} share the same orbits, \Cref{thm:orbit-stability} shows that they exhibit different stability properties around the orbit. The proof of \Cref{thm:orbit-stability} is deferred to appendix due to its length and technicality.

\begin{thm} \label{thm:orbit-stability}
Let $\mathcal{J}_{\hdm}$ and $\mathcal{J}_{\oh}$ denote the Jacobians of $\Fcal_{\hdm}$ and $\Fcal_\os$ respectively. For any $\eta \in (0, \frac{1}{L}]$, their spectral radius satisfy
\[ \rho(\mathcal{J}_{\hdm}(\alpha^\star, \hat{z}_1^\star)\mathcal{J}_{\hdm}(\alpha^\star, \hat{z}_2^\star)) < 1 \qquad \text{and} \qquad \rho(\mathcal{J}_{\os}(\alpha^\star, \hat{z}_1^\star)\mathcal{J}_{\os}(\alpha^\star, \hat{z}_2^\star)) > 1\]

In other words, the orbit $\{ (\alpha^{\star}, \hat{z}^{\star}_1),
  (\alpha^{\star}, \hat{z}^{\star}_2) \}$ is stable for
  ${\Fhdm}$ but unstable for
  ${\Fosgm}$.
\end{thm}

\Cref{thm:orbit-stability} indicates that $\{\alpha_k\}$ in {\hdmclassic} will get trapped around $\alpha^\star$ and the iterates $\{ x^k \}$ approach the optimal solution $x^{\star} = 0$ from two specific directions $d_1 \assign
\gamma q_1 + \delta q_n$ and $d_2 \assign - \gamma q_1 + \delta q_n$ (\Cref{fig:dynamics} (right)). This spectral property has been extensively discussed in the literature \cite{huang2020gradient,dai2006new}, and a common belief is that sticking to the orbit results in slow convergence \cite{zou2022delayed}. In particular, we have
\[ \tfrac{\mathd}{\mathd \alpha}h_{x^\star + \lambda d_1} (\alpha^\star) = \tfrac{\mathd}{\mathd \alpha} h_{x^\star + \lambda d_2} (\alpha^\star) = 0\]
for all $\lambda \neq 0$, which means $\alpha^\star$ is exactly the optimal stepsize for this trajectory. As a result, {\hdmclassic} will stop adapting $\{\alpha_k\}$: {\hdmclassic} drives the iterates to a trajectory with no regret. 
However, despite {\hdmclassic} learns the optimal stepsize on this trajectory, this trajectory is so bad that even the optimal stepsize still offers only slow convergence.
In this respect, instability of {\oh} keeps $\{\alpha_k\}$ oscillating between large and small stepsizes and balances different spectral components of $x^k$. The idea of alternating between large and small stepsizes to leave the orbit has been explicitly used in the literature \cite{de2013spectral,zou2021fast} by designing specific stepsize schedulers. In contrast, {\oh} automatically achieves this goal.

%% file: proofs/proof-lem-orbit.tex
To show the orbit, we simply plug in $(\alpha^{\star}, \hat{z}^{\star}_1)$ and
$(\alpha^{\star}, \hat{z}^{\star}_2)$ into the updates of both {\hdm}
and {\os}. Observe that $L^2 \gamma^2 = \mu^2 \delta^2$ and the
eigenvectors $q_1$ and $q_n$ are orthogonal to each other. Then
\begin{align}
  \langle \hat{z}^{\star}_1, \Lambda^3 \hat{z}^{\star}_1 \rangle & = \langle
  \gamma e_1 + \delta e_n, L^3 \gamma e_1 + \mu^3 \delta e_n \rangle = L^3
  \gamma^2 + \mu^3 \delta^2 = (L + \mu) L^2 \gamma^2, \label{eqn:pf-orbit-1}\\
  \langle \hat{z}^{\star}_2, \Lambda^3 \hat{z}^{\star}_2 \rangle & = \langle -
  \gamma e_1 + \delta e_n, - L^3 \gamma e_1 + \mu^3 \delta e_n \rangle = L^3
  \gamma^2 + \mu^3 \delta^2 = (L + \mu) L^2 \gamma^2, \label{eqn:pf-orbit-2}\\
  \langle \hat{z}^{\star}_1, \Lambda^2 \hat{z}^{\star}_1 \rangle & = \langle
  \gamma e_1 + \delta e_n, L^2 \gamma e_1 + \mu^2 \delta e_n \rangle = L^2
  \gamma^2 + \mu^2 \delta^2 = 2 L^2 \gamma^2, \label{eqn:pf-orbit-3}\\
  \langle \hat{z}^{\star}_2, \Lambda^2 \hat{z}^{\star}_2 \rangle & = \langle -
  \gamma e_1 + \delta e_n, - L^2 \gamma e_1 + \mu^2 \delta e_n \rangle = L^2
  \gamma^2 + \mu^2 \delta^2 = 2 L^2 \gamma^2, \label{eqn:pf-orbit-4}
\end{align}
and hence
\begin{equation} \label{eqn:pf-orbit-5}
\tfrac{\langle \hat{z}^{\star}_1, \Lambda^3 \hat{z}^{\star}_1
  \rangle}{\langle \hat{z}^{\star}_1, \Lambda^2 \hat{z}^{\star}_1 \rangle} =
  \tfrac{\langle \hat{z}^{\star}_2, \Lambda^3 \hat{z}^{\star}_2
  \rangle}{\langle \hat{z}^{\star}_2, \Lambda^2 \hat{z}^{\star}_2 \rangle} =
  \tfrac{L + \mu}{2} . 
\end{equation}
The $\alpha$-update of both systems ${\Fhdm}$ and
${\Fosgm}$ are the same. Starting from either state
$(\alpha^{\star}, \hat{z}^{\star}_1)$ and $(\alpha^{\star},
\hat{z}^{\star}_2)$, the new state of stepsize is
\begin{equation}\label{eqn:pf-orbit-6}
[\alpha^{\star}]^+ (\hat{z}^{\star}_i) = \alpha^{\star} + \eta -
  \alpha^{\star} \eta \tfrac{\langle \hat{z}^{\star}_i, \Lambda^3
  \hat{z}^{\star}_i \rangle}{\langle \hat{z}^{\star}_i, \Lambda^2
  \hat{z}^{\star}_i \rangle} = \alpha^{\star} + \eta - \alpha^{\star} \eta
  \tfrac{L + \mu}{2} = \alpha^{\star} .
\end{equation}
Since the new state of stepsize remains unchanged to be $\alpha^{\star} =
\tfrac{L + \mu}{2}$, both systems end up with the same new state:
\begin{align}
  {\Fhdm} (\alpha^{\star}, \hat{z}^{\star}_i) & =
  ( [\alpha^{\star}]^+ (\hat{z}^{\star}_i), \tfrac{(I -
  [\alpha^{\star}]^+ (\hat{z}^{\star}_i) \Lambda) \hat{z}^{\star}_i}{\| (I -
  [\alpha^{\star}]^+ (\hat{z}^{\star}_i) \Lambda) \hat{z}^{\star}_i \|}
  ) = ( \alpha^{\star}, \tfrac{(I - \alpha^{\star} \Lambda)
  \hat{z}^{\star}_i}{\| (I - \alpha^{\star} \Lambda) \hat{z}^{\star}_i \|}
  ) ; \nonumber\\
  {\Fosgm} (\alpha^{\star}, \hat{z}^{\star}_i) & =
  ( [\alpha^{\star}]^+ (\hat{z}^{\star}_i), \tfrac{(I - \alpha^{\star}
  \Lambda) \hat{z}^{\star}_i}{\| (I - \alpha^{\star} \Lambda)
  \hat{z}^{\star}_i \|} ) = ( \alpha^{\star}, \tfrac{(I -
  \alpha^{\star} \Lambda) \hat{z}^{\star}_i}{\| (I - \alpha^{\star} \Lambda)
  \hat{z}^{\star}_i \|} ) . \nonumber
\end{align}
The updates of direction before the normalization are
\begin{align}
  (I - \alpha^{\star} \Lambda) \hat{z}^{\star}_1 & = ( 1 - \tfrac{2 L}{L
  + \mu} ) \gamma e_1 + ( 1 - \tfrac{2 \mu}{L + \mu} ) \delta
  e_n = \tfrac{L - \mu}{L + \mu} (- \gamma e_1 + \delta e_n) = \tfrac{\kappa -
  1}{\kappa + 1} \hat{z}^{\star}_2 ; \label{eqn:pf-orbit-7}\\
  (I - \alpha^{\star} \Lambda) \hat{z}^{\star}_2 & = - ( 1 - \tfrac{2
  L}{L + \mu} ) \gamma e_1 + ( 1 - \tfrac{2 \mu}{L + \mu} )
  \delta e_n = \tfrac{L - \mu}{L + \mu} (\gamma e_1 + \delta e_n) =
  \tfrac{\kappa - 1}{\kappa + 1} \hat{z}^{\star}_1 . \label{eqn:pf-orbit-8}
\end{align}
Since $\hat{x}^{\star}_1$ and $\hat{x}^{\star}_2$ both have unit norm, after
the normalization, we have
\[ {\Fhdm} (\alpha^{\star}, \hat{z}^{\star}_1)
   ={\Fosgm} (\alpha^{\star}, \hat{z}^{\star}_1) =
   (\alpha^{\star}, \hat{z}^{\star}_2) \quad \text{and} \quad
   {\Fhdm} (\alpha^{\star}, \hat{z}^{\star}_2)
   ={\Fosgm} (\alpha^{\star}, \hat{z}^{\star}_2) =
   (\alpha^{\star}, \hat{z}^{\star}_1) . \]
This proves that $\{ (\alpha^{\star}, \hat{z}^{\star}_1),
(\alpha^{\star}, \hat{z}^{\star}_2) \}$ is an orbit with period $2$.

%% file: sec_practical_hbfeedback.tex
\subsection{Heavy-ball momentum and potential function-based feedback}
{\osgm} framework begins with the gradient-based update $x^{k + 1} = x^k - P_k
\nabla f (x^k)$. In contrast, modern first-order methods with good empirical
performance {\cite{kingma2014adam}} typically add a momentum component
\begin{equation} \label{eqn:gd-hb}
   x^{k + 1} = x^k - P_k \nabla f (x^k) + \beta_k (x^{k} - x^{k-1}),
\end{equation}
where $x^{k} - x^{k-1}$ is known as heavy-ball momentum
{\cite{danilova2020non,polyak1964some,polyak1987introduction}} and $\beta_k
\in \Bcal \subseteq \Rbb$ is the extrapolation coefficient. 
To adapt $(P, \beta)$ jointly, a straightforward approach defines feedback based on the objective
value {\cite{gunes2018online}} as in {\oh}:
\[ \tfrac{f (x - P \nabla f (x) + \beta (x - x^-)) - f (x)}{\| \nabla f (x)
   \|^2} . \]
However, this direct extension does not yield global convergence for two
reasons: 1) due to the non-monotone behavior of standard heavy-ball momentum
{\cite{danilova2020non}}, there is no universal $(P^{\star}, \beta^{\star})$
that ensures descent of the objective value; 2) a momentum method can take a large step even when the gradient norm is small, so it may happen that $\| x - x^- \| \gg \|
\nabla f (x) \|^2$ and hence $\tfrac{f (x - P \nabla f (x) + \beta (x - x^-))-f(x)}{\|
\nabla f (x) \|^2}$ can have an arbitrarily large Lipschitz constant $(P, \beta)$ depending on $x, x^{-}$. 
Feedback based on the function value alone does not faithfully reflect
the quality of $(P, \beta)$ and cannot guide the scheduler constructively. 
To identify a suitable feedback, recall that heavy-ball momentum admits a potential function \cite{danilova2020non} \begin{equation} \label{eqn:hb-potential}
 \varphi_{\omega} (x, x^{-}) \assign f (x) - f^{\star} +
   \tfrac{\omega}{2} \| x - x^{-} \|^2.
\end{equation}
This potential function monotonically decreases if $(P, \beta)$ are appropriately chosen. 
\begin{lem}[Heavy-ball potential {\cite[Theorem 2]{danilova2020non}}]\label{lem:heavyball-potential}
Suppose $f$ is $L$-smooth and convex. 
Given any $x, x^- \nin \mathcal{X}^{\star}$, denote $x^+ \assign x - P \nabla f(x) + \beta (x - x^-)$ with 
$P = \alpha I$, $\alpha \leq \tfrac{1}{L}$, and $\beta \in [ 0, \sqrt{1 - \alpha L} ]$.
Then the potential $\varphi_{\omega}$ with $\omega = \tfrac{1 - \alpha L}{2 \alpha}$ satisfies
\[ \varphi_{\omega} (x^{+}, x) \leq \varphi_{\omega} (x, x^{-})
   - \tfrac{\alpha}{2} \big[ \| \nabla f (x) \|^2 + \tfrac{1 - \alpha L -
   \beta^2}{\alpha^2} \| x - x^{-} \|^2 \big]. \]
\end{lem}

\input{mpproofs/proof-hb-potential.tex}

\Cref{lem:heavyball-potential} resembles the descent lemma: the progress of
each heavy-ball step in reducing $\varphi_{\omega}$ is proportional to $\|
\nabla f (x) \|^2 + \tfrac{\tau}{2} \| x - x^- \|^2$ for $\tau = \tfrac{1 - \alpha L -
   \beta^2}{\alpha^2}> 0$.
Therefore, we define the heavy-ball feedback as
\begin{equation} \label{eqn:hb-feedback}
h_{x, x^-} (P, \beta) \assign \tfrac{\varphi_{\omega} (x^+ (P, \beta), x)
   - \varphi_{\omega} (x, x^-)}{\| \nabla f (x) \|^2 + \frac{\tau}{2} \| x -
   x^- \|^2},
\end{equation}
where $x^+ (P, \beta) = x - P \nabla f (x) + \beta (x - x^-)$ and $\tau >0$ is a parameter. \\

The heavy-ball feedback \eqref{eqn:hb-feedback} depends on two points $x$ and $x^-$, which we denote by a state $z = (z_1, z_2) = (x, x^-)$. 
With slight abuse of notation, we denote $f(z):=f(z_1)$ and $\nabla f(z):= \nabla f(z_1)$. We will also slightly abuse the notation of the gradient by directly concatenating the gradients with respect to $P$ and $\beta$: $\nabla h_z(P, \beta) = ( \nabla_P h_z(P, \beta), \nabla_\beta h_z(P, \beta))$. The pushforward operation and the heavy-ball feedback are written with respect to state $z$ as
\begin{equation*}
   z^+ (P, \beta) \assign (z_1 - P \nabla f (z) + \beta (z_1 - z_2),~z_1),\qquad 
   h_{z} (P, \beta) \assign \tfrac{\varphi_{\omega} (z^+ (P, \beta)) - \varphi_{\omega} (z)}{\| \nabla f (z) \|^2 + \frac{\tau}{2} \| z_1 - z_2 \|^2}.
\end{equation*}
Using the superscript $k$ for the iteration count, the first component of the state $z^k \assign (z_1^k, z_2^k)$ satisfies $z^{k}_1 = x^k$;
however, the second component $z^k_2$ need not be equal to $x^{k-1}$.
We will utilize this flexibility to implement the monotone landscape action outlined below. In the product space $\Pcal \times \Bcal$ of stepsize and momentum, define the norm and its associated dual norm as
\begin{equation} \label{eqn:hb-norm}
   \|(P, \beta)\| \assign \sqrt{\|P\|_F^2 + \tfrac{1}{L^2} \beta^2}, \qquad 
   \|(P, \beta)\|_{\ast} \assign \sqrt{\| P \|_F^2 + L^2 \beta^2}.
\end{equation}
The weight $\tfrac{1}{L^2}$ ensures both $\|P\|_F^2$ and $\tfrac{1}{L^2} \beta^2$ have the same physical dimension since $P$ lies in the space of inverse Hessian and $\beta$ is dimensionless.
We prove properties of feedback $h_{z}(P, \beta)$ with respect to this norm.

\begin{lem}[Properties of heavy-ball feedback] \label{lem:hb-properties}
Suppose $f$ is convex and $L$-smooth and let $\tau \geq 2 L^2$. Then for any state $z = (z_1, z_2)$ with $z_1, z_2 \nin \Xcal^\star$, the heavy-ball feedback $h_{z} (P, \beta)$ is jointly convex in $(P, \beta)$ and $(L+\omega)$-smooth in $(P, \beta)$.
         Moreover, the gradient of $h_{z}(P, \beta)$ takes the form
   \[ \nabla h_{z} (P, \beta) = \Big(
   - \tfrac{[\nabla f (z^+(P, \beta)) + \omega (-P \nabla f(z) + \beta (z_1 - z_2))]
   \nabla f (z)^{\top}}{\| \nabla f (z) \|^2 + \frac{\tau}{2} \| z_1 - z_2
   \|^2}, ~
   \tfrac{\langle [\nabla f (z^+(P, \beta)) + \omega (-P \nabla f(z) + \beta (z_1 - z_2))],~z_1 - z_2 \rangle}{\| \nabla f (z) \|^2 + \frac{\tau}{2} \| z_1 - z_2 \|^2} \Big) . \]
\end{lem}

\input{mpproofs/proof-hb-properties.tex}

At each iteration, define the progress with respect to the heavy-ball feedback:
\[ b_k \assign \tfrac{\varphi_{\omega} (z^{k + 1}) - \varphi_{\omega}(z^k)}{\| \nabla f (z^k) \|^2 + \frac{\tau}{2} \| z_1^k - z_2^k \|^2} . \]
\Cref{thm:heavyball-reduction} shows a reduction from the function value gap to the cumulative progress $\sum^K_{k = 1} b_k$.

\begin{thm}[Heavy-ball reduction] \label{thm:heavyball-reduction}
Choose the parameters $\omega > 0$ and $\tau > 0$ such that $\tau \geq 2 L \omega$. 
Let $\{ z^k \}$ be any sequence of states such that $\varphi_{\omega} (z^{k+1}) \leq \varphi_{\omega} (z^k)$ and let the initial state $z^1 = (z_1^1, z_2^1)$ satisfy $z_1^1 = z_2^1$.
\begin{itemize}[leftmargin=15pt]
   \item If $f$ is convex, then
   \begin{equation} \label{eqn:hb-reduction-cvx}
      f (z^{K + 1}) - f^\star \leq \tfrac{f (z^1) - f^\star}{1 +
      V \sum^K_{k = 1} - b_k},
   \end{equation}
   where $V \assign \min \big\{ \tfrac{f (z^1) - f^\star}{4 \Delta^2}, \tfrac{L}{2} \big\}$ and $\Delta \assign \max_{x \in \{ x : f (x) \leq f (x^1) \}} \min_{x^{\star} \in \mathcal{X}^{\star}}  \| x - x^{\star} \|$.
   \item If $f$ is $\mu$-strongly convex, then
   \begin{equation} \label{eqn:hb-reduction-strongly-cvx}
   \textstyle f (z^{K + 1}) - f^\star \leq [f (z^1) - f^\star] ( 1 -
      \tfrac{\mu}{K} \sum_{k = 1}^K - b_k  )^K .
   \end{equation}
\end{itemize}
\end{thm}

\input{mpproofs/proof-hb-reduction.tex}

The heavy-ball reduction (\Cref{thm:heavyball-reduction}) requires monotonicity of the potential $\varphi_{\omega}(z^k)$, instead of function value $f(z^k)$. Hence, given the proposal state from the scheduler
\begin{equation} \label{eqn:hb-zproposal}
z^{k+1/2} \assign (z^k)^+ = ( z_1^k - P_k \nabla f (z^k) + \beta_k (z_1^k - z_2^k), z_1^k ),
\end{equation}
the monotone landscape action for heavy-ball feedback must ensure the potential is nonincreasing: 
\[
\varphi_{\omega}(z^{k+1}) \leq \min \{ \varphi_{\omega}(z^{k+1/2}), \varphi_{\omega}(z^k) \},
\]
where $z^{k+1}$ is the state returned by the monotone landscape action.
We may implement this monotone landscape action by the null step $z^{k+1} = \argmin \big\{ \varphi_{\omega} ( z^{k + 1/2}), \varphi_{\omega} (z^k) \big\}$, 
which simply returns whichever state has the smaller potential value. In practice, the lookahead action described in \Cref{sec:ob} is often preferred.

%% file: proofs/proof-hb-potential.tex
The proof is adapted from \cite[Theorem 2]{danilova2020non} for completeness. The relation $x^+ = x - \alpha \nabla f (x) + \beta (x - x^-)$ implies
\begin{align}
x^+ - x & = - \alpha \nabla f (x) + \beta (x - x^-), \label{eqn:pf-hb-potential-1}\\
\| x^+ - x \|^2 & = \alpha^2 \| \nabla f (x) \|^2 + \beta^2 \| x - x^-
\|^2 - 2 \alpha \beta \langle \nabla f (x), x - x^- \rangle . \label{eqn:pf-hb-potential-2}
\end{align}
The $L$-smoothness of $f$ implies
\begin{equation} \label{eqn:pf-hb-potential-3}
    f (x^+) \leq f (x) + \langle \nabla f (x), x^+ - x \rangle + \tfrac{L}{2}
    \| x^+ - x \|^2 .
\end{equation}
Plugging \eqref{eqn:pf-hb-potential-1} and \eqref{eqn:pf-hb-potential-2} into \eqref{eqn:pf-hb-potential-3} gives
\begin{align}
f (x^+) & \leq f (x) - \alpha \| \nabla f (x) \|^2 + \beta \langle \nabla
f (x), x - x^- \rangle \nonumber\\
& \quad + \tfrac{\alpha^2 L}{2} \| \nabla f (x) \|^2 + \tfrac{\beta^2
L}{2} \| x - x^- \|^2 - L \alpha \beta \langle \nabla f (x), x - x^-
\rangle \nonumber\\
& = f (x) + ( \tfrac{\alpha^2 L}{2} - \alpha ) \| \nabla f (x)
\|^2 + \beta (1 - L \alpha) \langle \nabla f (x), x - x^- \rangle +
\tfrac{\beta^2 L}{2} \| x - x^- \|^2 . \label{eqn:pf-hb-potential-4}
\end{align}
Multiplying \eqref{eqn:pf-hb-potential-2} by $\tfrac{1 - \alpha L}{2 \alpha}$ and adding to \eqref{eqn:pf-hb-potential-4}, we have
\begin{align}
f (x^+) + \tfrac{1 - \alpha L}{2 \alpha} \| x^+ - x \|^2 & \leq f (x) +
( \tfrac{\alpha^2 L}{2} - \alpha ) \| \nabla f (x) \|^2 + \beta
(1 - L \alpha) \langle \nabla f (x), x - x^- \rangle + \tfrac{\beta^2
L}{2} \| x - x^- \|^2 \nonumber\\
& \quad + \tfrac{\alpha (1 - \alpha L)}{2} \| \nabla f (x) \|^2 +
\tfrac{\beta^2 (1 - \alpha L)}{2 \alpha} \| x - x^- \|^2 - \beta (1 - L
\alpha) \langle \nabla f (x), x - x^- \rangle \nonumber\\
& = f (x) + \tfrac{\beta^2}{2 \alpha} \| x - x^- \|^2 - \tfrac{\alpha}{2}
\| \nabla f (x) \|^2 \nonumber\\
& = f (x) + \tfrac{1 - \alpha L}{2 \alpha} \| x - x^- \|^2 -
\tfrac{\alpha}{2} [ \| \nabla f (x) \|^2 + \tfrac{1 - \alpha L -
\beta^2}{\alpha^2} \| x - x^- \|^2 ]. \label{eqn:pf-hb-potential-5}
\end{align}
Let $\omega = \tfrac{1 - \alpha L}{2 \alpha}$. Then \eqref{eqn:pf-hb-potential-5} becomes
\begin{equation*}
    \varphi_{\omega} (x^{+}, x) \leq \varphi_{\omega} (x, x^{-})
   - \tfrac{\alpha}{2} \big[ \| \nabla f (x) \|^2 + \tfrac{1 - \alpha L -
   \beta^2}{\alpha^2} \| x - x^{-} \|^2 \big]
\end{equation*}
and this completes the proof.

%% file: proofs/proof-hb-properties.tex
\paragraph{Notation.} Recall that $z = (x, x^-)$. We explicitly use $(x, x^-)$ when establishing the analytic properties of $\varphi_\omega$.

\paragraph{Convexity.} To prove the convexity of $h_{x, x^-} (P, \beta)$ in $(P, \beta)$, it suffices to show the joint convexity of
\[ \varphi_{\omega} (x^+ (P, \beta), x) \assign f (x^+ (P, \beta)) +
   \tfrac{\omega}{2} \| x^+ (P, \beta) - x \|^2. \]
Since $x^+ (P, \beta) \assign x - P \nabla f (x) + \beta (x - x^-)$ is affine
in $(P, \beta)$ and $f$ and $\| \cdot \|^2$ are both convex, the composition
rule of convex functions implies the joint convexity of $\varphi_{\omega} (x^+ (P, \beta), x)$.

\paragraph{Gradient.}The chain rule implies
\begin{align}
  \nabla_P h_{x, x^-} (P, \beta) & = \tfrac{- [\nabla f (x^+ (P, \beta)) +
  \omega (x^+ (P, \beta) - x)] \nabla f (x)^{\top}}{\| \nabla f (x) \|^2 +
  \frac{\tau}{2} \| x - x^- \|^2}, \label{eqn:pf-property-ghb-P}\\
  \nabla_{\beta} h_{x, x^-} (P, \beta) & = \tfrac{\langle \nabla f (x^+ (P,
  \beta)) + \omega (x^+ (P, \beta) - x), ~x - x^- \rangle}{\| \nabla f (x)
  \|^2 + \frac{\tau}{2} \| x - x^- \|^2} . \label{eqn:pf-property-ghb-beta}
\end{align}

\paragraph{Smoothness.}The difference of the gradients in $P$ and $\beta$ are
respectively bounded by
\begin{align}
  \| \nabla_P h_{x, x^-} (P_1, \beta_1) - \nabla_P h_{x, x^-} (P_2,
  \beta_2) \|_F
  ={}& \tfrac{\| [\nabla f (x^+ (P_1, \beta_1)) - \nabla f (x^+ (P_2, \beta_2))
  + \omega (x^+ (P_1, \beta_1) - x^+ (P_2, \beta_2))] \nabla f (x)^{\top}
  \|_F}{\| \nabla f (x) \|^2 + \frac{\tau}{2} \| x - x^- \|^2} \nonumber\\
  ={}& \tfrac{\| \nabla f (x^+ (P_1, \beta_1)) - \nabla f (x^+ (P_2, \beta_2))
  + \omega (x^+ (P_1, \beta_1) - x^+ (P_2, \beta_2)) \| \| \nabla f (x) \|}{\|
  \nabla f (x) \|^2 + \frac{\tau}{2} \| x - x^- \|^2} \label{eqn:pf-property-1}\\
| \nabla_{\beta} h_{x, x^-} (P_1, \beta_1) - \nabla_{\beta} h_{x, x^-}
  (P_2, \beta_2) |
  ={}& \tfrac{| \langle \nabla f (x^+ (P_1, \beta_1)) - \nabla f (x^+ (P_2,
  \beta_2)) + \omega (x^+ (P_1, \beta_1) - x^+ (P_2, \beta_2)), x - x^-
  \rangle |}{\| \nabla f (x) \|^2 + \frac{\tau}{2} \| x - x^- \|^2}
  \nonumber\\
  \leq{}& \tfrac{\| \nabla f (x^+ (P_1, \beta_1)) - \nabla f (x^+ (P_2,
  \beta_2)) + \omega (x^+ (P_1, \beta_1) - x^+ (P_2, \beta_2)) \| \| x - x^-
  \|}{\| \nabla f (x) \|^2 + \frac{\tau}{2} \| x - x^- \|^2}, \label{eqn:pf-property-2}
\end{align}
where \eqref{eqn:pf-property-2} uses Cauchy's inequality. Observe that the common norm in the numerator can be bounded by
\begin{align}
  & \| \nabla f (x^+ (P_1, \beta_1)) - \nabla f (x^+ (P_2, \beta_2)) + \omega
  (x^+ (P_1, \beta_1) - x^+ (P_2, \beta_2)) \| \nonumber\\
  \leq{}& \| \nabla f (x^+ (P_1, \beta_1)) - \nabla f (x^+ (P_2, \beta_2)) \| +
  \omega \| x^+ (P_1, \beta_1) - x^+ (P_2, \beta_2) \| \tag{by $\|a + b\| \leq \|a\| + \|b\|$}\\
  \leq{}& (L + \omega) \| x^+ (P_1, \beta_1) - x^+ (P_2, \beta_2) \|
  \tag{by $L$-smooth of $f$}\\
  ={}& (L + \omega) \| (P_2 - P_1) \nabla f (x) + (\beta_1 - \beta_2) (x - x^-)
  \| \nonumber\\
  \leq{}& (L + \omega) \| P_1 - P_2 \|  \| \nabla f (x) \| + | \beta_1 -
  \beta_2 |  \| x - x^- \| \tag{by submultiplicativity}\\
  \leq{}& (L + \omega) \sqrt{\| P_1 - P_2 \|^2_F + \tfrac{1}{L^2} (\beta_1 -
  \beta_2)^2}  \sqrt{\| \nabla f (x) \|^2 + L^2 \| x - x^- \|^2} \tag{by Cauchy's inequality}\\
  ={}& (L + \omega) \sqrt{\| \nabla f (x) \|^2 + L^2 \| x - x^- \|^2} \| (P_1,
  \beta_1) - (P_2, \beta_2) \| . \tag{by defintion of $\|(P, \beta)\|$ norm}
\end{align}
Then \eqref{eqn:pf-property-1} and \eqref{eqn:pf-property-2} can be further bounded by
\begin{align}
  \| \nabla_P h_{x, x^-} (P_1, \beta_1) - \nabla_P h_{x, x^-} (P_2, \beta_2)
  \|_F & \leq (L + \omega)   \tfrac{\sqrt{\| \nabla f (x) \|^2 + L^2 \| x - x^- \|^2}  \| \nabla f (x)
  \|}{\| \nabla f (x) \|^2 + \frac{\tau}{2} \| x - x^- \|^2} \| (P_1,
  \beta_1) - (P_2, \beta_2) \| , \nonumber\\
  | \nabla_{\beta} h_{x, x^-} (P_1, \beta_1) - \nabla_{\beta} h_{x, x^-} (P_2,
  \beta_2) | & \leq (L + \omega)
  \tfrac{\sqrt{\| \nabla f (x) \|^2 + L^2 \| x - x^- \|^2}  \| x - x^- \|}{\|
  \nabla f (x) \|^2 + \frac{\tau}{2} \| x - x^- \|^2} \| (P_1,
  \beta_1) - (P_2, \beta_2) \| . \nonumber
\end{align}
Hence the dual norm $\| \cdot \|_{\ast}$ of the gradient difference is bounded by
\begin{align}
  & \| (\nabla_P h_{x, x^-} (P_1, \beta_1) - \nabla_P h_{x, x^-} (P_2,
  \beta_2), \nabla_{\beta} h_{x, x^-} (P_1, \beta_1) - \nabla_{\beta} h_{x,
  x^-} (P_2, \beta_2)) \|_{\ast}^2 \nonumber\\
  ={}& \| \nabla_P h_{x, x^-} (P_1, \beta_1) - \nabla_P h_{x, x^-} (P_2,
  \beta_2) \|^2_F + L^2 | \nabla_{\beta} h_{x, x^-} (P_1, \beta_1) -
  \nabla_{\beta} h_{x, x^-} (P_2, \beta_2) |^2 \nonumber\\
  \leq{}& (L + \omega)^2  \tfrac{[\|
  \nabla f (x) \|^2 + L^2 \| x - x^- \|^2]^2}{[ \| \nabla f (x) \|^2 +
  \frac{\tau}{2} \| x - x^- \|^2 ]^2} \| (P_1,
  \beta_1) - (P_2, \beta_2) \|^2 \nonumber\\
  \leq{}& (L + \omega)^2 \max \{
  \tfrac{4 L^4}{\tau^2}, 1 \} \| (P_1,
  \beta_1) - (P_2, \beta_2) \|^2 \tag{by $\tfrac{a^2 + \gamma b^2}{a^2 + \delta b^2} \leq \max\{ \tfrac{\gamma}{\delta}, 1 \}$}\\
  ={}& (L + \omega)^2 \| (P_1,
  \beta_1) - (P_2, \beta_2) \|^2, \tag{by $\tau \geq 2L^2$}
\end{align}
which implies $(L + \omega)$-smoothness of $h_{x, x^-} (P, \beta)$.  Finally, as a by-product relation that will be later used, we derive the dual norm $\| \cdot \|_{\ast}$ of the gradient
is bounded by
\begin{align}
  \| (\nabla_P h_{x, x^-} (P, \beta), \nabla_{\beta} h_{x, x^-} (P, \beta))
  \|^2_{\ast} & = \| \nabla_P h_{x, x^-} (P, \beta) \|^2 + L^2 \|
  \nabla_{\beta} h_{x, x^-} (P, \beta) \|^2 \nonumber\\
  & \leq \tfrac{\| \nabla f (x^+ (P, \beta)) + \omega (x^+ (P, \beta) - x)
  \|^2 [\| \nabla f (x) \|^2 + L^2 \| x - x^- \|^2]}{[ \| \nabla f (x)
  \|^2 + \frac{\tau}{2} \| x - x^- \|^2 ]^2} \tag{by \eqref{eqn:pf-property-ghb-P} and \eqref{eqn:pf-property-ghb-beta}}\\
  & \leq \max \{ \tfrac{2 L^2}{\tau}, 1 \} \tfrac{\| \nabla f (x^+
  (P, \beta)) + \omega (x^+ (P, \beta) - x) \|^2}{\| \nabla f (x) \|^2 +
  \frac{\tau}{2} \| x - x^- \|^2} . \label{eqn:pf-property-3}
\end{align}

%% file: proofs/proof-hb-reduction.tex
Note that the initial state $z^1$ with $z_2^1 = z_1^1$ implies that the initial potential is the initial suboptimality:
\begin{equation*}
    \varphi_{\omega} (z^1) = f (z^1) - f^{\star} + \tfrac{\omega}{2} \| z_1^1 - z_2^1 \|^2 = f (z^1) - f^{\star} .
\end{equation*}
\paragraph{Convex $f$.}Observe that
\begin{equation} \label{eqn:pf-reduction-1}
    \varphi_{\omega} (z^{K + 1}) 
    = \frac{1}{\sum^K_{k = 1} \big[ \frac{\varphi_{\omega}(z^k) - \varphi_{\omega}(z^{k + 1})}{\varphi_{\omega}(z^{k + 1}) \varphi_{\omega}(z^k)} \big] + \tfrac{1}{\varphi_{\omega} (z^1)}}
    = \frac{1}{\sum^K_{k = 1} - b_k \big[ \frac{\| \nabla f (z^k) \|^2 + \frac{\tau}{2} \| z_1^k - z_2^k \|^2}{\varphi_{\omega}(z^{k + 1}) \varphi_{\omega}(z^k)} \big] + \tfrac{1}{f (z^1) - f^{\star}}} .
\end{equation}
Note that $b_k \leq 0$ since the algorithm guarantees that the potential is monotone $\varphi_{\omega} (z^{k
+ 1}) \leq \varphi_{\omega} (z^k)$. To further upper bound $\varphi_{\omega}
(z^{K + 1})$, we need to lower bound the ratio
\[ \tfrac{\| \nabla f (z^k) \|^2 + \frac{\tau}{2} \| z_1^k -
   z_2^k \|^2}{\varphi_{\omega}(z^{k + 1}) \varphi_{\omega}(z^k)} \geq
   \tfrac{\| \nabla f (z^k) \|^2 + \frac{\tau}{2} \| z_1^k -
   z_2^k \|^2}{\varphi_{\omega}(z^k)^2} = \tfrac{\| \nabla f (z^k)
   \|^2 + \frac{\tau}{2} \| z_1^k - z_2^k \|^2}{\big[ f (z^k) - f^{\star} + \frac{\omega}{2} \| z_1^k - z_2^k \|^2 \big]^2} . \]
\textbf{Case 1.} $f (z^k) - f^{\star} \geq \tfrac{\omega}{2} \| z_1^k
- z_2^k \|^2$. Observe that $z_1^k \in \{ x : f (x) \leq f (z^1)
\}$ for all $k \geq 1$ since the monotonicity of potential guarantees
\[ f (z^k) - f^{\star} \leq f (z^k) - f^{\star} + \tfrac{\omega}{2}
   \| z_1^k - z_2^k \|^2 = \varphi_{\omega} (z^k) \leq
   \varphi_{\omega} (z^1) = f (z^1) - f^{\star} . \]
Then $\tfrac{f (z^k) - f^{\star}}{\| \nabla f (z^k) \|} \leq \| z_1^k -
x^{\star} \| \leq \Delta \assign \max_{x \in \{ x : f (x) \leq f (z^1) \}}
\min_{x^{\star} \in \mathcal{X}^{\star}}  \| x - x^{\star} \|$ and hence
\begin{equation} \label{eqn:pf-reduction-2}
\tfrac{\| \nabla f (z^k) \|^2 + \frac{\tau}{2} \| z_1^k -
z_2^k \|^2}{[ f (z^k) - f^{\star} +
\frac{\omega}{2} \| z_1^k - z_2^k \|^2 ]^2} \geq
\tfrac{\| \nabla f (z^k) \|^2}{4 [f (z^k) - f^{\star}]^2} \geq
\tfrac{1}{4 \Delta^2} .
\end{equation}
\textbf{Case 2.} $f (z^k) - f^{\star} \leq \tfrac{\omega}{2} \| z_1^k
- z_2^k \|^2$. Observe that
\begin{equation} \label{eqn:pf-reduction-2.5}
    \tfrac{\omega}{2} \| z_1^k - z_2^k \|^2 \leq
   \varphi_{\omega} (z^k) \leq \varphi_{\omega} (z^1) = f (z^1) - f^{\star}.
\end{equation}
Using \eqref{eqn:pf-reduction-2.5} and the choice of $\tau \geq 2 L \omega$, we have
\begin{equation} \label{eqn:pf-reduction-3}
\tfrac{\| \nabla f (z^k) \|^2 + \frac{\tau}{2} \| z_1^k -
z_2^k \|^2}{[ f (z^k) - f^{\star} +
\frac{\omega}{2} \| z_1^k - z_2^k \|^2 ]^2} \geq
\tfrac{\frac{\tau}{2} \| z_1^k - z_2^k \|^2}{\omega^2
\| z_1^k - z_2^k \|^4} = \tfrac{\tau}{2 \omega^2}
\tfrac{1}{\| z_1^k - z_2^k \|^2} \geq \tfrac{L}{2
[f (z^1) - f^{\star}]} .    
\end{equation}
Using \eqref{eqn:pf-reduction-1} and the two lower bounds \eqref{eqn:pf-reduction-2} and \eqref{eqn:pf-reduction-3}, the suboptimality can be bounded by
\begin{equation*}
    f(z^{K+1}) - f^{\star}
    \leq \varphi_{\omega} (z^{K + 1})
    \leq \tfrac{f (z^1) - f^{\star}}{1 + V \sum^K_{k = 1} - b_k},
\end{equation*}
where $V \assign \min \big\{ \tfrac{f (z^1) - f^{\star}}{4 \Delta^2}, \tfrac{L}{2} \big\}$.

\paragraph{$\mu$-strongly convex $f$.}Observe that
\begin{align}
  \tfrac{\varphi_{\omega} (z^{K + 1})}{\varphi_{\omega} (z^1)} = \textstyle \prod_{k =
  1}^K \tfrac{\varphi_{\omega} (z^{k + 1})}{\varphi_{\omega} (z^k)}
  & \leq \big( \tfrac{1}{K} \textstyle \sum_{k = 1}^K \tfrac{\varphi_{\omega} (z^{k +
  1})}{\varphi_{\omega} (z^k)} \big)^K \nonumber\\
  & = ( 1 + \tfrac{1}{K} \textstyle \sum_{k = 1}^K \tfrac{\varphi_{\omega} (z^{k +
  1}) - \varphi_{\omega} (z^k)}{\varphi_{\omega} (z^k)} )^K \nonumber\\
  & = ( 1 - \tfrac{1}{K} \textstyle \sum_{k = 1}^K - b_k  \tfrac{\| \nabla f (z^k)
  \|^2 + \frac{\tau}{2} \| z_1^k - z_2^k \|^2}{f (z^k) - f^{\star} + \frac{\omega}{2} \| z_1^k - z_2^k \|^2} )^K . \label{eqn:pf-reduction-4}
\end{align}

To further upper bound \eqref{eqn:pf-reduction-4}, we lower bound the ratio $\tfrac{\| \nabla f (z^k)\|^2 + \frac{\tau}{2} \| z_1^k - z_2^k \|^2}{f (z^k) - f^{\star} + \frac{\omega}{2} \| z_1^k - z_2^k \|^2}$ by case analysis.

\textbf{Case 1.} $f (z^k) - f^{\star} \leq \tfrac{\omega}{2} \| z_1^k
- z_2^k \|^2$. Then
\begin{equation} \label{eqn:pf-reduction-5}
    \tfrac{\| \nabla f (z^k) \|^2 + \frac{\tau}{2} \| z_1^k -
   z_2^k \|^2}{f (z^k) - f^{\star} + \frac{\omega}{2}
   \| z_1^k - z_2^k \|^2} \geq \tfrac{\frac{\tau}{2}
   \| z_1^k - z_2^k \|^2}{\omega \| z_1^k -
   z_2^k \|^2} = \tfrac{\tau}{2 \omega} .
\end{equation}

\textbf{Case 2.} $f (z^k) - f^{\star} \geq \tfrac{\omega}{2} \| z_1^k
- z_2^k \|^2$. Then using $\tfrac{\| \nabla f (z^k) \|^2}{f
(z_1^k) - f^{\star}} \geq 2 \mu$ to obtain
\begin{equation} \label{eqn:pf-reduction-6}
    \tfrac{\| \nabla f (z^k) \|^2 + \frac{\tau}{2} \| z_1^k -
   z_2^k \|^2}{f (z^k) - f^{\star} + \frac{\omega}{2}
   \| z_1^k - z_2^k \|^2} \geq \tfrac{\| \nabla f (z^k)
   \|^2}{2 [f (z^k) - f^{\star}]} \geq \mu .
\end{equation}
Since $\omega$ and $\tau$ are chosen such that $\tfrac{\tau}{2 \omega} \geq
L \geq \mu$, two lower bounds \eqref{eqn:pf-reduction-5}--\eqref{eqn:pf-reduction-6} together imply
\[ \tfrac{\| \nabla f (z^k) \|^2 + \frac{\tau}{2} \| z_1^k -
   z_2^k \|^2}{f (z^k) - f^{\star} + \frac{\omega}{2}
   \| z_1^k - z_2^k \|^2} \geq \min \{ \tfrac{\tau}{2
   \omega}, \mu \} = \mu, \]
and hence we conclude from \eqref{eqn:pf-reduction-4} that
\begin{align}
    f(z^{K+1}) - f^{\star} \leq \varphi_{\omega} (z^{K + 1}) \leq \varphi_{\omega} (z^1) \big( 1 - \tfrac{\mu}{K} \textstyle \sum_{k = 1}^K - b_k  \big)^K = [f (z^1) - f^{\star}] \big( 1 - \tfrac{\mu}{K} \textstyle \sum_{k = 1}^K - b_k \big)^K . \nonumber
\end{align}

%% file: sec_osgm_best.tex
\subsection{Algorithm design of {\ob}} \label{sec:ob}

We put together the insights above and propose {\ob}, a variant of {\os} that demonstrates competitive performance in practice.  {\ob} utilizes the lookahead landscape action, tailored to heavy-ball feedback, to suppress regret and thereby accelerate convergence.\\

Part I \cite{gao2025gradient} introduced the lookahead landscape action, which takes an additional gradient step on the proposed iterate $x^{k + 1 / 2}$ to suppress the regret of online gradient descent on feedback. 
Since the heavy-ball feedback is defined with respect to the potential function, the additional gradient step should be
taken with respect to the potential function $\varphi_{\omega}(\cdot, ~ z_2^{k + 1 / 2})$. 
As a result, the lookahead landscape action for heavy-ball feedback is
\begin{align}
  \Mlook (z^{k + 1 / 2}) & \assign{} ( z_1^{k + 1 / 2} - \tfrac{1}{L +
  \omega} \nabla_{1} \varphi_{\omega} (z^{k + 1 / 2}), ~ z_2^{k + 1 / 2} ) \nonumber\\
  & = ( z_1^{k + 1 / 2} - \tfrac{1}{L + \omega} [\nabla f (z^{k + 1 / 2}) +
  \omega (z_1^{k + 1 / 2} - z_2^{k + 1 / 2})], ~ z_2^{k + 1 / 2} ) \nonumber\\
  & = ( ( 1 - \tfrac{\omega}{L + \omega} ) z_1^{k + 1 / 2} +
  \tfrac{\omega}{L + \omega} z_2^{k + 1 / 2} - \tfrac{1}{L + \omega} \nabla f (z^{k + 1 / 2}), ~ z_2^{k + 1 / 2} ), \nonumber
\end{align}
where $\nabla_{1} \varphi_{\omega}(z)$ denotes the gradient with respect to $z_1$ and the stepsize $\frac{1}{L + \omega}$ is the inverse of the smoothness constant $L + \omega$ of the potential function $\varphi_{\omega}(z)$ with respect to $z_1$. 
{\ob} determines the next state $z^{k+1}$ by composing the lookahead and monotone landscape actions:
\begin{equation} \label{eqn:mono-look-action}
   \varphi_\omega (z^{k + 1}) \leq  \min \{\varphi_\omega ( \mathcal{M}_{\text{look}}(z^{k + 1 / 2})), \varphi_\omega(z^k) \},
\end{equation}
\Cref{alg:ohblook} summarizes {\ob} with monotone landscape action implemented by the null step.
We refer to \Cref{alg:ohblook} as {\ob} since it is empirically superior to other {\os} variants. The global convergence of {\os} is established in \Cref{thm:ohblook}.

\vspace{10pt}
\begin{algorithm}[H]
  \caption{{\ob}\label{alg:ohblook}}
  {\textbf{input} Initial state $z^1 = (z^1_1, z^1_2)$ with $z_1^1 = z_2^1$, initial stepsize $P_1 \in \mathbb{R}^{n \times n}$ and momentum parameter $\beta_1 \in \mathbb{R}$, online gradient descent stepsizes $\eta_P, \eta_{\beta} > 0$}
  
  \For{$k = 1, 2, \dots$}{
 $z^{k + 1 / 2} = ( z_1^k - P_k \nabla f (z^k) + \beta_k (z_1^k - z_2^{k}),~ z_1^k )$ \Comment{State proposal from scheduler}\\
 $z_{\text{look}}^{k + 1} = (z_1^{k + 1 / 2} -\tfrac{1}{L + \omega} \nabla_{1} \varphi_{\omega} (z^{k + 1 / 2}),~ z_2^{k+1/2} )$
          \Comment{Lookahead action}\\
Choose $z^{k+1}$ s.t. $\varphi_\omega(z^{k+1}) \leq \min \{ \varphi_{\omega} ( z_{\text{look}}^{k+1}), \varphi_{\omega} (z^k) \}$ \Comment{Monotone action by null step}\\
 $P_{k + 1} = P_k - \eta_P \nabla_P h_{z^k} (P_k, \beta_k)$
         \Comment{Online gradient update on $P$}\\
$\beta_{k + 1} = \beta_k - \eta_{\beta} \nabla_{\beta} h_{z^k} (P_k, \beta_k)$
          \Comment{Online gradient update on $\beta$}\\
}
\end{algorithm}
\vspace{10pt}

\begin{thm}[Global convergence] \label{thm:ohblook}
  Suppose $f$ is $L$-smooth. Choose $\omega = 3 L$, $\tau = 16 L^2$ and online gradient stepsizes $\eta_P = \tfrac{\eta_{\beta}}{L^2} = \frac{1}{2L}$. Then for any benchmark stepsize $\hat{P} \in \Rbb^{n \times n}$ and momentum parameter $\hat{\beta} \in \Rbb$, 
  \begin{align}
    f (x^{K + 1}) - f^{\star} & \leq \frac{f (x^1) - f^{\star}}{K V
    \max \{ \frac{1}{K} \textstyle \sum_{k = 1}^K - h_{z^k} (\hat{P}, \hat{\beta}) - \frac{L}{4K} \| (P_1 - \hat{P}, \beta_1 -
    \hat{\beta}) \|^2, 0 \} + 1}, \tag{convex}\\
    f (x^{K + 1}) - f^{\star} & \leq [f (x^1) - f^{\star}] ( 1 -
    \mu \max \{ \tfrac{1}{K} \textstyle \sum_{k = 1}^K - h_{z^k} (\hat{P},
    \hat{\beta}) - \tfrac{L}{4 K} \| (P_1 - \hat{P}, \beta_1
    - \hat{\beta}) \|^2, 0 \}  )^K . \tag{$\mu$-strongly convex}
  \end{align}
  where $V$ and $\Delta$ are constants defined in \Cref{thm:heavyball-reduction}.
In particular, with $(P_1, \beta_1) = (\frac{1}{4L} I, \frac{1}{2})$, we have
  \begin{align}
    f (x^{K + 1}) - f^{\star} & \leq \tfrac{8 L[f (x^1) - f^{\star}]}{K V
}, \tag{convex}\\
    f (x^{K + 1}) - f^{\star} & \leq [f (x^1) - f^{\star}] ( 1 -
    \tfrac{1}{8\kappa})^K. \tag{$\mu$-strongly convex}
  \end{align}
\end{thm}

\input{mpproofs/proof-hb-lookahead.tex}

%% file: proofs/proof-hb-lookahead.tex
Recall that the per iteration progress of {\ohb} is
\[ b_k \assign \tfrac{\varphi_{\omega} (z^{k + 1}) - \varphi_{\omega}
   (z^k)}{\| \nabla f (z^k) \|^2 + \frac{\tau}{2} \| z_1^k -
   z_2^k \|^2} . \]
The scheduler proposes the state
\[ z^{k + 1 / 2} \assign (z^k)^{+} = (z_1^k - P_k \nabla f (z^k) +
   \beta_k (z_1^k - z_2^{k}), z_1^k) . \]
{\ohblook} determines the next state $z^{k + 1}$ such that
\begin{equation} \label{eqn:mono-look-action}
   \varphi_\omega (z^{k + 1}) \leq  \min \{\varphi_\omega ( \mathcal{M}_{\text{look}}(z^{k + 1 / 2})), \varphi_\omega(z^k) \},
\end{equation}
where $\Mlook (z) \assign ( z_1 - \tfrac{1}{L + \omega} \nabla_1 \varphi_{\omega} ( z), z_2 )$
The next lemma quantifies the improvement on iteration progress by this landscape
action.

\begin{lem} \label{lem:hb-progress}
Suppose $\tau \geq 2 L^2$. Then the monotone lookahead landscape action \eqref{eqn:mono-look-action} guarantees
\begin{equation} \label{eqn:hb-progress}
    b_k \leq \min \{ h_{z^k} (P_k, \beta_k) - \tfrac{1}{2 (L + \omega)} \| \nabla h_{z^k} (P_k,
    \beta_k) \|_{\ast}^2, 0 \} .
\end{equation}
\end{lem}
\begin{proof}
The monotone landscape action guarantees
$\varphi_{\omega} (z^{k + 1}) \leq \min \{ \varphi_{\omega} (
    \Mlook (z^{k + 1 / 2}) ), \varphi_{\omega} (z^k)
    \},$
and hence
\begin{equation} \label{eqn:pf-hb-progress-1}
b_k \assign \tfrac{\varphi_{\omega} (z^{k + 1}) - \varphi_{\omega}
    (z^k)}{\| \nabla f (z^k) \|^2 + \frac{\tau}{2} \| z_1^k -
    z_2^k \|^2} \leq \min \{ \tfrac{\varphi_{\omega}
    ( \Mlook (z^{k + 1 / 2}) ) - \varphi_{\omega}
    (z^k)}{\| \nabla f (z^k) \|^2 + \frac{\tau}{2} \| z_1^k -
    z_2^k \|^2}, 0 \} .
\end{equation}
Descent lemma on $(L +
\omega)$-smooth potential $\varphi_{\omega} (\cdot, z_2^{k+1/2})$ implies
\begin{align}
\varphi_{\omega} \big( \Mlook (z^{k + 1 / 2}) \big) &
\leq \varphi_{\omega} (z^{k + 1 / 2}) - \tfrac{1}{2 (L + \omega)} \|
\nabla_1 \varphi_{\omega} (z^{k + 1 / 2}) \|^2
\nonumber\\
& = \varphi (z^{k + 1 / 2}) - \tfrac{1}{2 (L + \omega)} \| \nabla f (z^{k
+ 1 / 2}) + \omega (z_1^{k + 1 / 2} - z_2^{k + 1 / 2}) \|^2 . \nonumber
\end{align}
Subtracting both sides by $\varphi_{\omega} (z^k)$ and dividing by $\|
\nabla f (z^k) \|^2 + \tfrac{\tau}{2} \| z_1^k - z_2^k
\|^2$, we obtain
\begin{align}
\tfrac{\varphi_{\omega} ( \mathcal{M}_{\text{look}} (z^{k + 1 / 2})
) - \varphi_{\omega} (z^k)}{\| \nabla f (z^k) \|^2 + \frac{\tau}{2}
\| z_1^k - z_2^k \|^2} & \leq \tfrac{\varphi_{\omega}
(z^{k + 1 / 2}) - \varphi_{\omega} (z^k)}{\| \nabla f (z^k) \|^2 +
\frac{\tau}{2} \| z_1^k - z_2^k \|^2} - \tfrac{1}{2 (L
+ \omega)} \tfrac{\| \nabla f (z^{k + 1 / 2}) + \omega (z_1^{k + 1 / 2} - z_2^{k + 1 / 2}) \|^2}{\| \nabla f (z^k) \|^2 + \frac{\tau}{2} \| z_1^k -
z_2^k \|^2} \nonumber\\
& = h_{z^k} (P_k, \beta_k) - \tfrac{1}{2 (L + \omega)} \| \nabla h_{z^k} (P_k, \beta_k) \|_{\ast}^2,
\label{eqn:pf-hb-progress-2}
\end{align}
where \eqref{eqn:pf-hb-progress-2} uses the inequality of gradient norm from \eqref{eqn:pf-property-3} and that $\tau \geq 2L^2$:
\[ \| \nabla h_{z^k} (P_k, \beta_k) \|_{\ast}^2 \leq \tfrac{\| \nabla f (z^{k
+ 1 / 2}) + \omega (z_1^{k + 1 / 2} - z_2^{k + 1 / 2}) \|^2}{\| \nabla f (z^k) \|^2 + \frac{\tau}{2} \| z_1^k
    - z_2^k \|^2} . \]
Combining the two inequalities \eqref{eqn:pf-hb-progress-1} and \eqref{eqn:pf-hb-progress-2} completes the proof.
\end{proof}

To bound the cumulative feedback $\sum_{k = 1}^K h_{z^k} (P_k, \beta_k)$, we
use the analysis of online gradient descent.

\begin{lem} \label{lem:hb-ogd}
  Consider online gradient descent on heavy-ball feedback
  \begin{align}
    P_{k + 1} & = \Pi_{\mathcal{P}} [P_k - \eta_P \nabla h_{z^k} (P_k,
    \beta_k)], \nonumber\\
    \beta_{k + 1} & = \Pi_{\mathcal{B}} [\beta_k - \eta_{\beta} \nabla h_{z^k}
    (P_k, \beta_k)], \nonumber
  \end{align}
  with $\eta_P = \eta$ and $\eta_{\beta} = L^2 \eta$ for some $\eta > 0$. Then
  for any $(\hat{P}, \hat{\beta}) \in \mathcal{P} \times \mathcal{B}$, the
  sequence $\{ (P_k, \beta_k) \}$ satisfies
  \[ \textstyle \sum_{k = 1}^K h_{z^k} (P_k, \beta_k) \leq \sum_{k = 1}^K h_{z^k}
     (\hat{P}, \hat{\beta}) + \tfrac{1}{2 \eta} \| (P_1 - \hat{P}, \beta_1 -
     \hat{\beta}) \|^2 + \tfrac{\eta}{2} \sum_{k = 1}^K \| \nabla h_{z^k}
     (P_k, \beta_k) \|^2_{\ast} . \]
\end{lem}
\begin{proof}
  The proof is similar to the analysis of online gradient descent under
  Frobenius norm. Note that our algorithm finally takes $\Pcal = \Rbb^{n \times n}$, and the analysis here is more general.  Observe that
  \begin{align}
    \| (P_{k + 1} - \hat{P}, \beta_{k + 1} - \hat{\beta}) \|^2 & = \|
    \Pi_{\mathcal{P}} [P_k - \eta_P \nabla_P h_{z^k} (P_k, \beta_k)] - \hat{P}
    \|^2_F + \tfrac{1}{L^2} | \Pi_{\mathcal{B}} [\beta_k - \eta_{\beta}
    \nabla_{\beta} h_{z^k} (P_k, \beta_k)] - \hat{\beta} |^2 \nonumber\\
    & \leq \| P_k - \hat{P} - \eta_P \nabla_P h_{z^k} (P_k, \beta_k) \|^2_F +
    \tfrac{1}{L^2} | \beta_k - \hat{\beta} - \eta_{\beta} \nabla_{\beta}
    h_{z^k} (P_k, \beta_k) |^2 \nonumber\\
    & = \| P_k - \hat{P} \|^2_F - 2 \eta_P \langle P_k - \hat{P}, \nabla_P
    h_{z^k} (P_k, \beta_k) \rangle + \eta^2_P \| \nabla_P h_{z^k} (P_k,
    \beta_k) \|^2 \nonumber\\
    & \quad + \tfrac{1}{L^2} | \beta_k - \hat{\beta} |^2 - 2
    \tfrac{\eta_{\beta}}{L^2} \langle \beta_k - \hat{\beta}, \nabla_{\beta}
    h_{z^k} (P_k, \beta_k) \rangle + \tfrac{\eta^2_{\beta}}{L^2} |
    \nabla_{\beta} h_{z^k} (P_k, \beta_k) |^2 \nonumber\\
    & = \| (P_k - \hat{P}, \beta_k - \hat{\beta}) \|^2 + \eta^2 [\| \nabla_P
    h_{z^k} (P_k, \beta_k) \|^2 + L^2 | \nabla_{\beta} h_{z^k} (P_k, \beta_k)
    |^2] \nonumber\\
    & \quad - 2 \eta [\langle P_k - \hat{P}, \nabla_P h_{z^k} (P_k, \beta_k)
    \rangle + \langle \beta_k - \hat{\beta}, \nabla_{\beta} h_{z^k} (P_k,
    \beta_k) \rangle] \nonumber \\
    & \leq \| (P_k - \hat{P}, \beta_k - \hat{\beta}) \|^2 + \eta^2 \| \nabla
    h_{z^k} (P_k, \beta_k) \|_{\ast}^2 - 2 \eta [h_{z^k} (P_k, \beta_k) -
    h_{z^k} (\hat{P}, \hat{\beta})], \nonumber
  \end{align}
  which rearranges to
  \[ h_{z^k} (P_k, \beta_k) \leq h_{z^k} (\hat{P}, \hat{\beta}) + \tfrac{1}{2
     \eta} [\| (P_k - \hat{P}, \beta_k - \hat{\beta}) \|^2 - \| (P_{k + 1} -
     \hat{P}, \beta_{k + 1} - \hat{\beta}) \|^2] + \tfrac{\eta}{2} \| \nabla
     h_{z^k} (P_k, \beta_k) \|_{\ast}^2 . \]
  Telescope from $k = 1, \ldots, K$ to obtain
  \begin{align}
    \textstyle \sum_{k = 1}^K h_{z^k} (P_k, \beta_k) & \leq \textstyle \sum_{k = 1}^K h_{z^k}
    (\hat{P}, \hat{\beta}) + \tfrac{1}{2 \eta} \| (P_1 - \hat{P}, \beta_1 -
    \hat{\beta}) \|^2 + \tfrac{\eta}{2} \textstyle \sum_{k = 1}^K \| \nabla h_{z^k} (P_k,
    \beta_k) \|_{\ast}^2 . \nonumber
  \end{align}
  \end{proof}

\paragraph{Proof of \Cref{thm:ohblook}.}
Note that our choice of $\tau = 16 L^2$ satisfies the condition $\tau \geq 2 L^2$ in \Cref{lem:hb-progress}. 
By \Cref{lem:hb-progress} and \Cref{lem:hb-ogd}, the cumulative progress of
{\ohblook} with $\eta \leq \tfrac{1}{L + \omega}$ can be bounded by
\begin{align}
\textstyle \sum_{k = 1}^K b_k & \leq \textstyle \sum_{k = 1}^K \min
  \{ h_{z^k} (P_k, \beta_k) - \tfrac{1}{2 (L + \omega)} \| \nabla
  h_{z^k} (P_k, \beta_k) \|_{\ast}^2, 0 \} \tag{by \Cref{lem:hb-progress}}\\
  & \leq \min \{ \textstyle \sum_{k = 1}^K [ h_{z^k} (P_k,
  \beta_k) - \tfrac{\eta}{2} \| \nabla h_{z^k} (P_k, \beta_k) \|_{\ast}^2
  ], 0 \} \tag{by $\eta \leq \tfrac{\tau}{2 L^2 (L + \omega)}$}\\
  & \leq \min \{ \textstyle \sum_{k = 1}^K h_{z^k} (\hat{P},
  \hat{\beta}) + \tfrac{1}{2 \eta} \| (P_1 - \hat{P}, \beta_1 - \hat{\beta})
  \|^2, 0 \} . \tag{by \Cref{lem:hb-ogd}}
\end{align}
Choose the online gradient stepsize $\eta = \tfrac{1}{L + \omega}$ and denote
\begin{equation*}
    \rho_K(\hat{P}, \hat{\beta}) \assign \tfrac{1}{2 \eta} \| (P_1 - \hat{P}, \beta_1 - \hat{\beta}) \|^2
\end{equation*}
Then the bound on average progress of {\ohblook} becomes
\begin{equation} \label{eqn:pf-hb-progress-3}
    \tfrac{1}{K} \textstyle \sum_{k = 1}^K b_k \leq \min \{ \textstyle \sum_{k = 1}^K h_{z^k} (\hat{P}, \hat{\beta}) + \tfrac{\rho_K(\hat{P}, \hat{\beta})}{K}, 0 \} .
\end{equation}
Our choice of parameters $\omega = 3L$ and $\tau = 16 L^2$ also satisfies the condition $\tau \geq 2 L \omega$ in \Cref{thm:heavyball-reduction}.
Plug in the bound \eqref{eqn:pf-hb-progress-3} into the heavy-ball reductions \eqref{eqn:hb-reduction-cvx}--\eqref{eqn:hb-reduction-strongly-cvx} (\Cref{thm:heavyball-reduction}) to obtain
\begin{align}
\text{(convex)} \quad
f (z^{K + 1}) - f^{\star} & \leq \tfrac{f (z^1) - f^{\star}}{K V
\max \big\{ \frac{1}{K} \sum_{k = 1}^K - h_{z^k} (\hat{P}, \hat{\beta}) - \frac{\rho_K(\hat{P}, \hat{\beta})}{K}, 0 \big\} + 1}, \label{eqn:pf-hb-progress-4}\\
(\mu\text{-strongly convex)} \quad
f (z^{K + 1}) - f^{\star} & \leq [f (z^1) - f^{\star}] \big( 1 -
\mu \max \big\{ \tfrac{1}{K} \textstyle \sum_{k = 1}^K - h_{z^k} (\hat{P},
\hat{\beta}) - \frac{\rho_K(\hat{P}, \hat{\beta})}{K}, 0 \big\}  \big)^K . \label{eqn:pf-hb-progress-5}
\end{align}

\Cref{lem:heavyball-potential} with $\alpha = \tfrac{1}{4 L}$ and $\beta = \tfrac{1}{2}$ implies $\omega = 3 L$ and guarantees
\begin{align*}
    \varphi_w \big( (z^k)^+ (\tfrac{1}{4L} I, \tfrac{1}{2}) \big) 
    &= \varphi_w \big(z_1^k - \tfrac{1}{4L} \nabla f(z^k) + \tfrac{1}{2}(z_1^k - z_2^k), z_1^k \big) \\
    &\leq \varphi_w (z^k) - \tfrac{1}{8 L} [ \| \nabla f (z^k) \|^2 + 8 L^2 \| z_1^k - z_2^k \|^2 ] \\
    &= \varphi_w (z^k) - \tfrac{1}{8 L} [ \| \nabla f (z^k) \|^2 + \tfrac{\tau}{2} \| z_1^k - z_2^k \|^2 ],
    \tag{by $\tau = 16 L^2$}
\end{align*}
which rearranges to
\begin{equation} \label{eqn:pf-hb-hindsight} 
    h_{z^k} ( \tfrac{1}{4 L} I, \tfrac{1}{2} ) = \tfrac{\varphi_{\omega} ( (z^k)^+ (\frac{1}{4L} I, \frac{1}{2}) ) - \varphi_{\omega} (z^k)}{\| \nabla f (z^k) \|^2 + \frac{\tau}{2} \| z_1^k - z_2^k \|^2} \leq -\tfrac{1}{8 L} \text{, for all $k$.}
\end{equation}
Plug in $(\hat{P}, \hat{\beta}) = ( \tfrac{1}{4 L} I, \tfrac{1}{2}
)$ into \eqref{eqn:pf-hb-progress-4} and \eqref{eqn:pf-hb-progress-5} and further bound by \eqref{eqn:pf-hb-hindsight} completes the proof.

%% file: sec_nonconvex.tex
\section{{\os} for smooth nonconvex optimization} \label{sec:nonconvex}
{\os} can be extended to smooth nonconvex optimization problems, but requires modifying the feedback.
The analysis of {\os} in Part I \cite{gao2025gradient} and \Cref{sec:prac} assumes a convex objective $f$, which guarantees convexity of feedback and enables hypergradient \cite[Theorem 4.2]{gao2025gradient} and heavy-ball (\Cref{thm:heavyball-reduction}) progress reduction.
One possible approach to handle nonconvex $f$ is to apply the proximal point method {\cite{chen2024open}} with a sufficiently large regularizer $\gamma>L$: 
\begin{equation} \label{eqn:proxpt}
	x^{k + 1} = \argmin_x ~ \{ f (x) + \tfrac{\gamma}{2} \| x - x^k \|^2 \}.
\end{equation} 
When $f$ is nonconvex but $L$-smooth, the proximal subproblem \eqref{eqn:proxpt} is $(L - \gamma)$-weakly convex and becomes strongly convex if $\gamma > L$. 
Then {\os} can serve as a subroutine to solve the subproblem \eqref{eqn:proxpt}. However, this proximal point approach requires a double-loop algorithm: an outer loop updates the proximal center $x^k$ and an inner loop solves the proximal subproblem with {\os}. In contrast, this section introduces a different approach to smooth nonconvex optimization by stepsize regularization,
avoiding the need for a double-loop algorithm. Observe that $\tfrac{1}{L} I$ is a safe stepsize in smooth nonconvex optimization. 
Hence, we can regularize the feedback function $h_x (P)$ by  $\| P - \tfrac{1}{L} I \|_F^2$.
This regularizer vanishes when the stepsize $\hat{P}= \frac{1}{L} I$ is used as a benchmark, and an online learning method using this regularized feedback 
at least competes with the performance of $\tfrac{1}{L} I$.
Regularization in stepsize space requires only a fixed center $\tfrac{1}{L} I$ and avoids the second loop.

\paragraph{Structure of the section.}
The roadmap of this section resembles Part I: \Cref{sec:reg-feedback} introduces the regularized feedback and establishes its analytical properties; \Cref{sec:landscape-progress} introduces two progress reductions for general nonconvex and PL functions and several landscape actions. Finally, in \Cref{sec:nonconvex-algo} we design and analyze {\oh} variants for nonconvex optimization. For simplicity, we will again focus on the hypergradient feedback, while the results for ratio feedback hold similarly for nonconvex PL functions.

\subsection{Regularized feedback design and analysis} \label{sec:reg-feedback}

Given a regularization parameter $\lambda \geq 0$, define the regularized hypergradient feedback
\begin{equation} \label{eqn:def-regfeedback}
	h_x^{\lambda} (P) \assign h_x (P) + \tfrac{\lambda}{2} \| P -
  \tfrac{1}{L} I \|_F^2 = \tfrac{f (x - P \nabla f (x)) - f (x)}{\|
  \nabla f (x) \|^2} + \tfrac{\lambda}{2} \| P - \tfrac{1}{L}
  \|_F^2.
\end{equation}
The regularization parameter $\lambda$ is chosen to neutralize the \textit{local nonconvexity} of $h_x$ for $P \in \Pcal$. 
Since $h_x (P)$ is $L$-smooth, choosing $\lambda = L$ suffices since it guarantees $h^\lambda_x$ is globally convex.
However, this strong regularizer forces $P$ to stay close to the conservative stepsize $\tfrac{1}{L} I$. 
In practice, we want to choose $\lambda$ as small as possible to retain adaptivity, especially when $x$ enters a locally convex region. 
In other words, we need an estimate for the local weak convexity constant of $h_x (P)$. 
With additional assumptions on $f$, it is possible.

\begin{lem}[Nonconvexity of $h_x$]
  \label{lem:hx-nonconvexity}For any $x \nin \mathcal{S}^{\star} \assign \{ x
  : \| \nabla f (x) \| = 0 \}$, $h_x (P)$ has the following properties:
  \begin{enumerate}[leftmargin=15pt]
    \item If $f$ is $L$-smooth, then $h_x (P)$ has weak convexity constant
    $\rho \leq \rho_1 \assign L$.
    
    \item If $f$ is $L$-smooth, has $H$-Lipschitz Hessian and $\diam
    (\mathcal{P}) \leq D$, then $h_x (P)$ has weak convexity constant
    \[ \rho \leq \rho_2 \assign \max \{ - \lambda_{\min} (\nabla^2 h_x (0 \cdot I)) + H D \|
       \nabla f (x) \|, 0 \} \]
 	when $P \in \Pcal$. 
    \item Moreover, if $f$ is also $\mu$-PL, then
    $\rho \leq \rho_3 \assign (\tfrac{H}{\mu} + H
    D )\| \nabla f (x) \|$.
  \end{enumerate}
\end{lem}

The intuition behind \Cref{lem:hx-nonconvexity} is simple: when $x$ is close to stationarity and $\diam (\Pcal) < \infty$, then $\nabla^2 h_x(P)$ does not change much within $P \in \Pcal$. Hence, we can bound the worst-case weak convexity constant of $h_x(P)$ for any $P \in \Pcal$ by evaluating 
$\lambda_{\min} (\nabla^2 h_x (\bar{P}))$ for some fixed $\bar{P} \in \Pcal$. 
In particular, if $f$ is $\mu$-PL,  $\lambda_{\min} (\nabla^2 h_x (\bar{P}))$ can also be bounded using $\| \nabla f(x) \|$. 
\Cref{lem:hx-nonconvexity} takes $\bar{P} = 0 \cdot I$.
\begin{rem}
  Computing $\lambda_{\min} (\nabla^2 h_x (\bar{P}))$ can be computationally
  expensive in practice. One exception is when $P = \alpha I$ is a scalar
  multiple of the identity.
  In this case the second-order derivative of $h_x(\alpha I)$ with respect
  to $\alpha$ is $\langle \tfrac{\nabla f (x)}{\| \nabla f (x) \|},
  \nabla^2 f (x - \alpha \nabla f(x)) \tfrac{\nabla f (x)}{\| \nabla f (x)
  \|} \rangle$ and can be computed by evaluating the Hessian-vector product,
  without forming the Hessian matrix \cite{pearlmutter1994fast}.
\end{rem}

Since $h_x^{\lambda} (P)$ adds a strongly convex regularizer to $h_x(P)$, its analytic properties follow from those of $h_x (P)$.

\begin{lem}[Properties of $h_x^{\lambda}$]
  \label{lem:property-hxreg}For any $x \nin \mathcal{S}^{\star}$, the
  regularized feedback $h_x^{\lambda} (P)$ has gradient
  \[ \nabla h_x^{\lambda} (P) = - \tfrac{\nabla f (x - P \nabla f (x)) \nabla
     f (x)^{\top}}{\| \nabla f (x) \|^2} + \lambda ( P - \tfrac{1}{L} I
     ) . \]
  and $h_x^{\lambda} ( \tfrac{1}{L} I ) \leq - \tfrac{1}{2 L}$ for
  all $\lambda$.  Moreover, if $f$ is $L$-smooth, has $H$-Lipschitz Hessian and $\diam
  (\mathcal{P}) \leq D$, then $h_x^{\lambda} (P)$ is convex and $(L + \lambda)$-smooth
  for all $\lambda \geq \min \{ \rho_1, \rho_2, \rho_3 \}$.
\end{lem}

\input{../proofs/proof-property-hxreg.tex}

Finally, the definition of minimax stepsize does not immediately carry over to the nonconvex case; therefore, we use $\frac{1}{L} I$ as the benchmark stepsize in the theoretical analyses.

\subsection{Nonconvex progress reduction and landscape action} \label{sec:landscape-progress}
Recall the per iteration progress associated with the hypergradient feedback $h_k \assign \tfrac{f (x^{k + 1}) - f (x^k)}{\| \nabla f (x^k) \|^2}$. \Cref{thm:nonconvex-reduction} reduces the stationarity/optimality guarantees to the cumulative progress. Notably, these two reductions do not assume the convexity of $f$.
\begin{thm}[Nonconvex reductions] \label{thm:nonconvex-reduction}
  Let $\{ x^k \}$ be any sequence such that $x^k \nin \mathcal{S}^{\star}$ and
  $f (x^{k + 1}) \leq f (x^k)$. Then
\begin{equation} \label{eqn:nonconvex-reduction-gradnorm}
	\min_{1 \leq k \leq K}  \| \nabla f (x^k) \|^2 \leq \tfrac{f (x^1) -  f^{\star}}{K} \tfrac{1}{\frac{1}{K} \sum_{k = 1}^K - h_k}.
\end{equation}
  Moreover, if $f$ is $\mu$-PL, then
\begin{equation} \label{eqn:nonconvex-reduction-PL}
	\textstyle  f (x^{K + 1}) - f^{\star} \leq [f (x^1) - f^{\star}] ( 1 - \frac{2\mu}{K} \sum_{k = 1}^K -h_k )^K .
\end{equation}
\end{thm}

\input{../proofs/proof-thm-nonconvex-reduction.tex}

\paragraph{Landscape actions $\Mcal$ and progress.} As in Part I, the landscape action $\mathcal{M} (x^{k + 1 / 2}, x^k)$ determines the choice of the next iterate in {\os}. Recall the two landscape actions in {\oh}:
\begin{itemize}[leftmargin=15pt]
  \item \textit{Monotone.} $x^{k + 1}$ satisfies $f (x^{k + 1}) \leq \min \{
  f (x^{k + 1 / 2}), f (x^k) \}$.
  \item {\textit{Monotone lookahead}}. $x^{k + 1}$ satisfies $f (x^{k + 1}) \leq
  \min \{ f ( x^{k + 1 / 2} - \tfrac{1}{L} \nabla f (x^{k + 1 / 2})
  ), f (x^k) \}$.
\end{itemize}
For each action, we associate the progress and feedback in \Cref{lem:nonconvex-feedback-progress}.
\begin{lem}[Feedback and progress]  \label{lem:nonconvex-feedback-progress}
Let $f$ be $L$-smooth and suppose
  $\lambda_k \leq L$. Then each of the above two landscape actions guarantees the following relation between the feedback and per iteration progress:
  \begin{itemize}[leftmargin=15pt]
    \item Monotone. $h_k \leq \min \{ h^{\lambda_k}_{x^k} (P_k) -
    \tfrac{\lambda_k}{2} \| P_k - \tfrac{1}{L} \|_F^2, 0 \}$.
    
    \item Monotone lookahead. $h_k \leq \min \{ h^{\lambda_k}_{x^k} (P_k) -
    \tfrac{1}{4 L} \| \nabla h^{\lambda_k}_{x^k} (P_k) \|_F^2, 0 \}$.
  \end{itemize}
\end{lem}

\input{../proofs/proof-lem-nonconvex-feedback-progress.tex}

\subsection{Algorithm design and analysis} \label{sec:nonconvex-algo}

Putting the results together, we instantiate {\osgmhandsonhx} for nonconvex optimization with
\[\ell_x(P)\assign h_x^\lambda(P),\quad \text{Monotone Lookahead landscape: } f(x^{k+1}) \leq \min \{f(x^{k+1/2} - \tfrac{1}{L} \nabla f(x^{k+1/2})), f(x^k)\} \quad \Acal\assign \ogd. \]

The algorithm is {\osgmhandsonhx} with regularization (\Cref{alg:nonconvexosgm}).

\vspace{10pt}

\begin{algorithm}[H]
{\textbf{input} initial point $x^1$, initial stepsize $P_1 \in \Pcal$, online gradient stepsize $\eta > 0$}

\For{$k = 1, 2, \dots$}{
$x^{k + 1/2} = x^{k} - P_k \nabla f(x^k)$\\
Choose $x^{k+1}$ that satisfies $f(x^{k+1}) \leq \min \{ f(x^{k+1/2} - \frac{1}{L} \nabla f(x^{k+1/2})
   ), f(x^k) \}$\\
Choose $0 \leq \lambda_k \leq L$ based on \Cref{lem:hx-nonconvexity} such that $h_{x^k}^{\lambda_k}(P)$ is convex within $\Pcal$\\
$P_{k+1} = \Pi_\Pcal [P_k - \eta \nabla h_{x^k}^{\lambda_k}(P_k)]$ \label{alg:osgmhandsonhx-line} 
}
\caption{{\osgmhandsonhx} for nonconvex optimization\label{alg:nonconvexosgm}}
\end{algorithm}
\vspace{10pt}

Recall that the gradient of the regularized hypergradient feedback in \Cref{alg:osgmhandsonhx-line} of \Cref{alg:nonconvexosgm} takes the form
\[ \nabla h_x^{\lambda} (P) = - \tfrac{\nabla f (x - P \nabla f (x)) \nabla
     f (x)^{\top}}{\| \nabla f (x) \|^2} + \lambda ( P - \tfrac{1}{L} I
     ) \]
     
and $\Pcal$ can be chosen as a simple bounded set (e.g., an element-wise box) to allow inexpensive orthogonal projection. Due to the difficulty of the nonconvex setting, we restrict our analysis to global convergence only.

\begin{thm}[Convergence]
  \label{thm:nonconvex-global-conv}Let $f$ be $L$-smooth. For any benchmark
  stepsize $\hat{P} \in \mathcal{P}$, {\osgmhandsonhx} (\Cref{alg:nonconvexosgm}) with $\eta = 1 / (4 L)$ satisfies
  \begin{align}
    \min_{1 \leq k \leq K}  \| \nabla f (x^k) \|^2 
    \leq{} & \tfrac{f (x^1) - f^{\star}}{K} \tfrac{1}{\frac{1}{K} \sum_{k =
    1}^K [ - h_{x^k} (\hat{P}) - \frac{\lambda_k}{2} \| \hat{P} -
    \frac{1}{L} I \|_F^2 ] - \frac{L}{K} \| P_1 - \hat{P}
    \|_F^2} \tag{nonconvex} \\
     f (x^{K + 1}) - f^{\star}
    \leq{} & \textstyle [f (x^1) - f^{\star}] ( 1 - 2 \mu \max \{
    \tfrac{1}{K} \sum_{k = 1}^K [ - h_{x^k} (\hat{P}) -
    \tfrac{\lambda_k}{2} \| \hat{P} - \tfrac{1}{L} I \|_F^2
    ] - \tfrac{L}{K} \| P_1 - \hat{P} \|_F^2, 0 \} )^K \tag{$\mu$-PL}.
    \nonumber
  \end{align}  
In particular, if $P_1 = \tfrac{1}{L} I$, then
\begin{align}
\min_{1 \leq k \leq K}  \| \nabla f (x^k) \|^2 \leq{} & \tfrac{2 L [f (x^1)
    - f^{\star}]}{K} \tag{nonconvex}\\
f (x^{K + 1}) - f^{\star} \leq{} & [f (x^1) - f^{\star}] ( 1 -
    \tfrac{1}{\kappa} )^K \tag{$\mu$-PL}. \nonumber
  \end{align}
\end{thm}

\input{../proofs/proof-nonconvex-global-conv.tex}

\Cref{thm:nonconvex-global-conv} does not show an improved convergence rate compared to gradient descent using constant stepsize $\frac{1}{L}I$. However, when the iterates enter a locally convex region (e.g., approaching a local minimum), $\lambda_k$ will automatically vanish, and we recover adaptivity in the convex case.

%% file: proofs/proof-property-hxreg.tex
The expression of gradient follows by definition $\nabla h_x^{\lambda} (P) =
\nabla h_x (P) + \lambda ( P - \tfrac{1}{L} I )$. $h_x^{\lambda}
( \tfrac{1}{L} I ) = h_x ( \tfrac{1}{L} I ) \leq -
\tfrac{1}{2 L}$ follows from the descent lemma. Using \Cref{lem:hx-nonconvexity}, local convexity of $h_x^\lambda(P)$ follows immediately and the smoothness constant follows by adding $L$, the smoothness of $h_x(P)$, and $\lambda$.

%% file: proofs/proof-thm-nonconvex-reduction.tex
These reductions have appeared in \cite{gao2024gradient}, and we repeat the proofs for completeness. The first relation \eqref{eqn:nonconvex-reduction-gradnorm} follows by definition:
\begin{align}
  ( \textstyle \sum_{k = 1}^K - h_k ) \cdot \displaystyle \min_{1 \leq k \leq K} \|
  \nabla f (x^k) \|^2 
  \leq{} & \textstyle \sum_{k = 1}^K - h_k \| \nabla f (x^k) \|^2 \tag{by $-h_k \geq 0$}\\
  ={} & - [ \textstyle \sum_{k = 1}^K f (x^{k + 1}) - f (x^k) ] \tag{by definition of $h_k$}\\
  ={} & f (x^1) - f (x^{K + 1}) \nonumber\\
  \leq{} & f (x^1) - f^{\star} . \tag{by $f(x^{K+1}) \leq f^\star$}
\end{align}
and the second relation \eqref{eqn:nonconvex-reduction-PL} follows from:
\begin{align}
  \tfrac{f (x^{K + 1}) - f^{\star}}{f (x^1) - f^{\star}} \leq{} & (
  \tfrac{1}{K} \textstyle \sum_{k = 1}^K \tfrac{f (x^{k + 1}) - f^{\star}}{f (x^k) - f^{\star}} )^K \tag{by AM-GM inequality}\\
  ={} & ( 1 + \tfrac{1}{K} \textstyle \sum_{k = 1}^K \tfrac{f (x^{k + 1}) - f (x^k)}{f
  (x^k) - f^{\star}} )^K \nonumber\\
  ={} & ( 1 + \tfrac{1}{K} \textstyle \sum_{k = 1}^K \tfrac{h_k \| \nabla f (x^k)
  \|^2}{f (x^k) - f^{\star}} )^K \tag{by definition of $h_k$}\\
  \leq{} & ( 1 - \tfrac{2 \mu}{K} \textstyle \sum_{k = 1}^K - h_k )^K \tag{by $h_k \leq 0$ and $\frac{\|\nabla f(x^k)\|^2}{f(x^k) - f^\star} \geq 2\mu $}
\end{align}
and this completes the proof.

%% file: proofs/proof-lem-nonconvex-feedback-progress.tex
Recall that $x^{k + 1} =\mathcal{M} (x^{k + 1 / 2}, x^k)$ and $x^{k + 1 / 2} =
x^k - P_k \nabla f (x^k)$. The per iteration progress and (regularized)
feedback are
\[ h_k = \tfrac{f (x^{k + 1}) - f (x^k)}{\| \nabla f (x^k) \|^2}, \quad
   h_{x^k} (P_k) = \tfrac{f (x^{k + 1 / 2}) - f (x^k)}{\| \nabla f (x^k)
   \|^2}, \quad h_{x^k}^{\lambda} (P_k) = h_{x^k} (P_k) + \tfrac{\lambda}{2}
   \| P_k - \tfrac{1}{L} I \|_F^2. \]
\begin{itemize}[leftmargin=10pt]
  \item {\textit{Monotone}} landscape satisfies $f (x^{k + 1}) \leq \min \{ f
  (x^{k + 1 / 2}), f (x^k) \}$ and thus
  \[ h_k = \tfrac{f (x^{k + 1}) - f (x^k)}{\| \nabla f (x^k) \|^2} \leq \min
     \{ h_{x^k} (P_k), 0 \} = \min \{ h_{x^k}^{\lambda} (P_k) -
     \tfrac{\lambda}{2} \| P_k - \tfrac{1}{L} I \|_F^2, 0 \}
     . \]
  \item {\textit{Monotone lookahead}} landscape satisfies $f (x^{k + 1}) \leq
  \min \{ f ( x^{k + 1 / 2} - \tfrac{1}{L} \nabla f (x^{k + 1 / 2})
  ), 0 \}$ and by descent lemma, $f ( x^{k + 1 / 2} -
  \tfrac{1}{L} \nabla f (x^{k + 1 / 2}) ) \leq f (x^{k + 1 / 2}) -
  \tfrac{1}{2 L} \| \nabla f (x^{k + 1 / 2}) \|^2$ and
  \begin{align}
    h_k ={} & \tfrac{f (x^{k + 1}) - f (x^k)}{\| \nabla f (x^k) \|^2}
    \nonumber\\
    \leq{} & \tfrac{f (x^{k + 1 / 2}) - f (x^k)}{\| \nabla f (x^k) \|^2} -
    \tfrac{1}{2 L} \tfrac{\| \nabla f (x^{k + 1 / 2}) \|^2}{\| \nabla f (x^k)
    \|^2} \tag{by descent lemma}\\
    ={} & h_{x^k} (P_k) - \tfrac{1}{2 L} \| \nabla h_{x^k} (P_k) \|_F^2
    \tag{by $\| \nabla h_{x^k} (P_k) \|_F^2 = \tfrac{\| \nabla f (x^{k + 1 / 2}) \|^2}{\| \nabla f (x^k) \|^2}$}\\
    ={} & h_{x^k}^{\lambda} (P_k) - \tfrac{\lambda}{2} \| P_k -
    \tfrac{1}{L} I \|_F^2 - \tfrac{1}{2 L} \| \nabla h_{x^k} (P_k)
    \|_F^2 \nonumber\\
    ={} & h_{x^k}^{\lambda} (P_k) - \tfrac{1}{4 L} \| \nabla h_{x^k}^{\lambda}
    (P_k) \|_F^2 + \tfrac{1}{4 L} [ \| \nabla h_{x^k}^{\lambda} (P_k)
    \|_F^2 - 2 \lambda L \| P_k - \tfrac{1}{L} I \|_F^2 - 2 \|
    \nabla h_{x^k} (P_k) \|_F^2 ] . \nonumber
  \end{align}
Notice that when $\lambda \leq L$,  
  \begin{align}
    \| \nabla h_{x^k}^{\lambda} (P_k) \|_F^2 ={} & \| \nabla h_{x^k} (P_k)
    + \lambda ( P_k - \tfrac{1}{L} I ) \|_F^2 \nonumber\\
    \leq{} & 2 \| \nabla h_{x^k} (P_k) \|_F^2 + 2 \lambda^2 \| P_k -
    \tfrac{1}{L} I \|_F^2 \tag{by $\|P + Q\|_F^2 \leq 2\|P\|_F^2 + 2\|Q\|_F^2$}\\
    \leq{} & 2 \| \nabla h_{x^k} (P_k) \|_F^2 + 2 \lambda L \| P_k -
    \tfrac{1}{L} I \|_F^2 \tag{by $\lambda \leq L$}
  \end{align}
  
which implies $h_k \leq h_{x^k}^{\lambda} (P_k) - \tfrac{1}{4 L} \| \nabla
  h_{x^k}^{\lambda} (P_k) \|_F^2$ and completes the proof.
\end{itemize}

%% file: proofs/proof-nonconvex-global-conv.tex
Recall the regret guarantee of online gradient descent: for any $\hat{P} \in
\mathcal{P}$, we successively deduce that
\begin{align}
  \| P_{k + 1} - \hat{P} \|_F^2 ={} & \| \Pi_\Pcal [P_k - \eta \nabla
  h_{x^k}^{\lambda_k} (P_k)] - \hat{P} \|_F^2 \nonumber \\
  \leq{} & \| P_k - \eta \nabla h_{x^k}^{\lambda_k} (P_k) - \hat{P} \|_F^2 
  \tag{by non-expansiveness}\\
  ={} & \| P_k - \hat{P} \|_F^2 - 2 \eta \langle \nabla h_{x^k}^{\lambda_k}
  (P_k), P_k - \hat{P} \rangle + \eta^2 \| \nabla h_{x^k}^{\lambda_k} (P_k)
  \|^2_F \nonumber\\
  \leq{} & \| P_k - \hat{P} \|_F^2 - 2 \eta [h_{x^k}^{\lambda_k} (P_k) -
  h_{x^k}^{\lambda_k} (\hat{P})] + \eta^2 \| \nabla h_{x^k}^{\lambda_k} (P_k)
  \|^2_F, \nonumber
\end{align}
where the last inequality holds since $P_k, \hat{P} \in \mathcal{P}$ and our choice of $\{ \lambda_k \}$ ensures $h_{x^k}^{\lambda_k} (P)$ is locally convex. Next, we rearrange the terms and obtain
\begin{equation} \label{eqn:proof-nonconvex-global-conv-1}
	h_{x^k}^{\lambda_k} (P_k) - h_{x^k}^{\lambda_k} (\hat{P}) \leq \tfrac{1}{2
   \eta} [\| P_k - \hat{P} \|_F^2 - \| P_{k + 1} - \hat{P} \|_F^2] +
   \tfrac{\eta}{2} \| \nabla h_{x^k}^{\lambda_k} (P_k) \|^2_F .
\end{equation}
Back to the algorithm analysis. First, according to \Cref{thm:nonconvex-reduction}, the 
convergence of {\osgmhandsonhx} follows from bounding the cumulative progress
\begin{align}
  \min_{1 \leq k \leq K}  \| \nabla f (x^k) \|^2 \leq{} & \tfrac{f (x^1) - f^{\star}}{K} \tfrac{1}{\frac{1}{K}  \sum_{k = 1}^K - h_k}, \label{eqn:proof-nonconvex-global-conv-2} \\
  f (x^{K + 1}) - f^{\star} \leq{} & [f (x^1) - f^{\star}] ( 1 -
  \tfrac{2 \mu}{K}  \textstyle \sum_{k = 1}^K - h_k )^K.  \label{eqn:proof-nonconvex-global-conv-3}
\end{align}
Second, the cumulative progress $\textstyle \sum_{k = 1}^K h_k$ under monotone lookahead
action is bounded by \Cref{lem:nonconvex-feedback-progress}
\begin{equation} \label{eqn:proof-nonconvex-global-conv-4}
	\tfrac{1}{K} \textstyle \sum_{k = 1}^K h_k \leq \tfrac{1}{K} \textstyle \sum_{k = 1}^K \min
   \{ h^{\lambda_k}_{x^k} (P_k) - \tfrac{1}{4 L} \| \nabla
   h^{\lambda_k}_{x^k} (P_k) \|_F^2, 0 \} \leq \min \{ \tfrac{1}{K}
   \textstyle \sum_{k = 1}^K h^{\lambda_k}_{x^k} (P_k) - \tfrac{1}{4 L} \| \nabla
   h^{\lambda_k}_{x^k} (P_k) \|_F^2, 0 \} .
\end{equation}
Third, the regret guarantee of online gradient descent \eqref{eqn:proof-nonconvex-global-conv-1} bounds the cumulative
feedback $\textstyle \sum_{k = 1}^K h^{\lambda_k}_{x^k} (P_k)$ by
\begin{equation} \label{eqn:proof-nonconvex-global-conv-5}
	\textstyle \sum_{k = 1}^K h^{\lambda_k}_{x^k} (P_k) \leq \textstyle \sum_{k = 1}^K
   h_{x^k}^{\lambda_k} (\hat{P}) + \tfrac{1}{2 \eta} \| P_1 - \hat{P} \|_F^2 +
   \tfrac{\eta}{2} \textstyle \sum_{k = 1}^K \| \nabla h_{x^k}^{\lambda_k} (P_k) \|^2_F .
\end{equation}
Putting \eqref{eqn:proof-nonconvex-global-conv-4} and \eqref{eqn:proof-nonconvex-global-conv-5} together, the cumulative progress is bounded by
\[ \tfrac{1}{K} \textstyle \sum_{k = 1}^K h_k \leq \min \{ \textstyle \sum_{k = 1}^K
   h_{x^k}^{\lambda_k} (\hat{P}) + \tfrac{1}{2 \eta} \| P_1 - \hat{P} \|_F^2 +
   ( \tfrac{\eta}{2} - \tfrac{1}{4 L} ) \textstyle \sum_{k = 1}^K \| \nabla
   h_{x^k}^{\lambda_k} (P_k) \|^2_F, 0 \} . \]
When $\eta \leq \tfrac{1}{2 L}$, the gradient norm term vanishes. With
$\eta = \tfrac{1}{2 L}$ the relation simplifies to
\[ \tfrac{1}{K} \textstyle \sum_{k = 1}^K h_k \leq \min \{ \textstyle \sum_{k = 1}^K
   h_{x^k}^{\lambda_k} (\hat{P}) + L \| P_1 - \hat{P} \|_F^2, 0 \} =
   \min \{ \textstyle \sum_{k = 1}^K h_{x^k} (\hat{P}) + \tfrac{\lambda}{2} \|
   P - \tfrac{1}{L} \|_F^2 + L \| P_1 - \hat{P} \|_F^2, 0 \} . \]
Plugging the relation back into  \eqref{eqn:proof-nonconvex-global-conv-2} and \eqref{eqn:proof-nonconvex-global-conv-3} completes the proof. Also notice that when $P_1 = \hat{P} = \tfrac{1}{L} I$, $h_x^{\lambda} (\hat{P}) \leq - \tfrac{1}{2 L}$ and $\|
P_1 - \hat{P} \|_F^2 = 0$.

%% file: sec_exp.tex
\section{Numerical experiments} \label{sec:exp}

This section conducts numerical experiments to validate the empirical
performance of hypergradient descent. We compare {\ob} with 
different adaptive optimization algorithms.

\subsection{Experiment setup}

\paragraph{Dataset.} We use 38 datasets from \texttt{LIBSVM} for convex optimization tasks and 47 unconstrained problems from \texttt{CUTEst} dataset \cite{gould2015cutest}.

\paragraph{Algorithm benchmark.}
We benchmark the following algorithms.\\
\begin{itemize}[leftmargin=10pt,itemsep=2pt,topsep=0pt]
    \item \texttt{GD}. Vanilla gradient descent.
    \item \texttt{GD-HB}. Gradient descent with heavy-ball momentum. \cite{polyak1964some}
    \item \texttt{AGD-CVX}. The smooth convex version of accelerated gradient descent (Nesterov momentum). \cite{d2021acceleration}
    \item \texttt{AGD-SCVX}. The smooth strongly convex version of accelerated gradient descent. \cite{d2021acceleration}
    \item \texttt{Adam}. Adaptive momentum estimation. \cite{kingma2014adam}
    \item \texttt{AdaGrad}. Adaptive (sub)gradient method. \cite{duchi2011adaptive}
    \item \texttt{BFGS}. {\bfgs} from \texttt{scipy} \cite{nocedal1999numerical,virtanen2020scipy}.
    \item \texttt{L-BFGS-Mk}. {\lbfgs} with memory size \texttt{k} in \texttt{scipy}.
    \item Practical variant {\ob} uses as memory $7$ vectors of size $n$, comparable to memory for \texttt{L-BFGS-M1}.
\end{itemize}

\paragraph{Algorithm configuration.} 
\begin{itemize}[leftmargin=10pt]
  
  \item Stepsize in \texttt{GD}, \texttt{GD-HB},
  \texttt{AGD-CVX}, and \texttt{AGD-SCVX} are all set to $1 / L$.
  
  \item The momentum parameter in \texttt{GD-HB} is chosen within the set $\{
  0.1, 0.5, 0.9, 0.99 \}$.
  
  \item The \texttt{Adam} stepsize is chosen within the set $\{ 1 / L, 10^{- 3},
  10^{- 2}, 10^{- 1}, 1, 10 \}$. $\beta_1 = 0.9, \beta_2 = 0.999$.
  
  \item  The \texttt{AdaGrad} stepsize is chosen within the set $\{ 1 / L, 10^{- 3},
  10^{- 2}, 10^{- 1}, 1, 10 \}$.
  
  \item {\bfgs}, \texttt{L-BFGS-Mk} use default parameters in
  \texttt{scipy}.
  
  \item {\ob} uses the default parameters without tuning.
\end{itemize}

\paragraph{Testing configurations.}

\begin{enumerate}[leftmargin=15pt,label=\textbf{\arabic*)}]
  \item {\textit{Maximum oracle access.}} For convex problems, we allow a maximum of 1000 gradient
  oracles for each algorithm. For nonconvex problems, we allow a maximum of 2000 gradient oracles.
  
  \item {\textit{Initial point}}. For convex problems, All the algorithms are initialized from the
  same starting point generated from the normal distribution $\mathcal{N} (0,
  I_n)$ and normalized to have unit length. For nonconvex problems, we use the starting point provided by \texttt{CUTEst}.
  
  \item {\textit{Stopping criterion.}} Algorithms stop if $\| \nabla f 
  \|_\infty \leq 10^{-3}$.
\end{enumerate}

\begin{table}[h]
\centering
\caption{Number of solved problems for each algorithm. \label{table:stats}}
\begin{tabular}{cccc}
\toprule
    Algorithm/Problem & SVM (38) $\uparrow$ & Logistic Regression (38) $\uparrow$ & CUTEst (47) $\uparrow$ \\
\midrule
    \texttt{GD} & 5 & 2 & 13 \\
    \texttt{GD-HB} & 13 & 10 & 20 \\
    \texttt{AGD-CVX} & 13 & 6 & 20 \\
    \texttt{AGD-SCVX} & 10 & 9 & 15 \\
    \texttt{Adam} & 31 & 16 & 24 \\
    \texttt{AdaGrad} & 11 & 11 & 18 \\
    \texttt{L-BFGS-M1} & 27 & 18 & 35 \\
    \texttt{L-BFGS-M3} & 33 & 19 & 37 \\
    \texttt{L-BFGS-M5} & 35 & 26 & 37 \\
    \texttt{L-BFGS-M10} & 35 & 31 & 36 \\
    \texttt{BFGS} & 37 & 36 & 38 \\
    {\ob} & 35 & 35 & 34 \\
\bottomrule
\end{tabular}

\end{table}

\subsection{Convex optimization}
\begin{figure}[!h]
\centering
\includegraphics[scale=0.17]{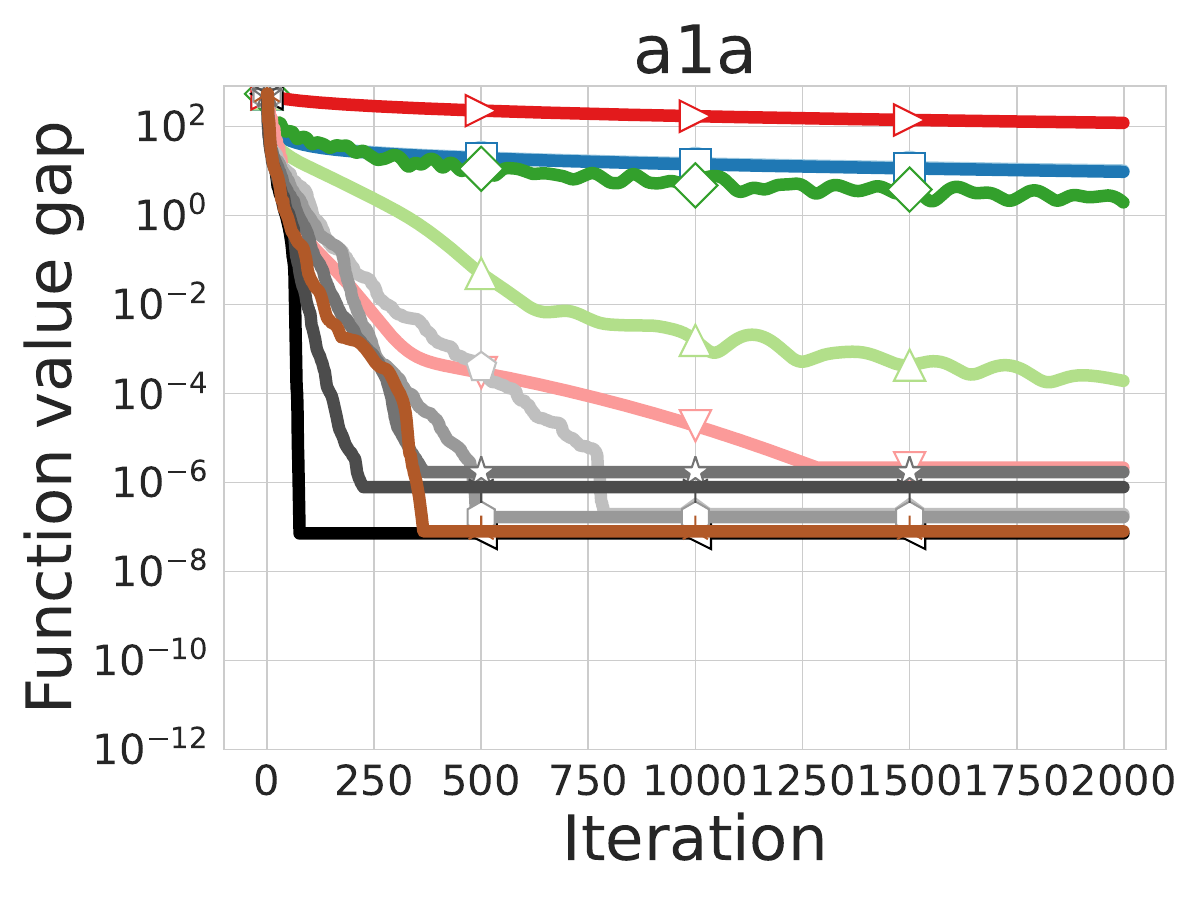}
\includegraphics[scale=0.17]{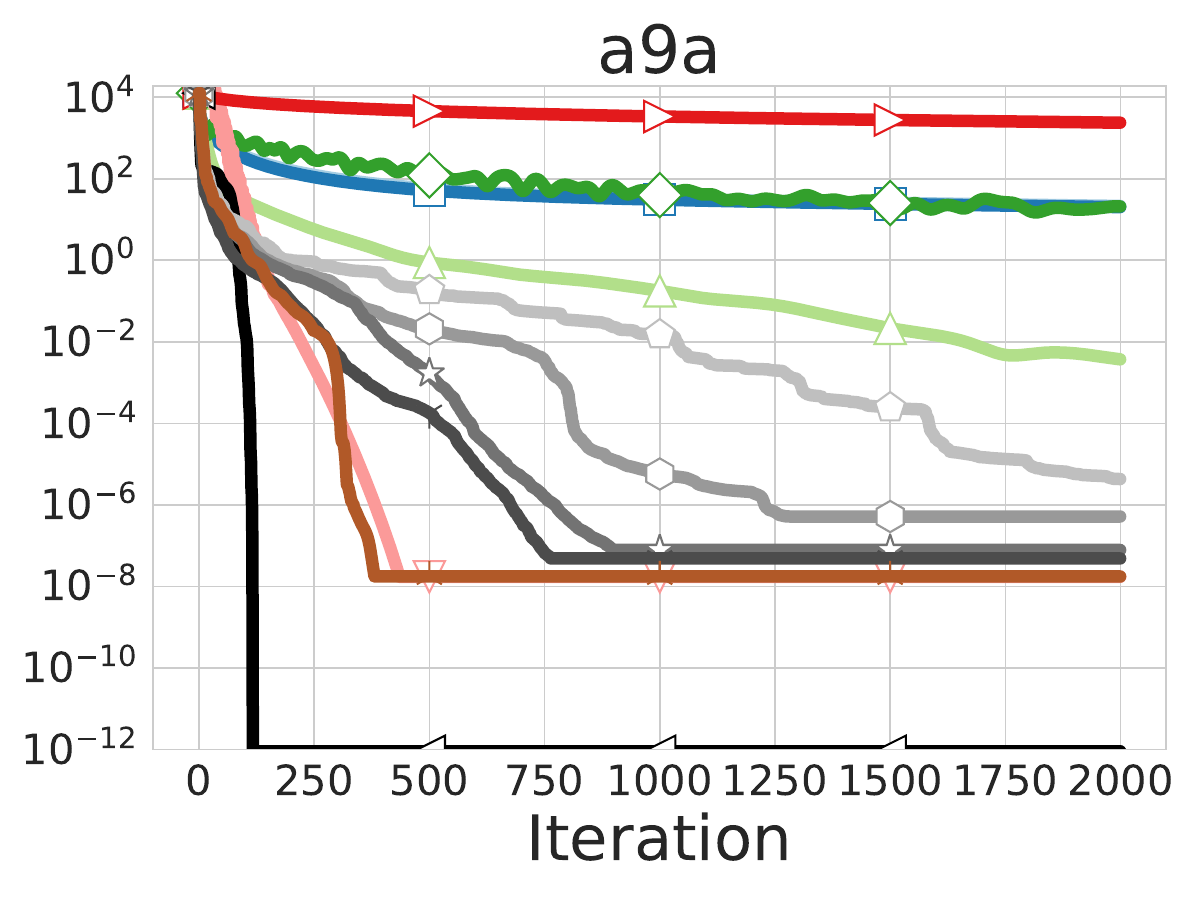}
\includegraphics[scale=0.17]{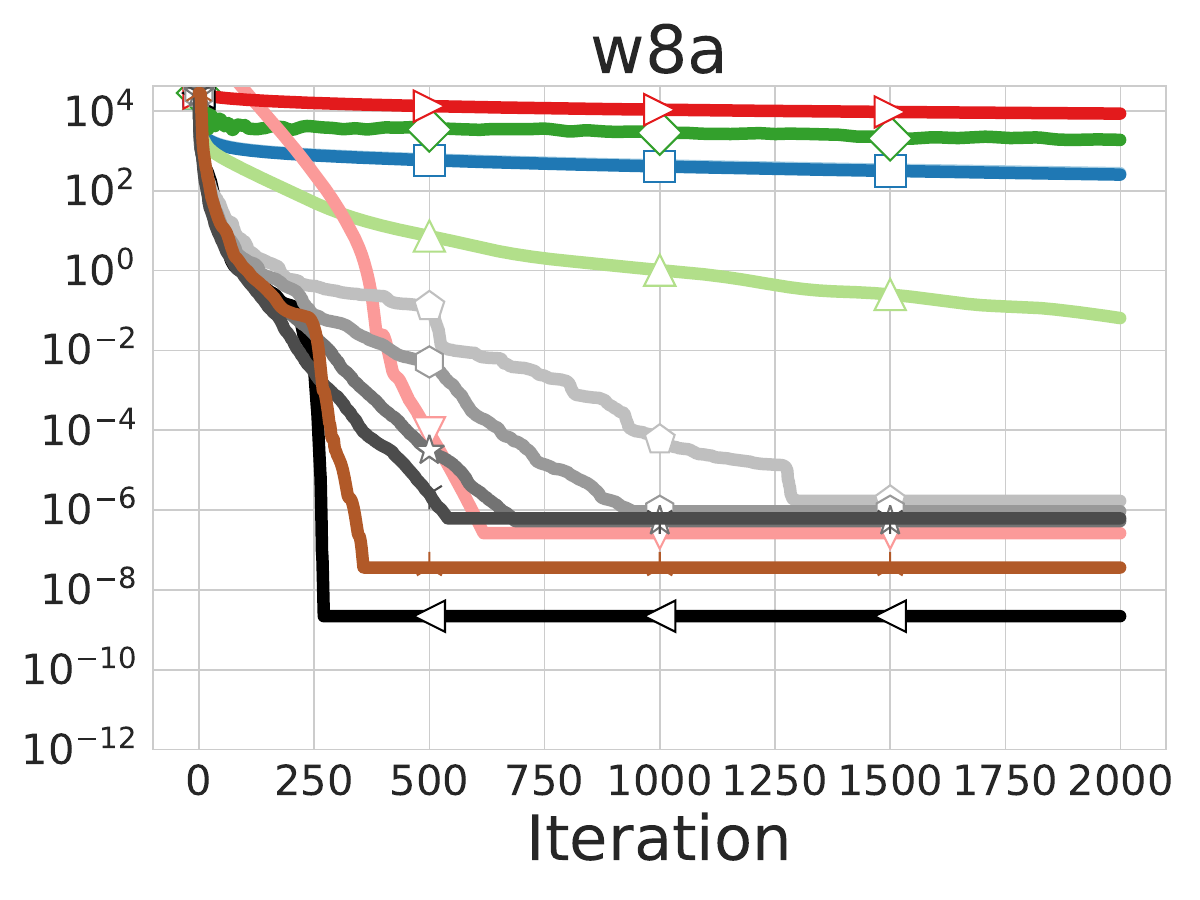}
\includegraphics[scale=0.17]{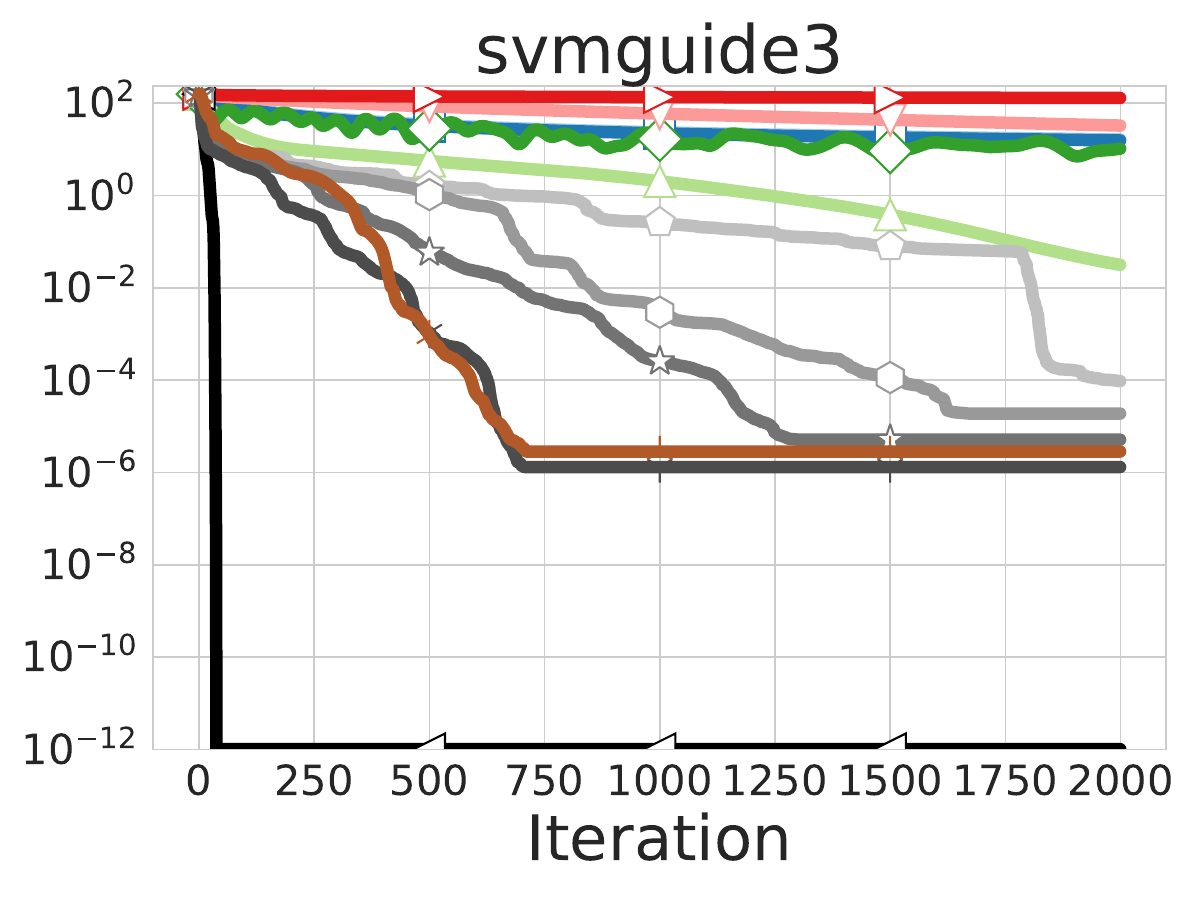}
\includegraphics[scale=0.17]{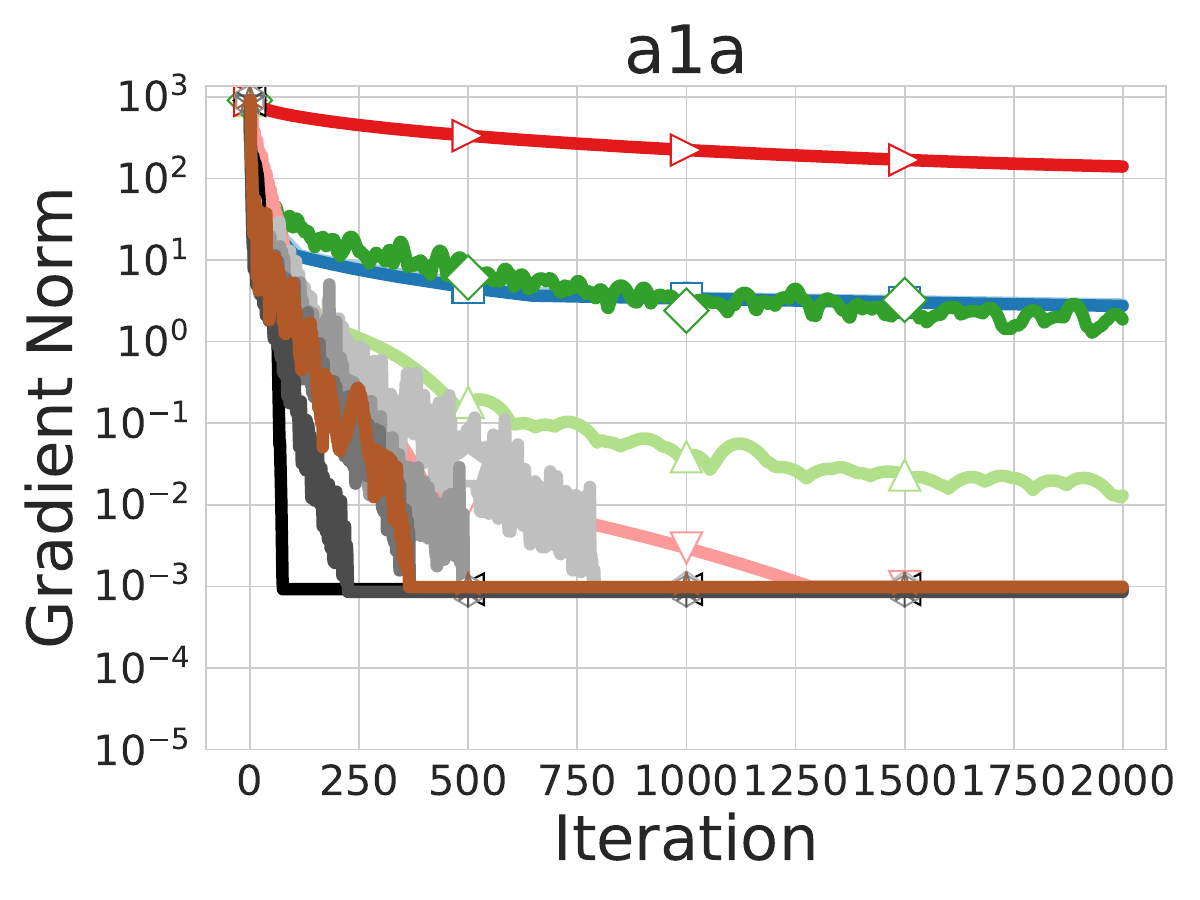}
\includegraphics[scale=0.17]{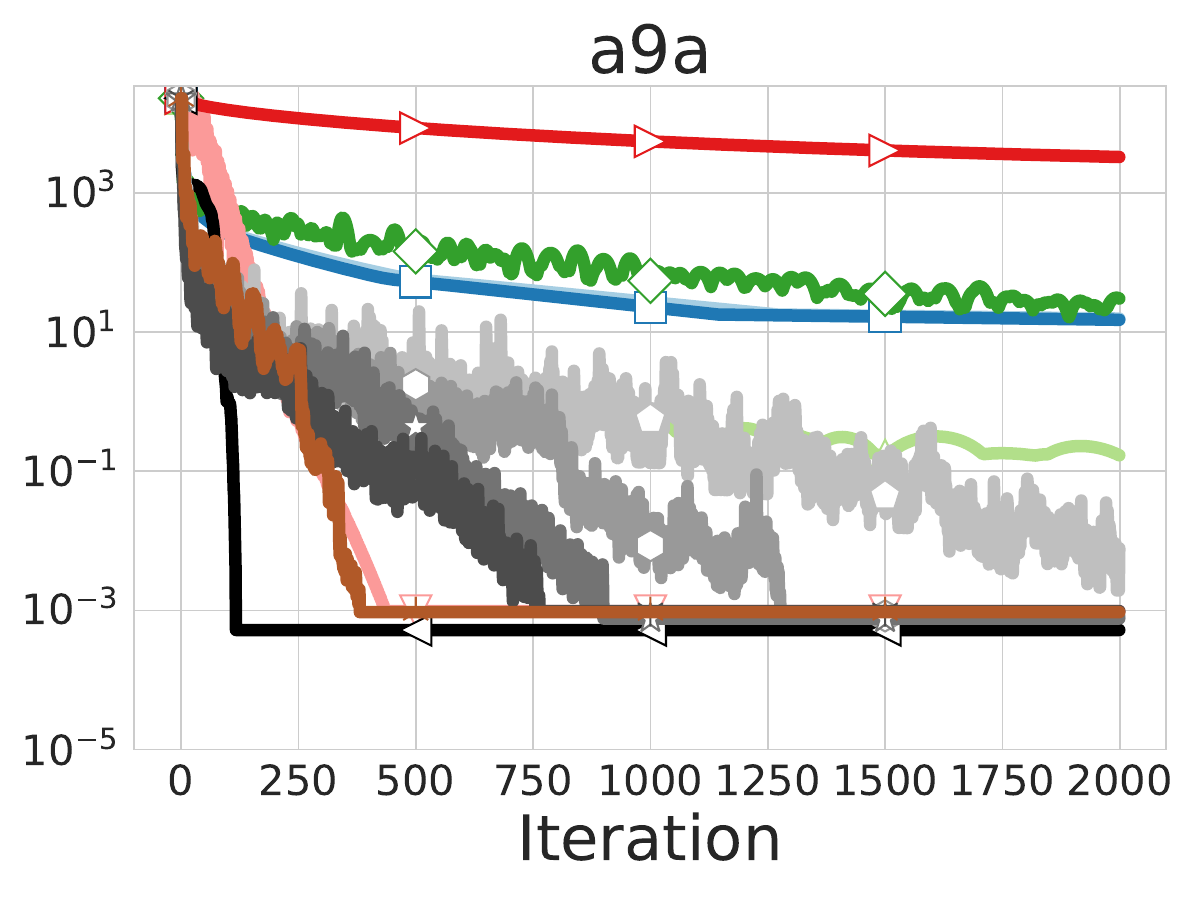}
\includegraphics[scale=0.17]{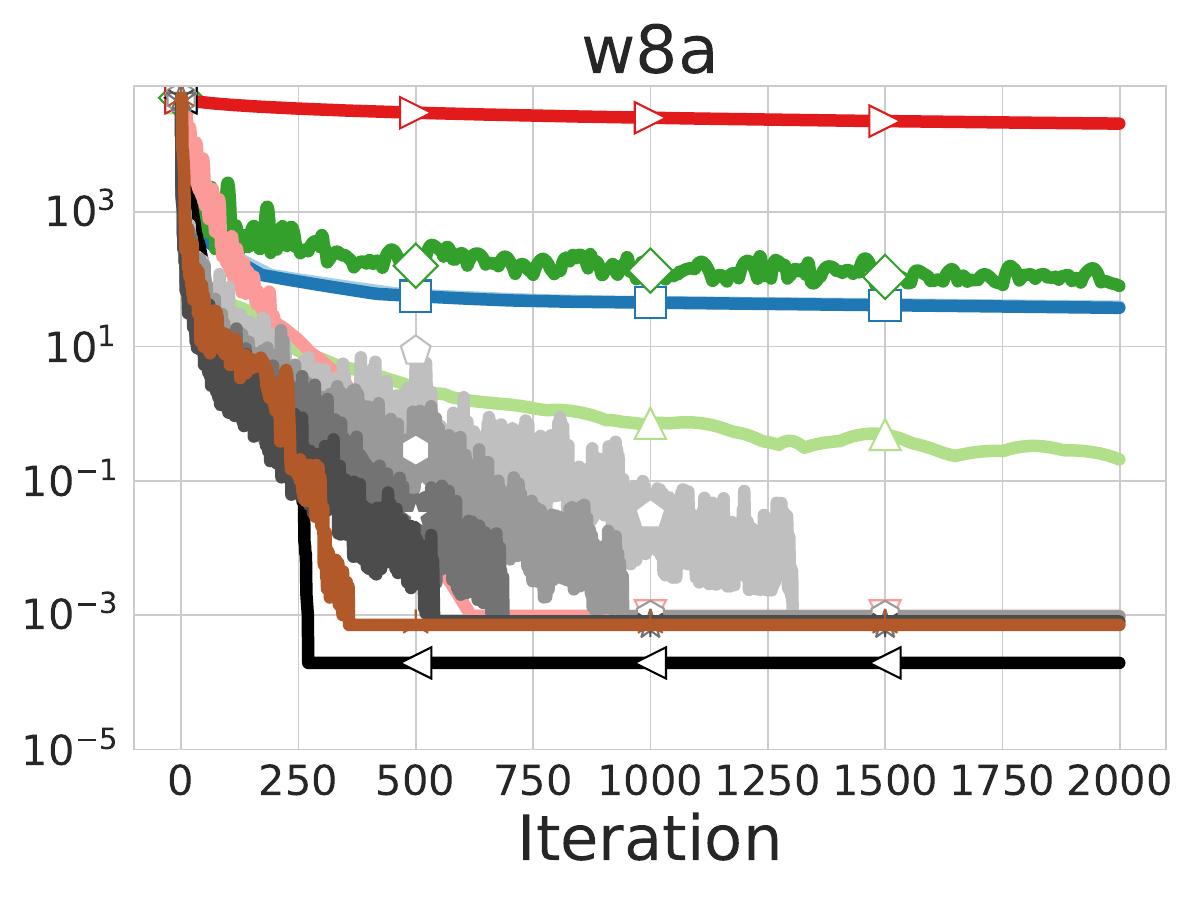}
\includegraphics[scale=0.17]{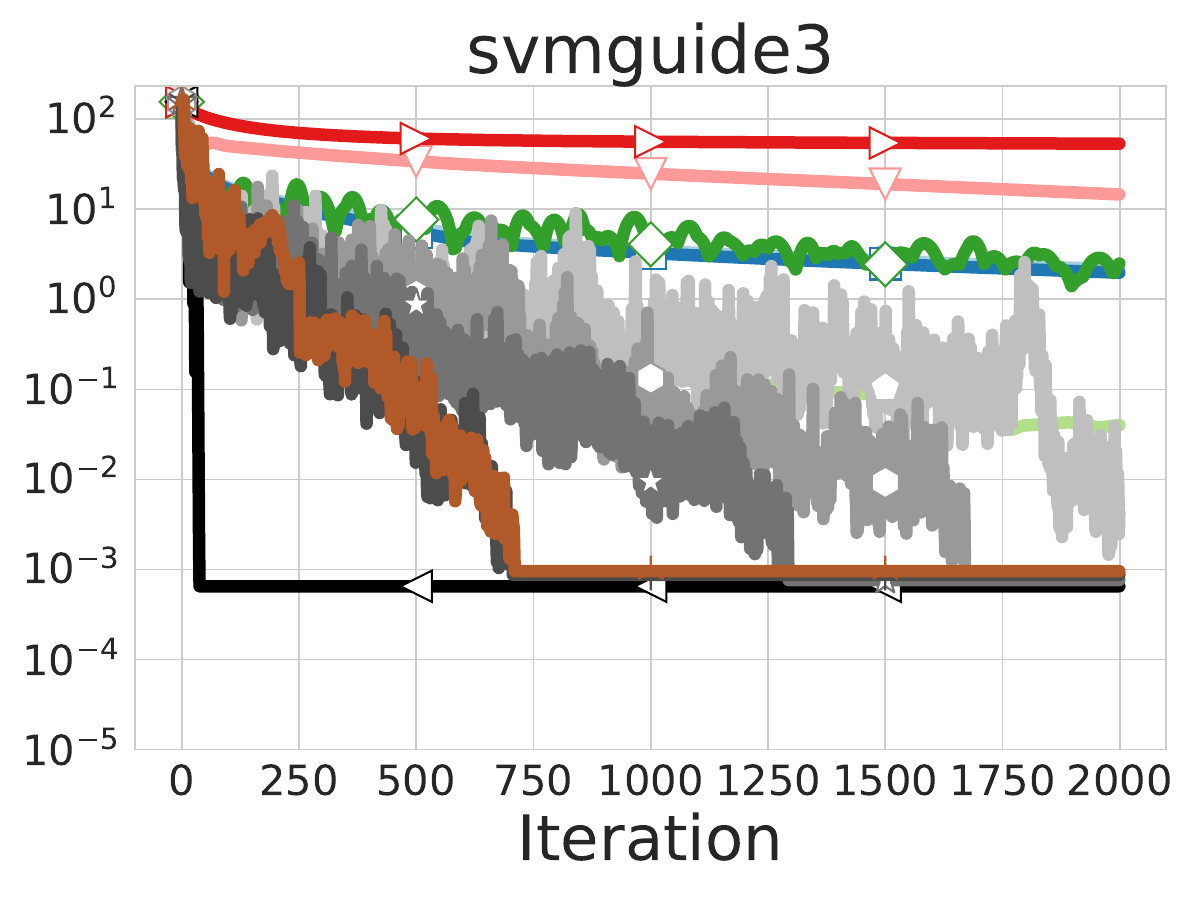}
\includegraphics[scale=0.35]{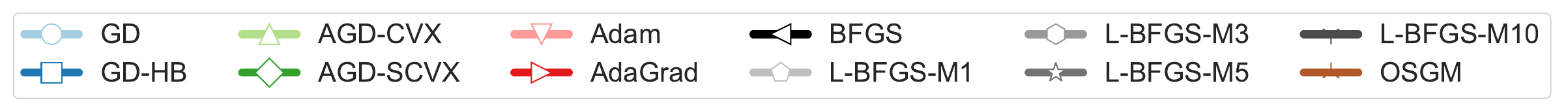}
\caption{Support vector-machine problems. First row: function value gap. Second row: gradient norm. \label{fig:svm}}
\end{figure}

\begin{figure}[!h]
\centering
\includegraphics[scale=0.17]{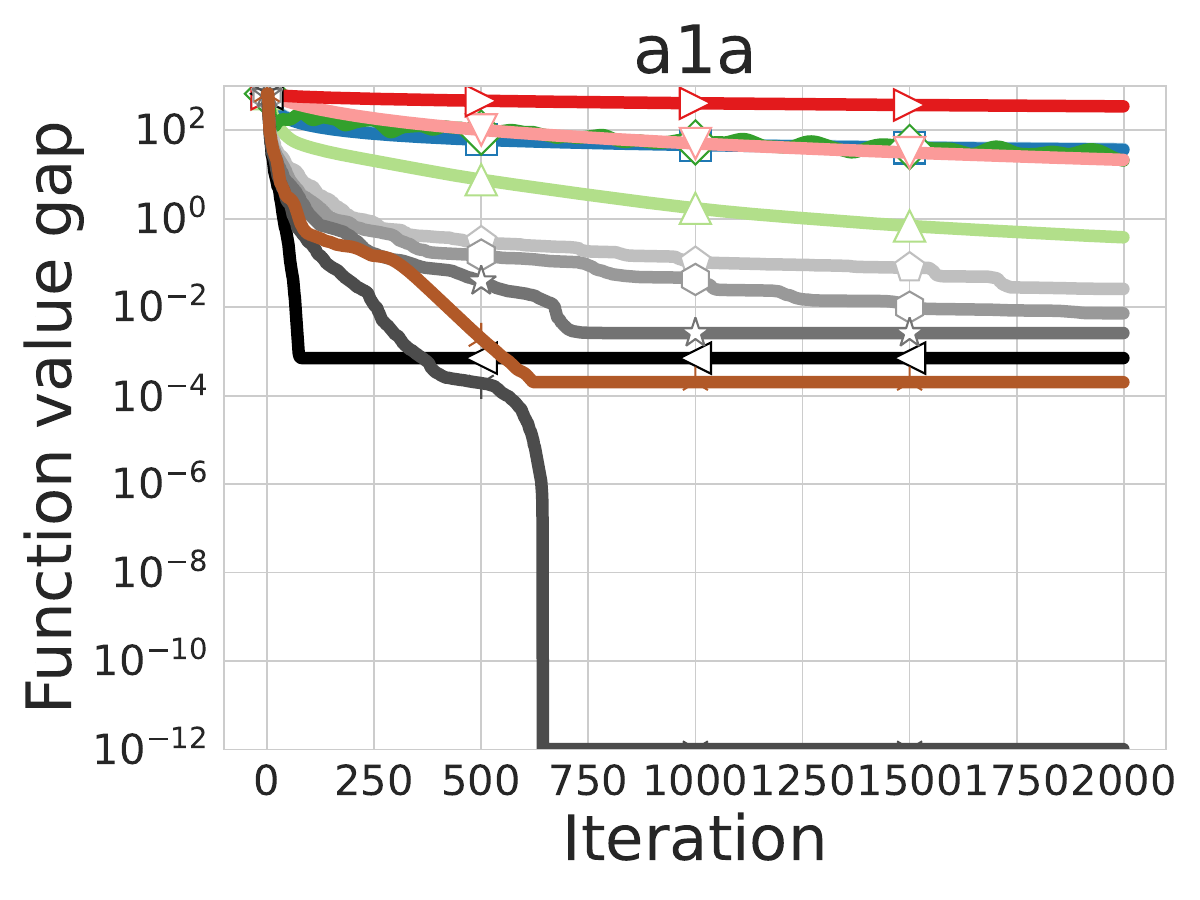}
\includegraphics[scale=0.17]{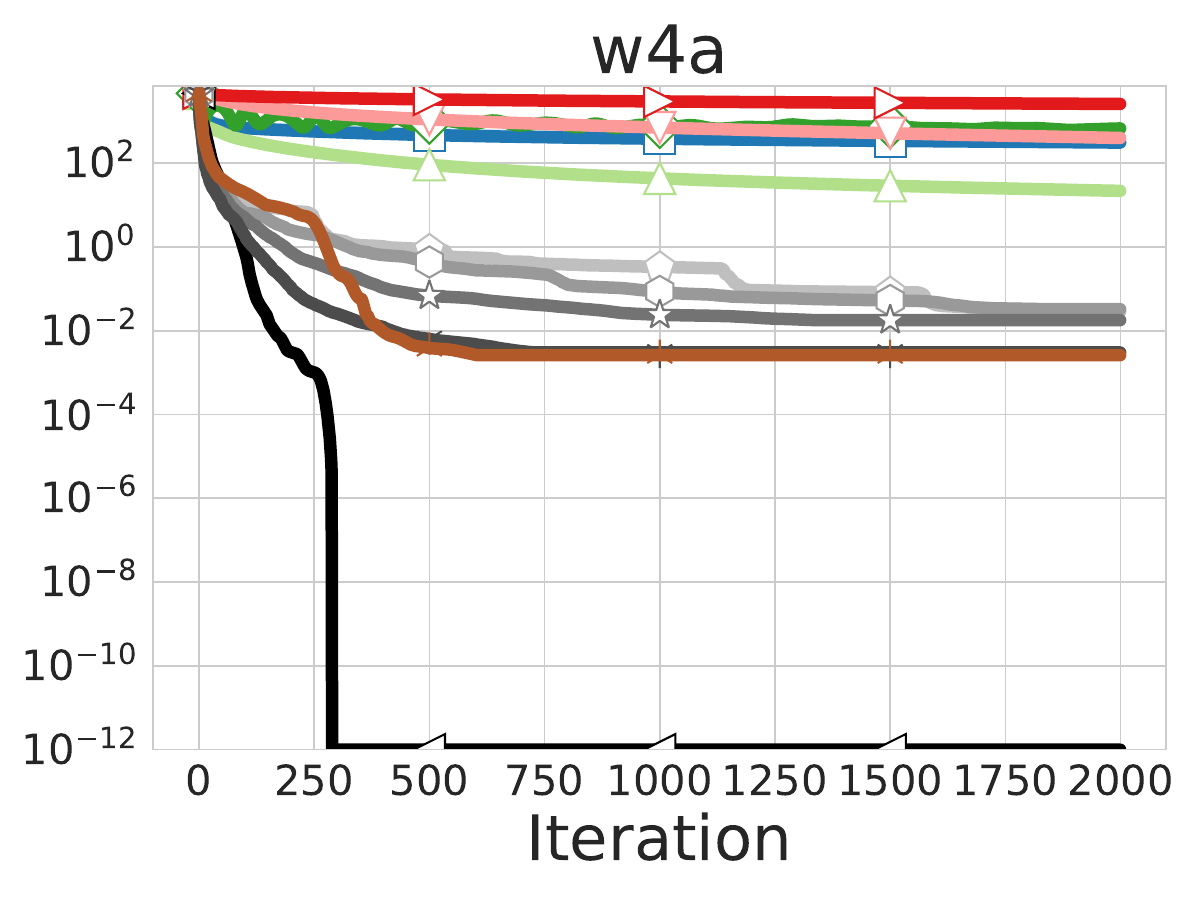}
\includegraphics[scale=0.17]{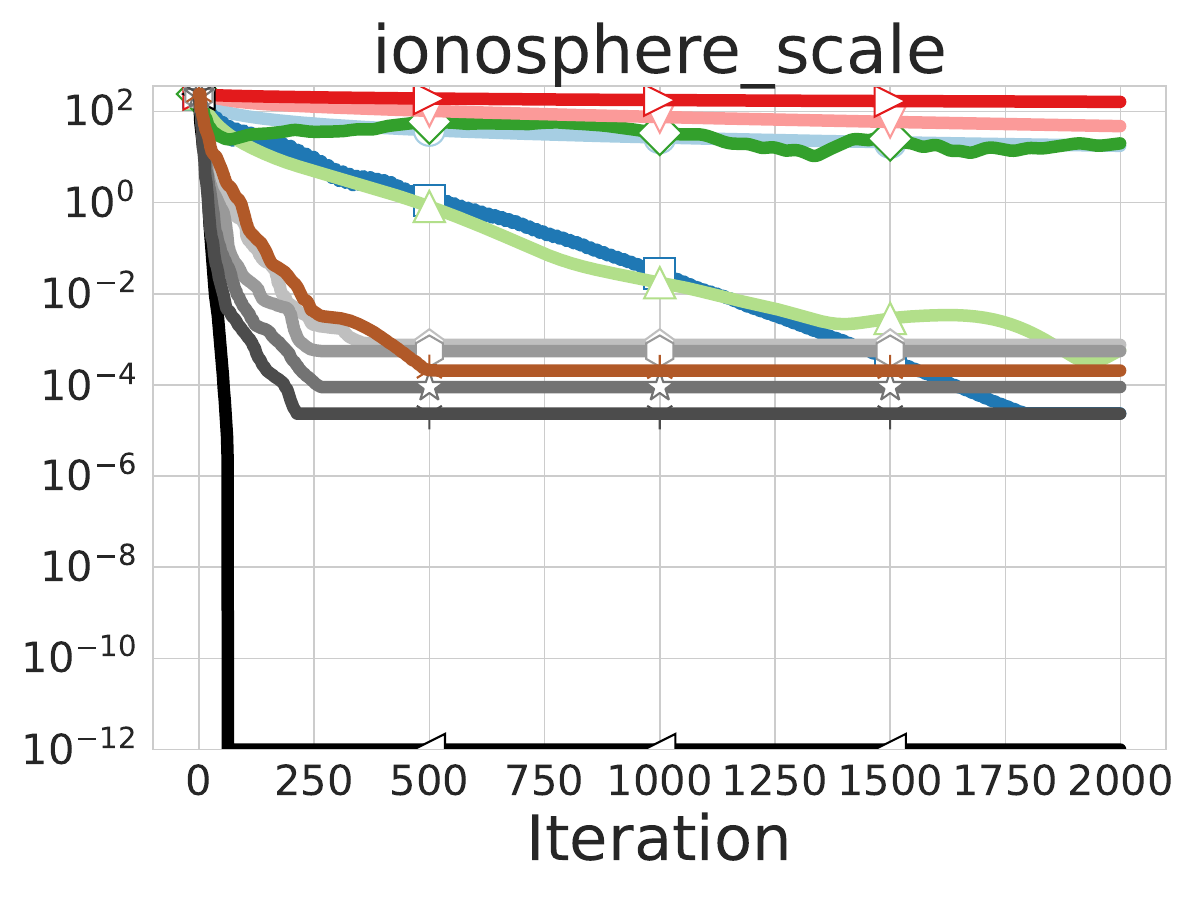}
\includegraphics[scale=0.17]{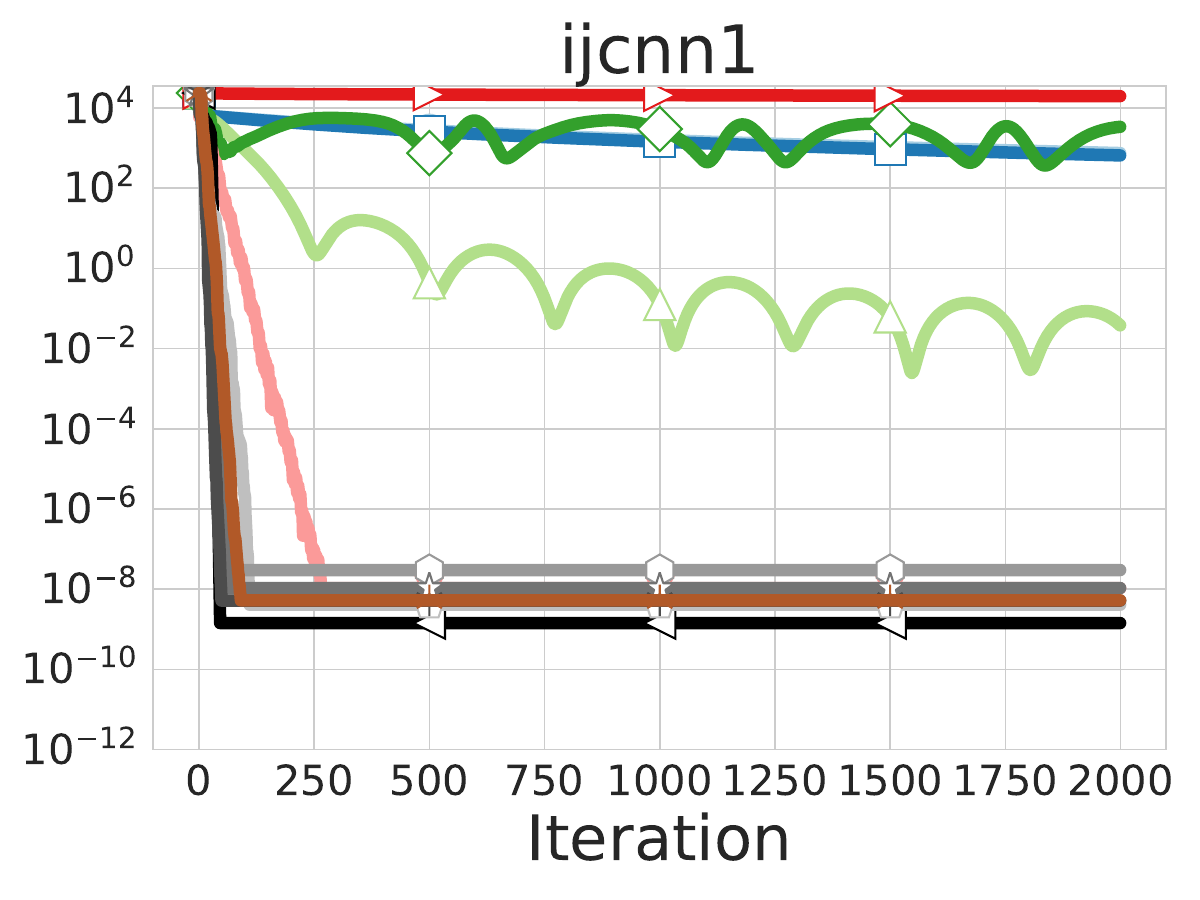}
\includegraphics[scale=0.17]{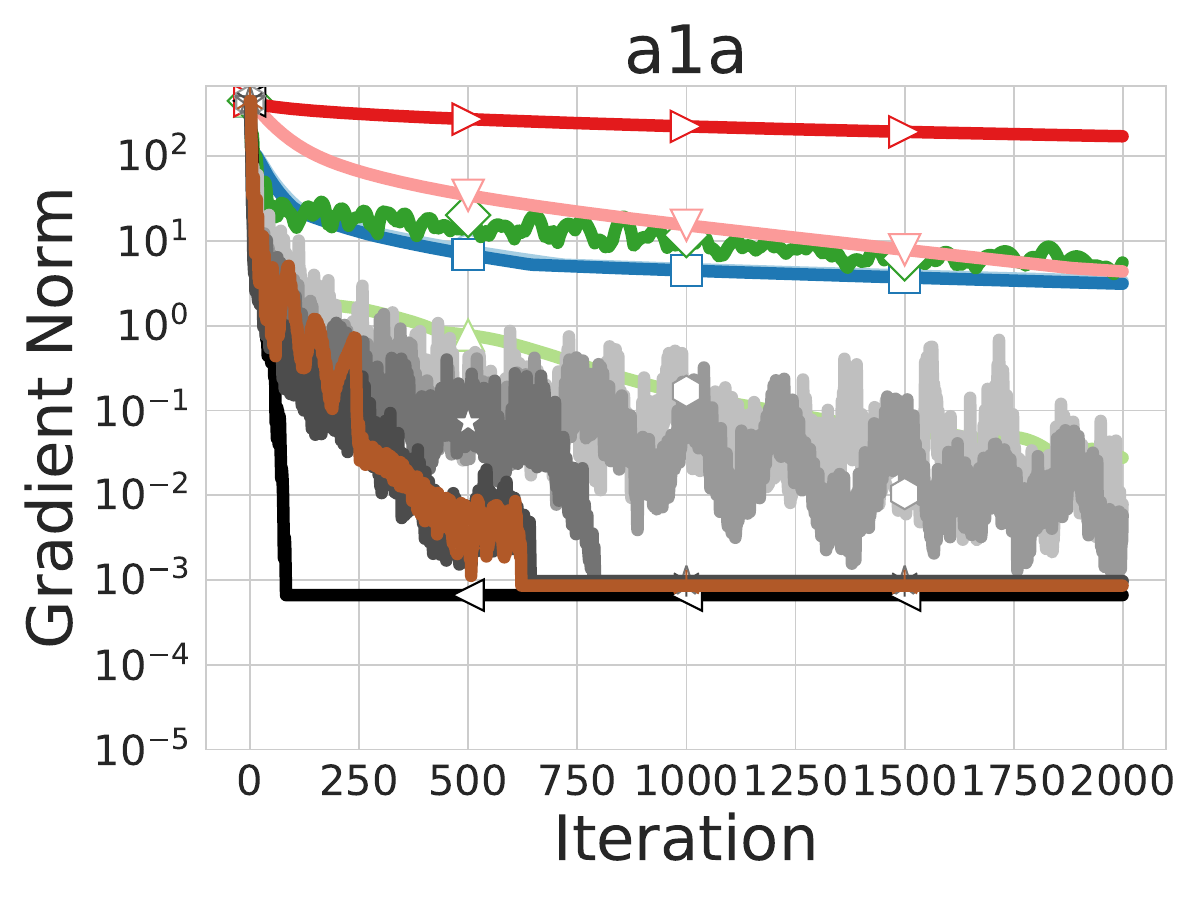}
\includegraphics[scale=0.17]{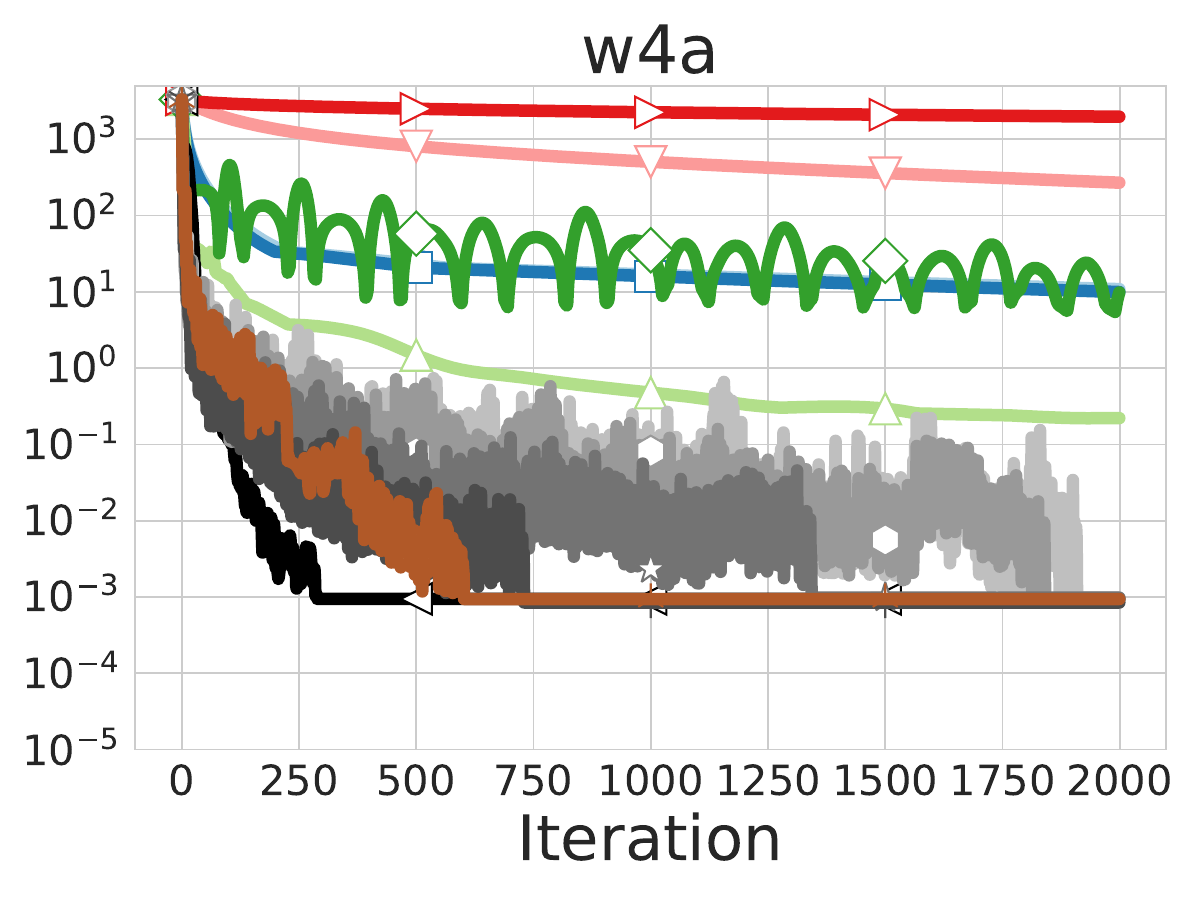}
\includegraphics[scale=0.17]{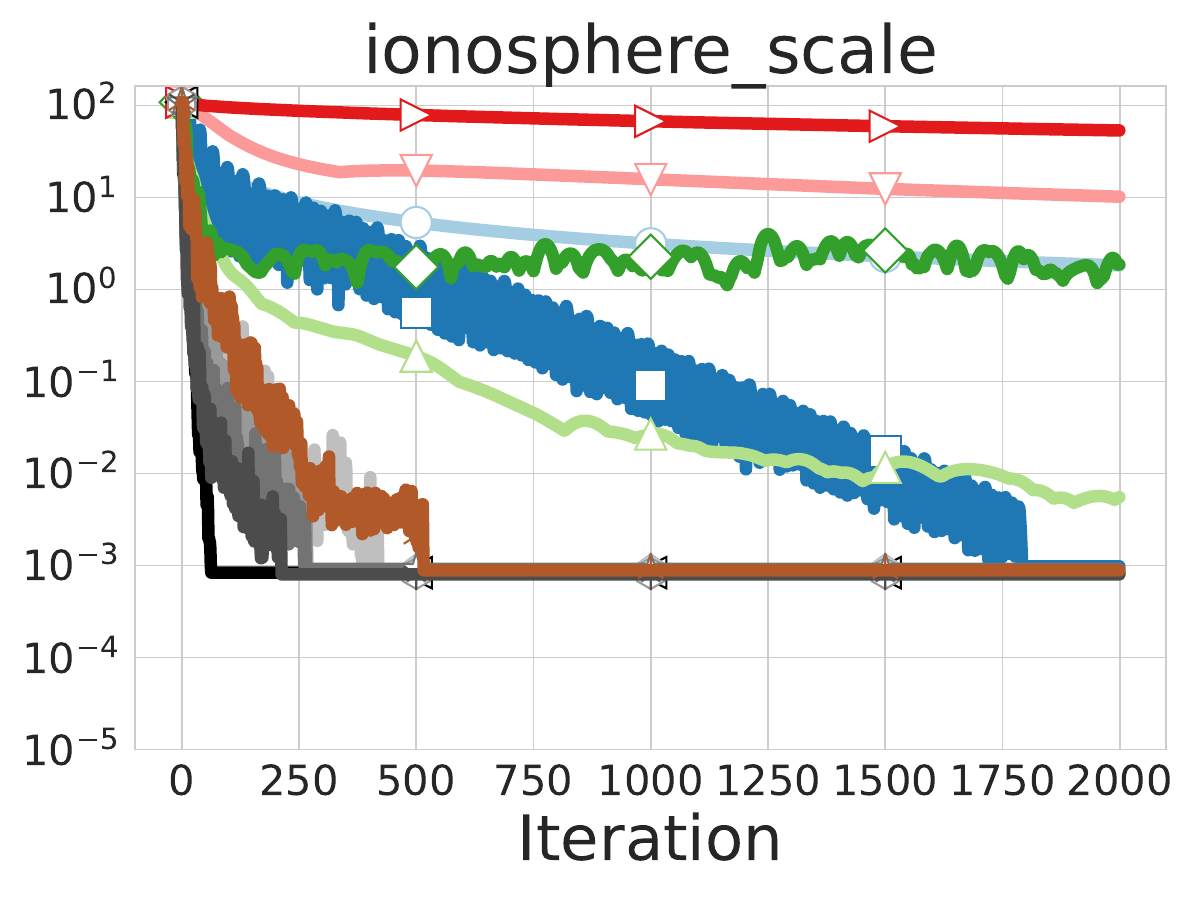}
\includegraphics[scale=0.17]{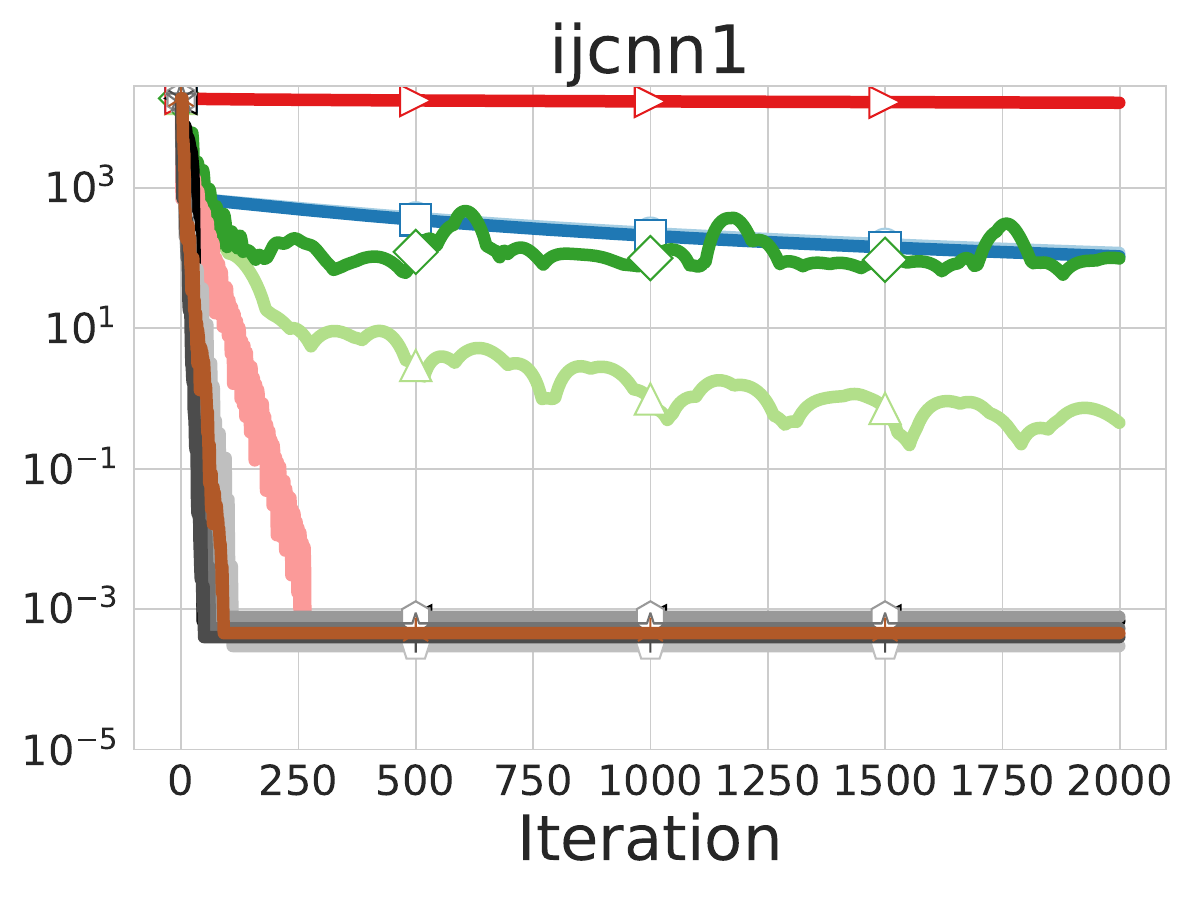}
\caption{Logistic regression problems. First row: function value gap. Second row: gradient norm.  \label{fig:logistic}}
\end{figure}

For each algorithm, we record the number of successfully solved instances ($\|\nabla f\|_\infty \leq 10^{-3}$ within 1000 gradient oracles). \Cref{table:stats} summarizes the detailed statistics. The number of instances solved by {\ob} is comparable to that of \texttt{L-BFGS-M10}.

\paragraph{Support Vector Machine.} \Cref{fig:svm} shows the function value gap and gradient norm plots on sample test instances on support vector machine problems. The optimal value for each instance is obtained by running {\bfgs} until $\|\nabla f \|_\infty \leq 10^{-3}$.  We see that the practical variant of {\ob} achieves a significant speedup over other adaptive first-order methods. In particular, {\ob} often matches \texttt{L-BFGS-M5} and \texttt{L-BFGS-M10}, while its memory usage is closer to \texttt{L-BFGS-M1}. Notably, {\adam} also achieves competitive performance in several instances.

\paragraph{Logistic Regression.} In logistic regression (\Cref{fig:logistic}), {\ob} still compares well with \texttt{L-BFGS-M5} and is significantly faster than other adaptive first-order methods.\\

\subsection{Nonconvex optimization}
On nonconvex problems, {\ob} still significantly outperforms standard adaptive first-order methods, and its performance is comparable to {\lbfgs}.

\begin{figure}[!h]
\centering
\includegraphics[scale=0.17]{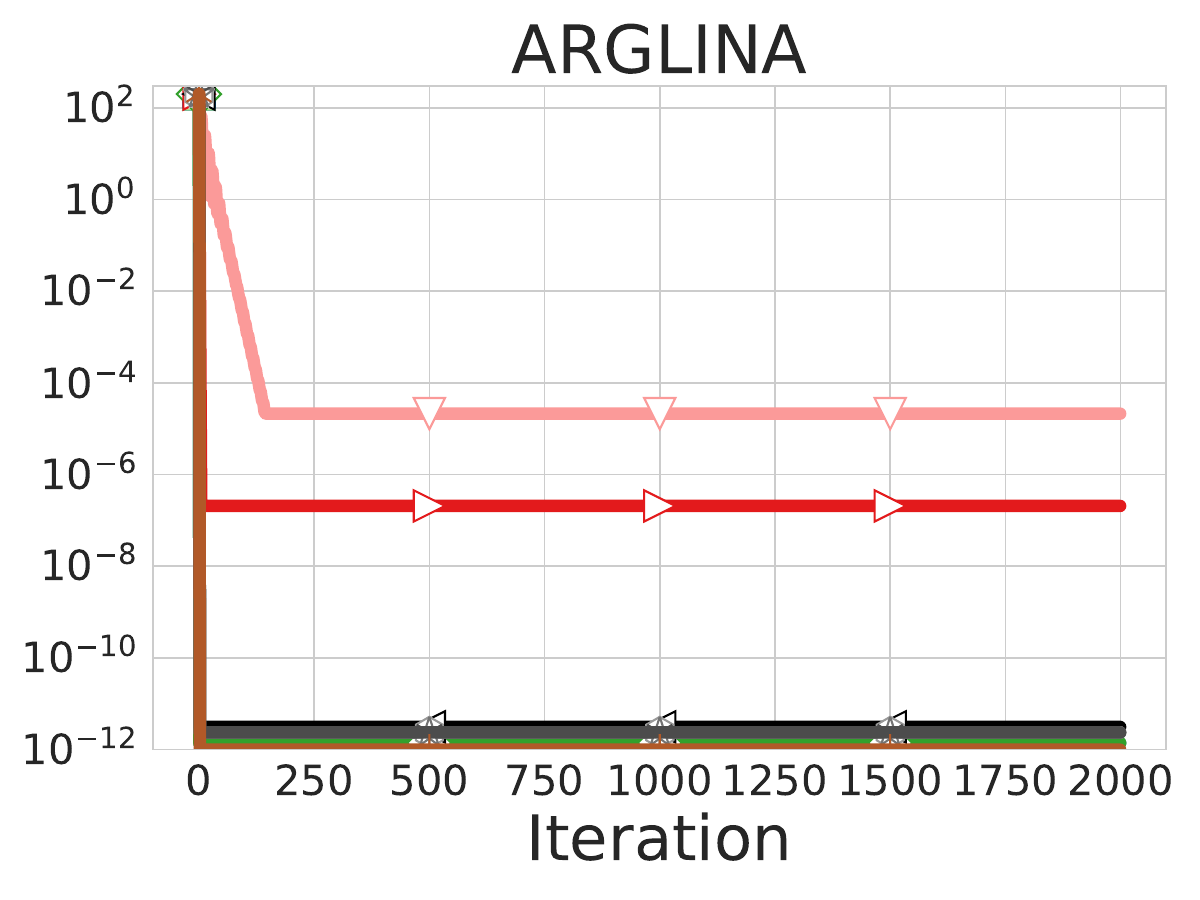}
\includegraphics[scale=0.17]{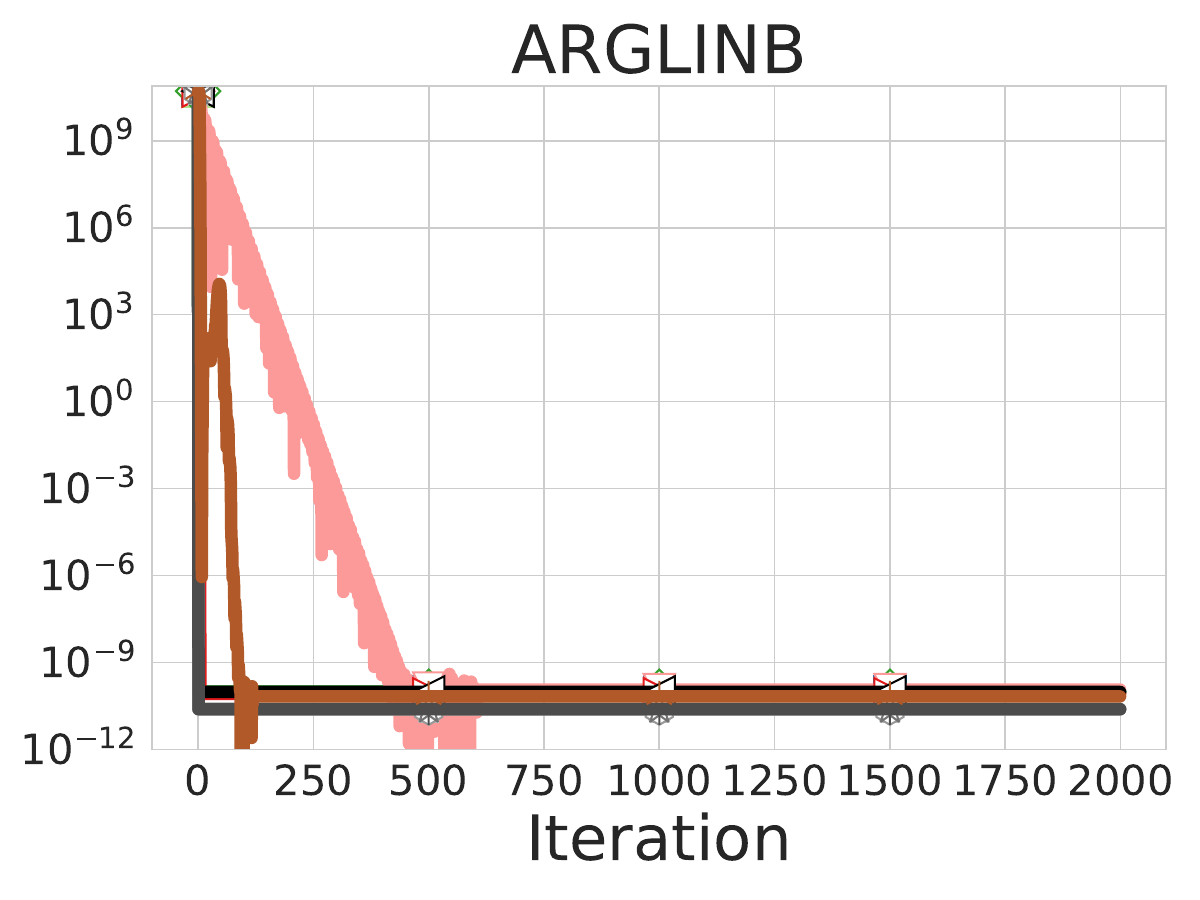}
\includegraphics[scale=0.17]{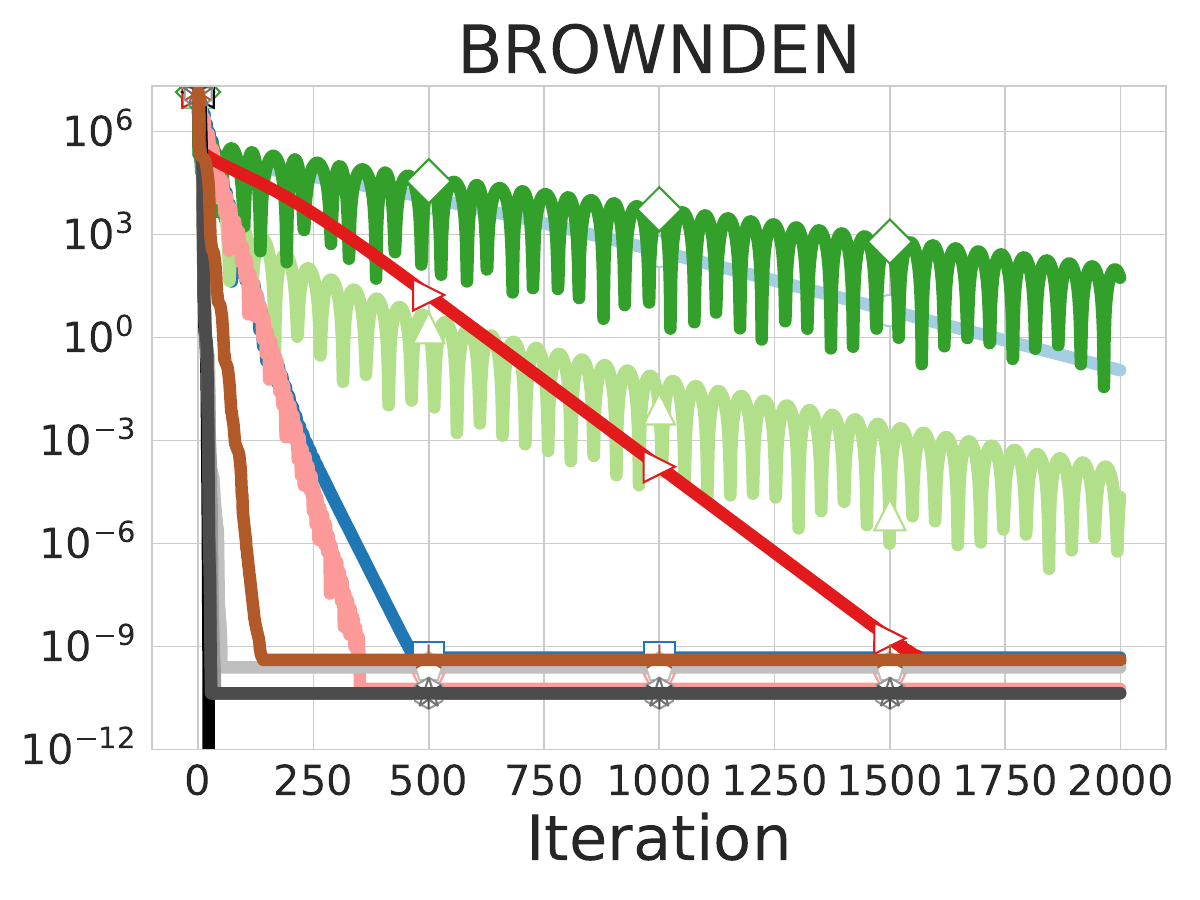}
\includegraphics[scale=0.17]{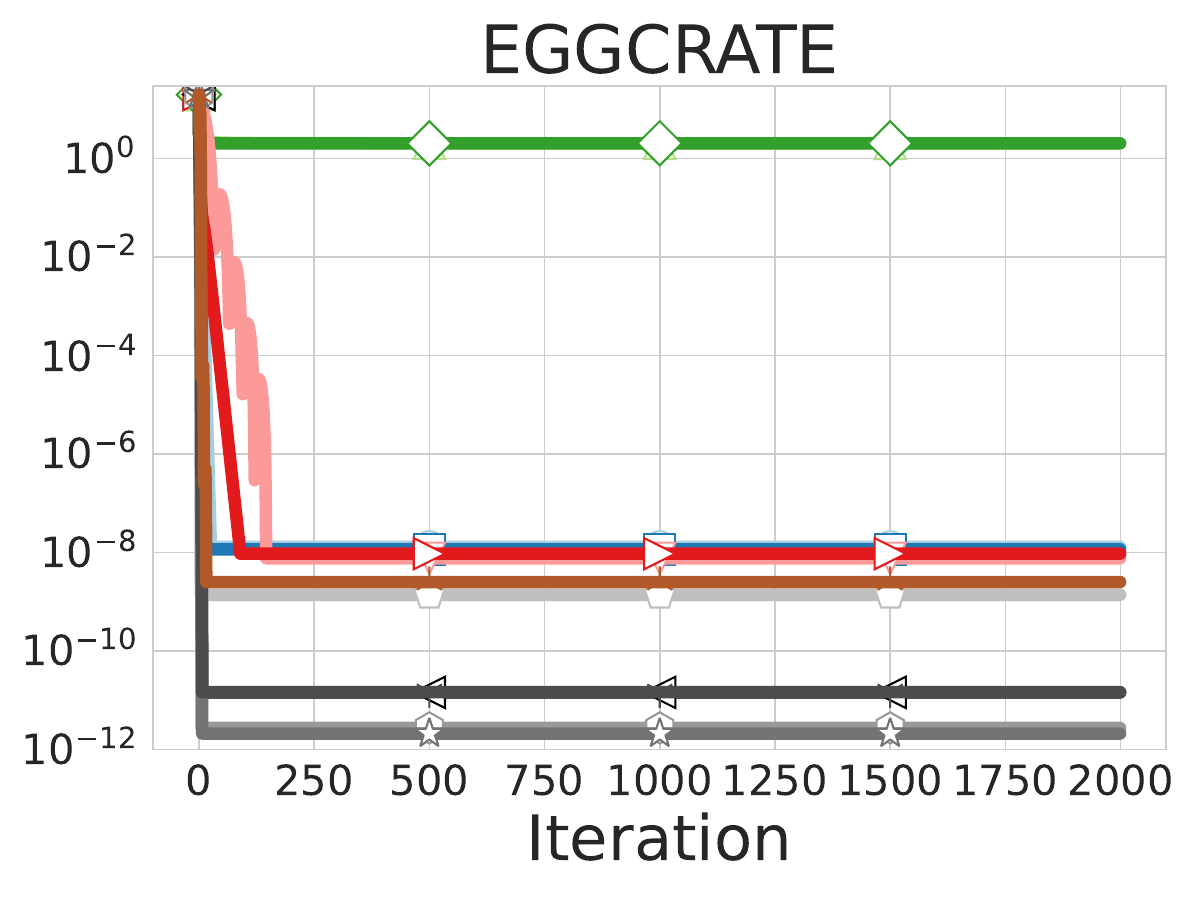}
\includegraphics[scale=0.17]{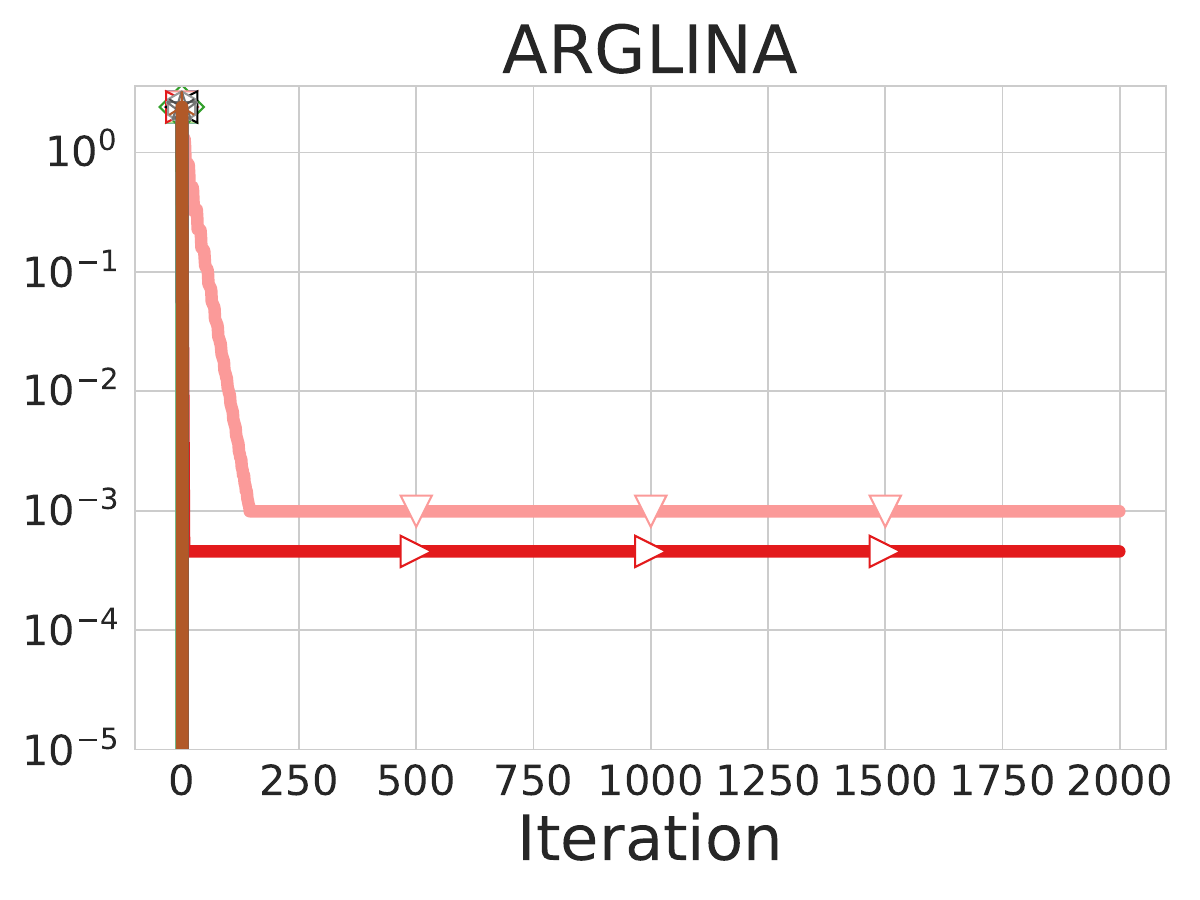}
\includegraphics[scale=0.17]{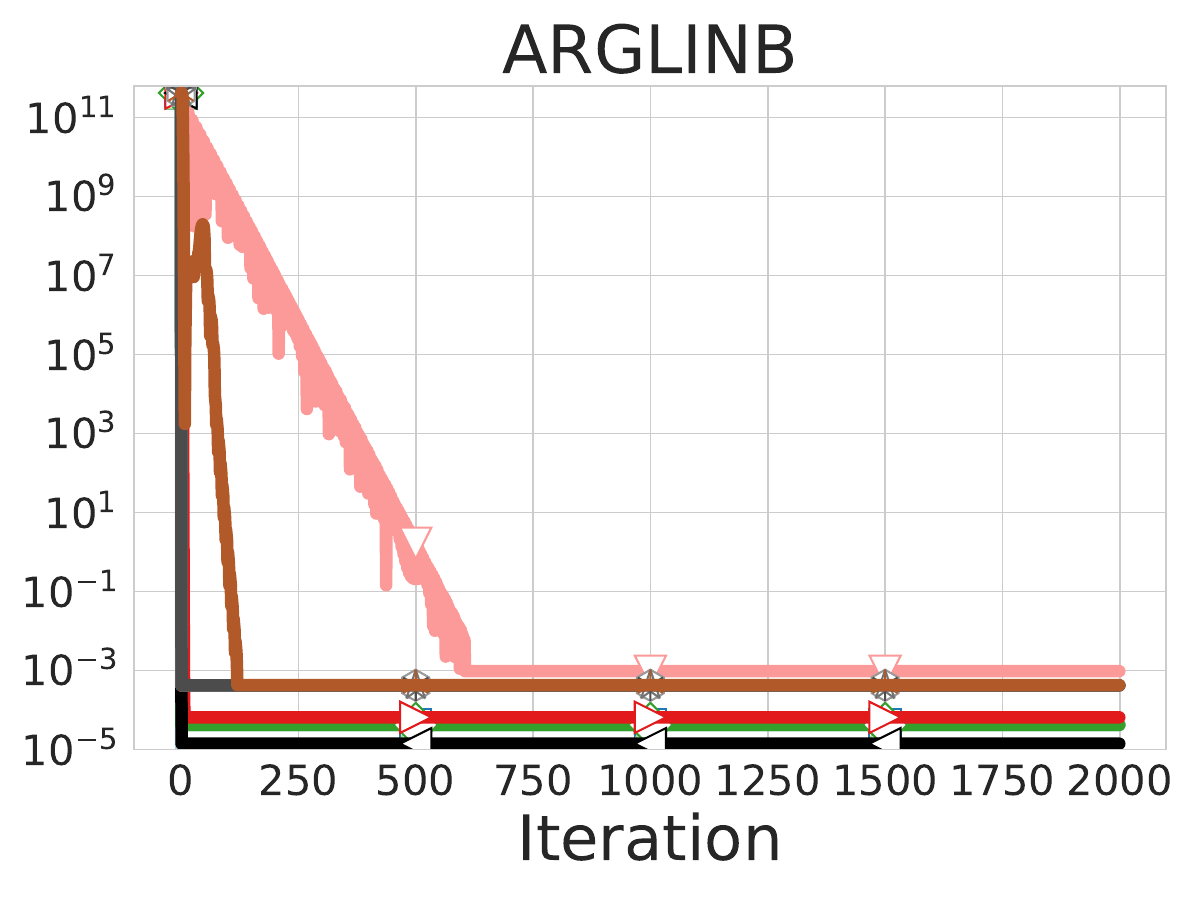}
\includegraphics[scale=0.17]{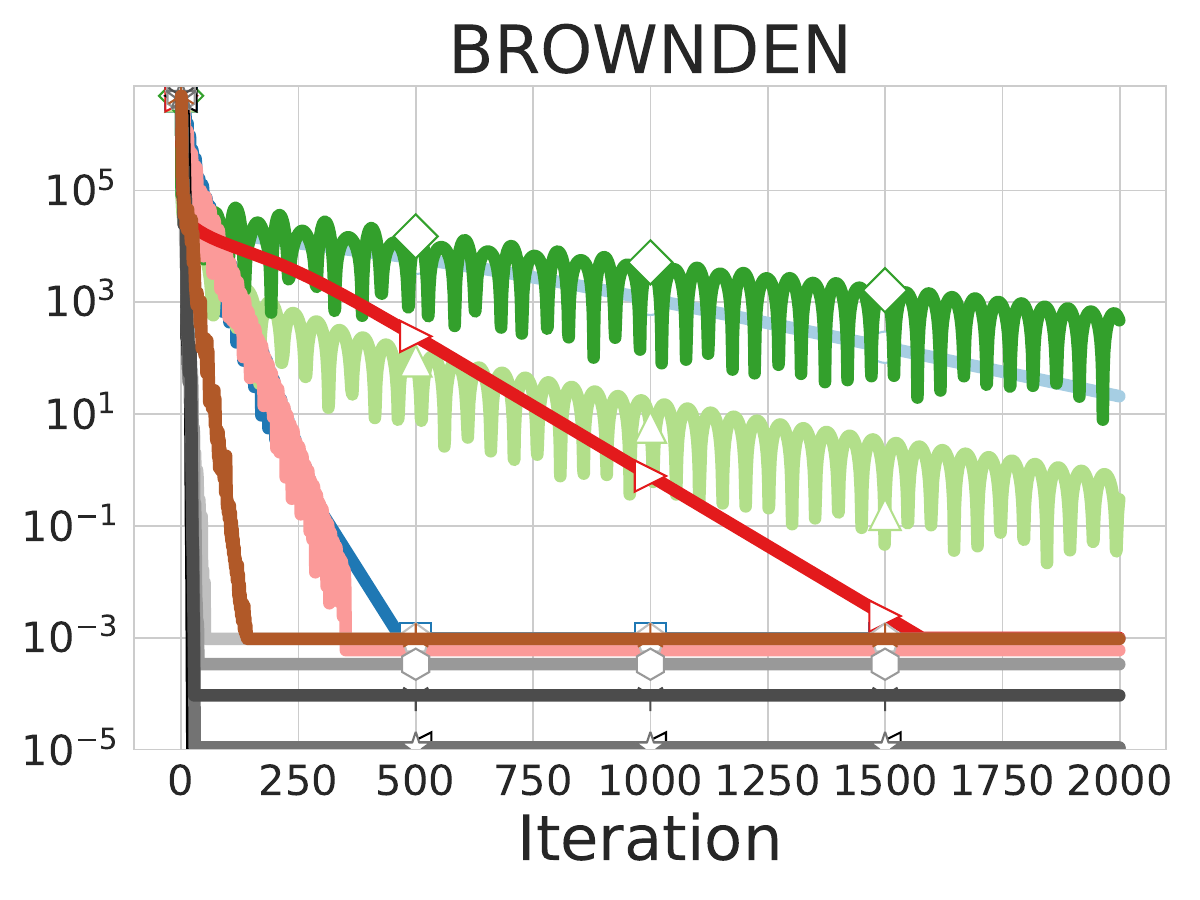}
\includegraphics[scale=0.17]{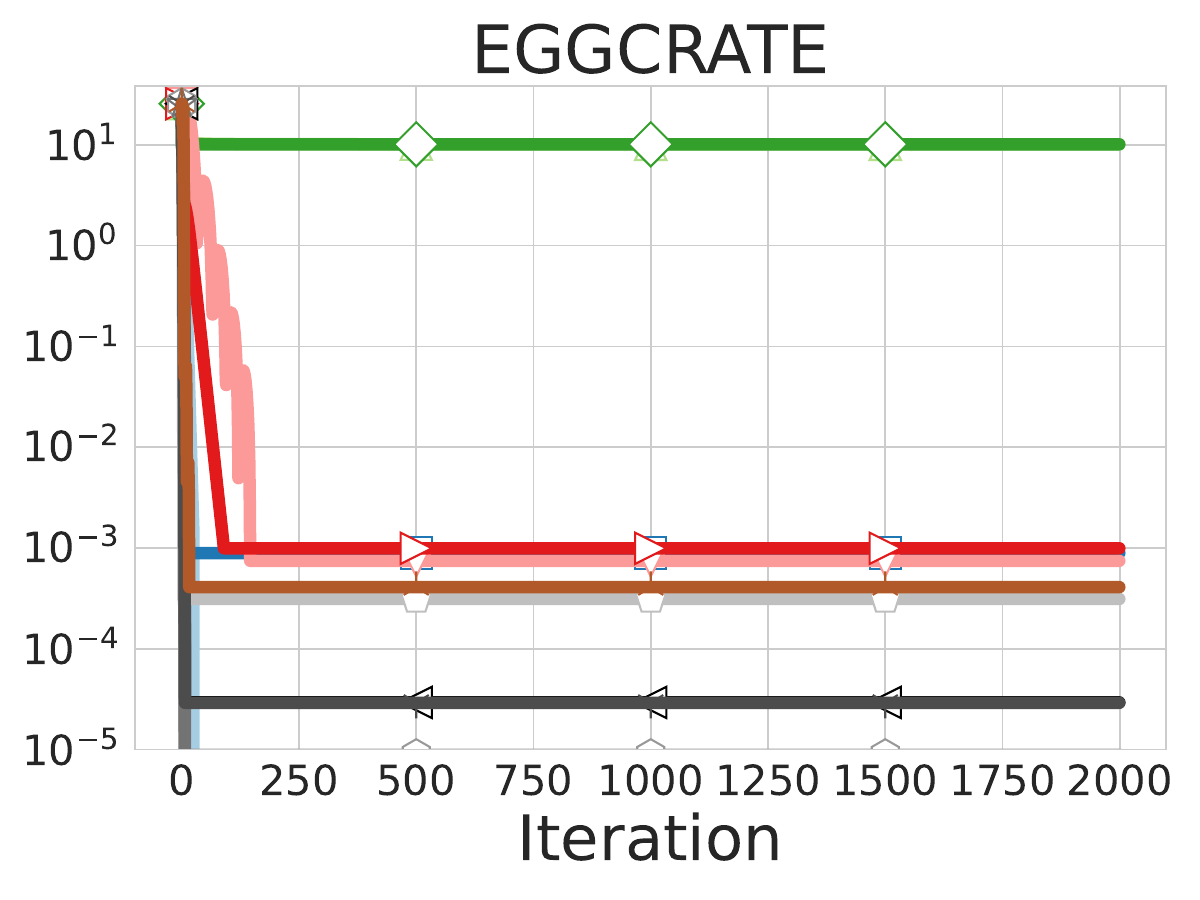}
\includegraphics[scale=0.35]{figs/legend.pdf}
\caption{\texttt{CUTEst} problems. First row: function value gap. Second row: gradient norm. \label{fig:svm}}
\end{figure}

%% file: sec_directions.tex
\section{Future directions and open problems} \label{sec:extensions}

In this section, we sketch extensions of {\os} and provide several possible future directions.

\subsection{Barzilai-Borwein stepsize and online proximal point} \label{sec:bb}

The Barzilai-Borwein stepsize (also known as the BB step) is a widely used practical stepsize selection strategy \cite{barzilai1988two}. There are two types of BB steps:
\begin{equation}
\alpha_k \assign \tfrac{\| x^k - x^{k - 1} \|^2}{\langle \nabla f (x^k) -
  \nabla f (x^{k - 1}), x^k - x^{k - 1} \rangle} \quad \text{and} \quad \alpha'_k \assign \tfrac{\langle \nabla f (x^k) - \nabla f (x^{k - 1}),
  x^k - x^{k - 1} \rangle}{\| \nabla f (x^k) - \nabla f (x^{k - 1}) \|^2}. \nonumber
\end{equation}
In the context of convex quadratic minimization, it is known that $\alpha_k$ is
the steepest descent stepsize {\emph{at the previous iteration}} in terms of
function value \cite[Section 2.1]{zou2022delayed},
\[ \alpha_{k} = \argmin_{\alpha} f (x^{k-1} - \alpha \nabla f (x^{k-1})), \]
and similarly, $\alpha'_k$ is the steepest descent stepsize at the previous iteration in terms of the gradient norm
\[ \alpha_{k}' = \argmin_{\alpha}  \| \nabla f (x^{k-1} - \alpha \nabla f (x^{k-1}))
   \| . \]
Both BB steps fit into {\os} framework by rewriting them as the minimizers of feedback functions:
\begin{equation}
\alpha_{k+1} = \argmin_{\alpha} ~h_{x^k} (\alpha) \assign \tfrac{f (x^k - \alpha \nabla f (x^k)) -f(x^k)}{\| \nabla f(x^k) \|^2} \quad \text{and} \quad \alpha'_{k+1} = \argmin_{\alpha} ~g_{x^k} (\alpha) \assign \tfrac{\| \nabla f (x^k - \alpha \nabla f (x^k)) \|}{\| \nabla f(x^k) \|}, \nonumber
\end{equation}
where $h_x(\alpha)$ is the hypergradient feedback and $g_x(\alpha)$ is the gradient norm feedback discussed in \cite{gao2024gradient}. 
The BB steps can be derived from the minimization of one of the feedback functions, and we may further develop this perspective to create variants of the BB step powered by online learning.
For example, we can write the BB step as
as the online proximal point update 
(also known as implicit online learning \cite{kulis2010implicit})
with $\eta \rightarrow \infty$
\begin{equation} \label{eqn:bb-osgm-proxpt}
	\alpha_{k + 1} = \argmin_{\alpha} ~ \{ h_{x^k} (\alpha) + \tfrac{1}{2
   \eta} (\alpha - \alpha_k)^2 \} .
\end{equation}
To our knowledge, the connection between the BB step and online learning has not been previously explored. Related ideas have primarily appeared in the literature under the guise of ``retarded'' or ``delayed'' gradient updates \cite{friedlander1998gradient,zou2022delayed}. We believe this new perspective can motivate extensions or improvements on the BB step.  
For example, using finite values of $\eta$ above would yield stabilized algorithms similar to the BB step, with a pathway towards a provable convergence guarantee through this connection to the online proximal point method and {\os}.
Alternatively, this perspective clarifies how to develop momentum or diagonal variants of the BB step. In \Cref{app:bb}, we illustrate a convergence analysis of {\oh} based on \eqref{eqn:bb-osgm-proxpt}. 

\begin{ques}
The BB step is guaranteed to converge on quadratics, but no theoretical guarantees are available for general functions without complex linesearch \cite{zhang2004nonmonotone}. Is it possible to use this online learning perspective to design a competitive variant of the BB step with a theoretically guaranteed convergence rate for general functions?
\end{ques}

Finally, for quadratic $f$, the BB step uses the steepest descent stepsize from the previous iteration. In contrast, the steepest descent method knows the future ahead of time and applies the ``optimal'' stepsize. Although the BB step can be unstable when $f$ has a large condition number, it typically outperforms steepest descent. This phenomenon, where knowing the future hinders convergence, aligns with our observation in \Cref{sec:prescient}. To formally state the analogy,  steepest descent is to BB step what {\hdmclassic} is to {\oh}. 

\subsection{{\osr} without knowing $f^{\star}$} \label{sec:nofstar}

Although {\oh} has more appealing practical performance, only
{\osr} theoretically improves the complexity of gradient descent.
However, {\osr} has a drawback of requiring knowledge of $f^{\star}$. 

\begin{ques}
  Is it possible to obtain $\mathcal{O} (\kappa^{\star} \log (1 /
  \varepsilon))$ complexity without knowledge of $f^\star$?
\end{ques}

This issue also appears in well-known methods such as the Polyak
stepsize \cite{hazan2019revisiting}. Two strategies from the Polyak stepsize literature can
be borrowed: a double-loop algorithm that searches for $f^\star$ can maintain a provable convergence rate, but introduces an additional $\log(1/\varepsilon)$ factor in the bound \cite{hazan2019revisiting}. On the other hand, methods that dynamically estimate $f^\star$, such as the variable target or moving target methods \cite{gower2021stochastic,kim1990variable}, typically cannot recover a nonasymptotic convergence rate in the smooth strongly convex case.
Here, we demonstrate an additional (partial) solution to this question in the special case when a dual bound on the primal objective $f$ can be obtained using the Fenchel trick.
Suppose the original optimization problem takes the form
\[ \min_{x \in \mathbb{R}^n} ~f (x) \assign f_1 (x) + f_2 (A x) \]
where $A \in \mathbb{R}^{m \times n}$. Then, using the Fenchel dual trick, we
can write the dual problem as
\[ \max_{y \in \mathbb{R}^m} ~g(y) \assign - f_1^{\ast} (- A^{\top} y) - f_2^{\ast} (y). \]
Suppose $f: \Rbb^n \rightarrow \Rbb$ and $g: \Rbb^m \rightarrow \Rbb$ are both smooth strongly convex, and that first-order oracles $f(x), \nabla f (x), g(y), \nabla g (y)$ can be efficiently evaluated. Examples satisfying these two conditions include ridge regression, where $f (x)
= \tfrac{1}{2} \| A x - b \|^2 + \tfrac{\lambda}{2} \| x \|^2, g (y) = -
\tfrac{1}{2 \lambda} \| A^{\top} y \|^2 - \langle b, y \rangle - \tfrac{1}{2}
\| y \|^2$ and certain regularized generalized linear models \cite{shalev2013stochastic}. Then strong duality holds  $\min_x f (x) = f^{\star} = \max_y g (y)$, and we define the duality gap function $F(z) = F(x, y) \coloneqq f(x) - g(y)$. Since $F(z^\star) = 0$, we can apply {\os} to the duality gap function $F$ to achieve $\mathcal{O} (\kappa^{\star}_F \log (1 / \varepsilon))$ complexity, where $\kappa^{\star}_F$ is the optimal condition number of the gap function. This idea can also be applied to remove $f^\star$ in the Polyak stepsize.

\subsection{{\os} for accelerated gradient descent}
{\os} works well with heavy-ball momentum. Can it similarly improve the performance of Nesterov momentum, for example, by using {\agd} as our base algorithm instead of gradient descent? 

\begin{ques}
Is there a method that (with knowledge of $f^\star$) achieves (asymptotic) $\Ocal(\sqrt{\kappa^\star} \log(1/\varepsilon))$ complexity with only access to a first-order oracle for $f$?
\end{ques}

To be concrete, consider the {\agd} and preconditioned {\agd} iterations (\Cref{alg:agd}, \Cref{alg:agdpre}),
\begin{figure}[h]
\begin{minipage}[t]{0.49\textwidth}
\begin{algorithm}[H]
\caption{Accelerated GD \label{alg:agd}}
{\textbf{input} initial point $x^1=z^1$}\\
\For{$k = 1, 2, \dots$}{
$y^k = x^k + \tfrac{1}{\sqrt{\kappa} + 1} (z^k - x^k)$\\
$x^{k + 1} = y^k - \tfrac{1}{L} \nabla f (y^k)$ \label{line-agd-update} \\
$z^{k + 1} = ( 1 - \tfrac{1}{\sqrt{\kappa}} ) z^k +
  \tfrac{1}{\sqrt{\kappa}} ( y^k - \tfrac{1}{\mu} \nabla f (y^k) )$ 
}
\end{algorithm}
\end{minipage}
\hfill
\begin{minipage}[t]{0.49\textwidth}
\begin{algorithm}[H]
\caption{Preconditioned {\agd}  \label{alg:agdpre}}
{\textbf{input} initial point $x^1=z^1$}\\
\For{$k = 1, 2, \dots$}{
$y^k = x^k + \tfrac{1}{\sqrt{\kappa} + 1} (z^k - x^k)$\\
$x^{k + 1} = y^k - {P} \nabla f (y^k)$\\
$z^{k + 1} = ( 1 - \tfrac{1}{\sqrt{\kappa}} ) z^k +
  \tfrac{1}{\sqrt{\kappa}} ( y^k - \kappa {P} \nabla f (y^k) )$ }
\end{algorithm}
\end{minipage}
\end{figure}
where the standard $x$ update in \Cref{line-agd-update} is replaced by $x^{k + 1} = y^k - P_k \nabla f (y^k)$. We believe that the biggest challenge in combining {\os} with {\agd} lies in feedback
design: as in the heavy-ball method, {\agd} does not guarantee monotonicity in $f (x)$. To faithfully reflect the quality of a stepsize $P$, we likely require another potential function to measure the quality of a stepsize $P$. Unfortunately, potential functions from the literature are all impossible to evaluate efficiently. In fact, they all require knowing the optimum $x^\star$ in advance \cite{d2021acceleration}!
For example, $\phi (x, z) \assign f (x) - f^\star + \tfrac{\mu}{2} \| z - x^{\star} \|^2$ is a valid potential function \cite{d2021acceleration}:
\[ [ f (x^{k + 1}) - f^\star + \tfrac{\mu}{2} \| z^{k + 1} -
   x^{\star} \|^2 ] \leq ( 1 - \tfrac{1}{\sqrt{\kappa}} )
   [ f (x^k) - f^\star + \tfrac{\mu}{2} \| z^k - x^{\star} \|^2
   ] . \]
In view of this problem, the challenge is to design a potential function for {\agd} that does not depend on $x^\star$. For quadratic functions, such a potential function exists (\Cref{thm:agd-quadpot}).
\begin{thm}[Potential function for {\agd} on quadratic] \label{thm:agd-quadpot}
  Suppose $f$ is a strongly-convex quadratic, and that the iterates $\{(x^k,z^k)\}$ are generated by {\agd} with a fixed value of the momentum $\beta$.  Define the potential
  \[ \varphi_{\mu} (x, z) \assign f (z) - f^\star + \tfrac{1}{2 \mu} \|
     \nabla f (x) \|^2 . \]
  Then
  \[ \varphi_{\mu} (x^{k + 1}, z^{k + 1}) \leq \Big[ \tfrac{(\kappa -
     1)^2}{\kappa ( \sqrt{\kappa} + 1 )^2} ( 1 + \tfrac{1}{2
     \kappa} + \tfrac{\sqrt{4 \kappa + 1}}{2 \kappa} ) \Big]~
     \varphi_{\mu} (x^k, z^k) 
     \leq (1 - \tfrac{1}{\sqrt{\kappa}}) \varphi_{\mu} (x^k, z^k). \]
\end{thm}
This potential still requires knowledge of $f^\star$, but as we have seen in \Cref{sec:nofstar}, it is often possible to relax this requirement as well. Then following the same techniques as {\osgmrx} \cite{gao2025gradient}, we can define ratio feedback $r_{x, z} (P) \assign \tfrac{\varphi_\mu(x^+(P), z^+(P))}{\varphi_\mu(x, z)}$, where $x^+(P)$ and $z^+(P)$ are the iterates after one iteration of \Cref{alg:agdpre}. Define $\hat{\kappa}$ such that
\[1 - \tfrac{1}{\sqrt{\hat{\kappa}}} = \min_P \max_{x, z} ~r_{x, z}(P),\]
where $\hat{\kappa} \leq \kappa$ since $r_{x, z}(\frac{1}{L}) \leq 1-\frac{1}{\sqrt{\kappa}}$ for any $x, z$. Then an $\Ocal(\sqrt{\hat{\kappa}} \log(1/\varepsilon))$ complexity can be obtained. However, whether it is possible to design an algorithm that achieves $\Ocal(\sqrt{\kappa^\star} \log(1/\varepsilon))$ with only blackbox access to $f$ remains unclear.

\paragraph{Feedback design through performance estimation.} 
{\os} provides a general mechanism to accelerate iterative methods when a computable potential (Lyapunov) function is known. Although such potential functions can often be obtained from literature (for example, \Cref{sec:prac} demonstrates this approach in the case of heavy-ball momentum), it can be challenging to find a computable potential function that does not rely on $x^\star$. 

\begin{ques}
Is there a systematic way to design computable potential functions for iterative methods (that do not rely on the exact solution $x^\star$)?
\end{ques}

Recent advances in first-order methods have shown that finding a potential function can be reduced to a small-scale semidefinite optimization problem through performance estimation (PEP) \cite{taylor2018lyapunov,upadhyaya2025automated}. The advantage of the PEP framework is that it can numerically find a tight worst-case potential function for some problem class (e.g., for some given smoothness constant $L$ and strongly convexity constant $\mu$). If PEP can find a potential function that does not depend on $x^\star$, then {\os} can be directly applied since it only requires a numerically computable potential function. It is appealing to accelerate first-order methods in practice by applying {\os} to the potential functions generated by performance estimation.

\subsection{Convex composite optimization and proximal gradient} \label{sec:proxgrad}

{\os} works primarily in the smooth unconstrained setting. \Cref{sec:nonconvex} shows one way to extend {\os} to smooth nonconvex optimization. 
Here, we discuss challenges and opportunities for {\os} in nonsmooth optimization.

\begin{ques} \label{ques:prox}
Can {\os} be used to accelerate the solution of convex composite optimization or general nonsmooth convex optimization problems?
\end{ques}

 For example, consider the convex composite optimization problem $\min_x \varphi (x) \assign f (x) + w (x)$, where $f$ is smooth convex and $w$ is possibly nonsmooth but
prox-friendly: that is, for any $\gamma>0$, the proximal mapping
\[ \tmop{prox}_{w / \gamma} (x) \assign \argmin_z~  \{ w (z) +
   \tfrac{\gamma}{2} \| z - x \|^2 \} \]
can be efficiently evaluated. 
The operator $\mathcal{G}_{\gamma} (x) \assign \gamma
(x - \tmop{prox}_{w / \gamma} ( x - \tfrac{1}{\gamma} \nabla f (x)))$, known as the gradient map,  generalizes gradient to this composite setting, and the parameter $\gamma>0$ generalizes the stepsize. 
The proximal gradient method updates the iterate as
\[ x^{k + 1} = x^k - \tfrac{1}{\gamma} \mathcal{G}_{\gamma} (x^k) . \]

One natural idea to accelerate convergence is to adapt the parameter $\gamma$ using {\os}.
However, a basic implementation of this idea does not work, as the ratio or hypergradient feedback is not a convex function of $\gamma$ (or $1/\gamma$) due to the complex dependence from $\mathcal{G}_{\gamma} (x)$. A partial solution in this context separates the stepsize multiplying the gradient map from the parameter inside the map: fix $\gamma = L$ in the gradient map and consider the update
\[ x^{k + 1} = x^k - P_k \mathcal{G}_L (x^k)\]
with stepsize $P_k$.
This iteration directly corresponds to gradient descent.
We can show an improved convergence rate in the following two settings:

\paragraph{Case 1.} If $f (x) = 0$, then $\mathcal{G}_L (x)$ is the gradient of the Moreau envelope
$\varphi^{1 / L} (x) = w^{1 / L} (x) \assign \min_z  \{w(z) + \tfrac{L}{2} \| z - x \|^2
\}$, which is smooth and convex. The update is essentially gradient descent on the envelope $w^{1 / L}$, and all the {\os} results apply. If $w$ is not prox-friendly, we can solve the proximal subproblem inexactly with an additional subroutine and accelerate gradient descent on the envelope.

\paragraph{Case 2.} When $f$ is a general smooth convex function, we can define the proximal hypergradient feedback
\[ h_x (P) \assign \tfrac{\varphi (x - P\mathcal{G}_L (x)) - \varphi (x)}{\|
   \mathcal{G}_L (x) \|^2} = \tfrac{f (x - P\mathcal{G}_L (x)) + w (x -
   P\mathcal{G}_L (x)) - [f (x) + w (x)]}{\| \mathcal{G}_L (x) \|^2} . \]
Since $w$ is nonsmooth, $h_x (P)$ will not be Lipschitz continuous
when $\| \mathcal{G}_L (x) \| \rightarrow 0$ and online subgradient method
cannot guarantee sublinear regret. To address this issue, we consider online proximal gradient update
\[ P_{k + 1} = \argmin_P ~ \Big\{ \tfrac{f(x^k - P_k\mathcal{G}_L (x^k))+ \langle \nabla f (x^k - P_k
   \mathcal{G}_L (x^k)) \mathcal{G}_L (x^k)^{\top},\, P - P_k \rangle}{\|
   \mathcal{G}_L (x^k) \|^2} + \tfrac{w (x^k - P\mathcal{G}_L (x^k))}{\|
   \mathcal{G}_L (x^k) \|^2} + \tfrac{1}{2 \eta} \| P - P_k \|_F^2 \Big\},
\]
where $h_x (P)$ has been partially linearized. 
This subproblem is a proximal update in $P$ with respect to the (scaled, translated) 
function $\tilde w(P) = w(x-P\mathcal G_L(x))$.
This function often, but not always, inherits from $w$ a proximal mapping that can be efficiently evaluated.
If $P$ is a diagonal matrix and $w$ is separable in the coordinates,
the proximal operator of $\tilde w$ can be evaluated with the same complexity as the proximal operator of $w$. 
Combined with the monotone landscape action, we can obtain a convergence guarantee in this case in terms of $\| \mathcal{G}_L(x) \|$: $\min_{1 \leq k \leq K} \| \mathcal{G}_L (x^k) \|^2 =\mathcal{O}
( \tfrac{1}{K})$. 

\begin{thm}[Informal] \label{thm:proxgrad}
Suppose $f$ is $L$-smooth convex and $w$ is convex. Then given any benchmark stepsize $\hat{P} \in \Pcal$, a variant of {\os} achieves 
\[ \min_{1 \leq k \leq K} \| \mathcal{G}_L (x^k) \|^2\leq \tfrac{\varphi(x^1) - \varphi(x^\star)}{K} \tfrac{1}{\frac{1}{K} \max \{\sum_{k=1}^K h_{x^k}(\hat{P}) - \frac{L}{2}\|P_1 - \hat{P}\|_F^2, 0\}}.\]
In particular, taking $P_1 = \hat{P} = \frac{1}{L} I$ yields $\min_{1 \leq k \leq K} \| \mathcal{G}_L (x^k) \|^2 \leq 
 \tfrac{2L [\varphi(x^1) - \varphi^\star] }{K}$.
\end{thm}

A detailed analysis is available in the appendix.

%% file: conclusions.tex
\section{Conclusions}

The second part of our paper explores the practical aspects of {\os} and outlines several promising directions for future research. We design and analyze {\ob}, a practical variant of {\os} that integrates insights into algorithmic behavior with heavy-ball momentum. We also extend the {\os} framework to smooth nonconvex optimization problems. Numerical experiments demonstrate that {\os} achieves competitive performance across a range of real-world optimization tasks.

%% file: app_prac.tex
\section{Proof of results in \Cref{sec:prac}}

\subsection{Proof of Lemma \ref{lem:osgmhx-divergence}}
\input{proofs/proof-lem-osgmhx-divergence.tex}

\subsection{Proof of Theorem \ref{thm:orbit-stability}}
\input{proofs/proof-stability.tex}

\subsubsection{Proof of Lemma \ref{lem:orbit}}
\input{proofs/proof-lem-orbit.tex}

\subsubsection{Proof of Lemma \ref{lem:jacobian-blkdiag}}
\input{proofs/proof-jacob-blkdiag.tex}

\subsubsection{Proof of Lemma \ref{lem:jacobian-spec-radius}}
\input{proofs/proof-jacob-specrad.tex}
\subsection{Proof of Lemma \ref{lem:heavyball-potential}}
\input{proofs/proof-hb-potential.tex}

\subsection{Proof of Lemma \ref{lem:hb-properties}}
\input{proofs/proof-hb-properties.tex}

\subsection{Proof of Theorem \ref{thm:heavyball-reduction}}
\input{proofs/proof-hb-reduction.tex}

\subsection{Proof of Theorem \ref{thm:ohblook}}
\input{proofs/proof-hb-lookahead.tex}


%% file: proofs/proof-stability.tex
The proof of \Cref{thm:orbit-stability} requires several lemmas. The proof of each lemma is deferred to the end of the section.

\begin{lem} \label{lem:jacobian-blkdiag}
  Let $J_{{\hdm}} (\alpha, \hat{z})$ and
  $J_{\os} (\alpha, \hat{z})$ denote the Jacobians of the
  dynamical systems $\mathcal{F}_{{\hdm}}$ and
  $\mathcal{F}_{\os}$ respectively with $\eta = \frac{1}{L}$. Then each of
  $J_{{\texttt{HDM}}} \left( \alpha^{\star}, {\hat{z}^{\star}_i} 
  \right)$ and $J_{\os} (\alpha, \hat{z}^{\star}_i)$, $i = 1,
  2$, is similar to the block diagonal matrix by a common permutation matrix:
  \[ J_{{\hdm}} \left( \alpha^{\star}, {\hat{z}^{\star}_i} 
     \right) \sim \tmop{Diag} (1 - \alpha^{\star} \lambda_2, \ldots, 1 -
     \alpha^{\star} \lambda_{n - 1}, J_{{\hdm}}^i) \quad
     \text{and} \quad J_{\os} \left( \alpha^{\star},
     {\hat{z}^{\star}_i}  \right) \sim \tmop{Diag} (1 - \alpha^{\star}
     \lambda_2, \ldots, 1 - \alpha^{\star} \lambda_{n - 1},
     J_{\os}^i), \]
  where $J_{{\hdm}}^i$ and $J_{\os}^i$, $i = 1,
  2$, are $3 \times 3$ matrices defined by
  \begin{align}
    J_{{\hdm}}^1 \assign D \left(\begin{smallmatrix}
      \frac{\kappa - 1}{2 \kappa} & - \sqrt{\kappa^2 + 1} & \pm
      \frac{\sqrt{\kappa^2 + 1}}{\kappa}\\
      \frac{- \kappa (\kappa + 1)}{(\kappa^2 + 1)^{3 / 2}} & \frac{\kappa}{\kappa^2 + 1} & \frac{\pm 1}{\kappa^2 + 1} \\
      \frac{\mp (\kappa + 1)}{(\kappa^2 + 1)^{3 / 2}} & \frac{\mp 1}{\kappa^2 + 1} & \frac{-1}{\kappa(\kappa^2 + 1)}
    \end{smallmatrix}\right) D^{-1}, & \quad J_{{\hdm}}^2 \assign
    D \left(\begin{smallmatrix}
      \frac{\kappa - 1}{2 \kappa} & \sqrt{\kappa^2 + 1} & \pm
      \frac{\sqrt{\kappa^2 + 1}}{\kappa}\\
      \frac{\kappa (\kappa + 1)}{(\kappa^2 + 1)^{3 / 2}} & \frac{\kappa}{\kappa^2 + 1} & \frac{\mp 1}{\kappa^2 + 1} \\
      \frac{\mp (\kappa + 1)}{(\kappa^2 + 1)^{3 / 2}} & \frac{\pm 1}{\kappa^2 + 1} & \frac{-1}{\kappa(\kappa^2 + 1)}
    \end{smallmatrix}\right) D^{-1} ; \nonumber\\
    J_{\os}^1 \assign D \left(\begin{smallmatrix}
      \frac{\kappa - 1}{2 \kappa} & - \sqrt{\kappa^2 + 1} & \pm
      \frac{\sqrt{\kappa^2 + 1}}{\kappa}\\
      \frac{- \kappa (\kappa + 1)}{(\kappa^2 + 1)^{3 / 2}} & \frac{-
      \kappa^2}{\kappa^2 + 1} & \frac{\pm \kappa}{\kappa^2 + 1}\\
      \frac{\mp (\kappa + 1)}{(\kappa^2 + 1)^{3 / 2}} & \frac{\mp
      \kappa}{\kappa^2 + 1} & \frac{1}{\kappa^2 + 1}
    \end{smallmatrix}\right) D^{-1}, & \quad J_{\os}^2 \assign
    D \left(\begin{smallmatrix}
      \frac{\kappa - 1}{2 \kappa} & \sqrt{\kappa^2 + 1} & \pm
      \frac{\sqrt{\kappa^2 + 1}}{\kappa}\\
      \frac{\kappa (\kappa + 1)}{(\kappa^2 + 1)^{3 / 2}} & \frac{-
      \kappa^2}{\kappa^2 + 1} & \frac{\mp \kappa}{\kappa^2 + 1}\\
      \frac{\mp (\kappa + 1)}{(\kappa^2 + 1)^{3 / 2}} & \frac{\pm
      \kappa}{\kappa^2 + 1} & \frac{1}{\kappa^2 + 1}
    \end{smallmatrix}\right) D^{-1}, \nonumber
  \end{align} 
  and $D = \diag(\tfrac{\kappa - 1}{L (\kappa + 1)}, 1, 1)$ is a $3 \times 3$ invertible diagonal matrix.
\end{lem}

The stability of the period-$2$ orbit is determined by the spectral radius of the products
$J_{{\hdm}} \left( \alpha^{\star}, {\hat{z}^{\star}_1}  \right)
J_{{\hdm}} \left( \alpha^{\star}, {\hat{z}^{\star}_2}  \right)$
and $J_{\os} \left( \alpha^{\star}, {\hat{z}^{\star}_1} 
\right) J_{\os} \left( \alpha^{\star}, {\hat{z}^{\star}_2} 
\right)$, which equals the maximum absolute value of the eigenvalues. By
\Cref{lem:jacobian-blkdiag}, the spectrum of these two matrices are
\begin{align}
  \sigma \left( J_{{\hdm}} \left( \alpha^{\star},
  {\hat{z}^{\star}_1}  \right) J_{{\hdm}} \left( \alpha^{\star},
  {\hat{z}^{\star}_2}  \right) \right) & = \{ (1 - \alpha^{\star}
  \lambda_2)^2, \ldots, (1 - \alpha^{\star} \lambda_{n - 1})^2 \} \cup \sigma
  (J_{{\hdm}}^1 J_{{\hdm}}^2) ; \nonumber\\
  \sigma \left( J_{\os} \left( \alpha^{\star},
  {\hat{z}^{\star}_1}  \right) J_{\os} \left( \alpha^{\star},
  {\hat{z}^{\star}_2}  \right) \right) & = \{ (1 - \alpha^{\star}
  \lambda_2)^2, \ldots, (1 - \alpha^{\star} \lambda_{n - 1})^2 \} \cup \sigma
  (J_{\os}^1 J_{\os}^2) . \nonumber
\end{align}

Since the eigenvalues of $A$ are distinct, i.e., $L = \lambda_1 > \lambda_2 >
\cdots > \lambda_{n - 1} > \lambda_n = \mu$, and $\alpha^{\star} = \tfrac{2}{L
+ \mu}$, the first $n - 2$ eigenvalues satisfy
\[ - \tfrac{\kappa - 1}{\kappa + 1} = 1 - \alpha^{\star} L < 1 -
   \alpha^{\star} \lambda_i < 1 - \alpha^{\star} \mu = \tfrac{\kappa -
   1}{\kappa + 1} \quad \text{for } i = 2, \ldots, n - 1. \]
Then $(1 - \alpha^{\star} \lambda_i)^2 < 1$ for all $i = 2, \ldots, n - 1$ and
the stability of the orbit $\{ (\alpha^{\star}, \hat{z}^{\star}_1),
(\alpha^{\star}, \hat{z}^{\star}_2) \}$ is determined by the spectral radius
of $J_{{\hdm}}^1 J_{{\hdm}}^2$ and
$J_{\os}^1 J_{\os}^2$ for the corresponding
system.

\begin{lem} \label{lem:jacobian-spec-radius}
  The $3 \times 3$ matrices $J_{{\hdm}}^1
  J_{{\hdm}}^2$ and $J_{\os}^1
  J_{\os}^2$ have spectral radius 
\begin{align*}
\rho(J_{{\hdm}}^1 J_{{\hdm}}^2) ={} & \max \Big\{ -\tfrac{3 \kappa ^3-7 \kappa -9}{8 (\kappa +1) \kappa ^2}\pm \Big[ \tfrac{\sqrt{-(\kappa -3)^2 \kappa ^2 (\kappa +1)^3 (7 \kappa -9) (\kappa ^2+1)^5}}{8 (\kappa ^2+1)^{5/2} (\kappa +1) \kappa ^3}-\tfrac{5}{8 (\kappa +1)}\Big] \Big\} \leq 1	\\
  \rho (J_{\os}^1 J_{\os}^2) ={} & \tfrac{3 \kappa + 1}{2 \kappa} > 1   
\end{align*}
  

\end{lem}

Observe that the spectral radius $\rho(J_{{\hdm}}^1 J_{{\hdm}}^2)$ asymptotically increases to $\tfrac{1}{2}$ when $\kappa \rightarrow \infty$, implying that the eigenvalues of the $3 \times 3$ submatrix $J_{{\hdm}}^1 J_{{\hdm}}^2$ are not active for the spectral radius of the whole Jacobian when $\kappa$ is large. 
In this case, the spectral radius of $J_{{\hdm}} \left( \alpha^{\star}, {\hat{z}^{\star}_1}  \right) J_{{\hdm}} \left( \alpha^{\star}, {\hat{z}^{\star}_2}  \right)$ is the maximum of $(1-\alpha^\star \lambda_2)^2$ and $(1-\alpha^\star \lambda_{n-1})^2$. 

\begin{rem}
	While our result is established for $\eta = \frac{1}{L}$, the same quantitative result holds for $\eta = \frac{\theta}{L}, \theta \in (0, 1]$.
\end{rem}

%% file: proofs/proof-jacob-blkdiag.tex
The Jacobians of ${\Fhdm}$ and ${\Fosgm}$ are
\begin{align}
  J_{{\hdm}} (\alpha, \hat{z}) & = \left(\begin{array}{cc}
    1 - \eta \tfrac{\langle \hat{z}, \Lambda^3 \hat{z} \rangle}{\langle
    \hat{z}, \Lambda^2 \hat{z} \rangle} & v (\alpha, \hat{z})^{\top}\\[5pt]
    \tfrac{- 1}{\| u (\alpha, \hat{z}) \|} \Pi(\alpha, \hat{z})
    \Big( 1 - \eta \tfrac{\langle \hat{z}, \Lambda^3 \hat{z} \rangle}{\langle
    \hat{z}, \Lambda^2 \hat{z} \rangle} \Big) \Lambda \hat{z} & \tfrac{1}{\|
    u (\alpha, \hat{z}) \|} \Pi(\alpha, \hat{z}) [(I - \alpha^+
    (\hat{z}) \Lambda) - \Lambda \hat{z} v (\alpha, \hat{z})^{\top}]
  \end{array}\right), \nonumber\\
  J_{{\os}} (\alpha, \hat{z}) & = \left(\begin{array}{cc}
    1 - \eta \tfrac{\langle \hat{z}, \Lambda^3 \hat{z} \rangle}{\langle
    \hat{z}, A^2 \hat{z} \rangle} & v (\alpha, \hat{z})^{\top}\\[5pt]
    \tfrac{- 1}{\| u (\alpha, \hat{z}) \|} \Pi(\alpha, \hat{z})
    \Lambda \hat{z} & \tfrac{1}{\| u (\alpha, \hat{z}) \|} \Pi(\alpha, \hat{z}) (I - \alpha \Lambda)
  \end{array}\right), \nonumber
\end{align}
where $\Pi(\alpha, \hat{z}) \assign I - \tfrac{u(\alpha, \hat{z}) u (\alpha, \hat{z})^{\top}}{\| u (\alpha, \hat{z}) \|^2}$, 
$u (\alpha, \hat{x}) \assign (I - \alpha \Lambda) \hat{z}$,
and $v(\alpha, \hat{z}) \assign \tfrac{2 \alpha \eta}{\langle \hat{z}, \Lambda^2
\hat{z} \rangle} \Lambda^2 ( \tfrac{\langle \hat{z}, \Lambda^3 \hat{z}
\rangle}{\langle \hat{z}, \Lambda^2 \hat{z} \rangle} I - \Lambda )
\hat{z}$. 
To evaluate the Jacobians at $(\alpha^{\star}, \hat{z}^{\star}_1)$ and $(\alpha^{\star}, \hat{z}^{\star}_2)$, we calculate each block one by one.
Note that $J_{{\hdm}}$ and $J_{\os}$ only differ in the last row.

\paragraph{$(1,1)$-block of $J_{{\hdm}}(\alpha^{\star}, \hat{z}^{\star}_i)$ and $J_{\os}(\alpha^{\star}, \hat{z}^{\star}_i)$.}
Equation \eqref{eqn:pf-orbit-5} and $\eta = \tfrac{1}{L}$
imply
\begin{equation} \label{eqn:pf-jacobian-sim-0}
1 - \eta \tfrac{\langle \hat{z}^{\star}_i, \Lambda^3 \hat{z}^{\star}_i
\rangle}{\langle \hat{z}^{\star}_i, \Lambda^2 \hat{z}^{\star}_i \rangle} =
1 - \tfrac{L + \mu}{2 L} = \tfrac{\kappa - 1}{2 \kappa} \quad \text{for } i
= 1, 2. 
\end{equation}

\paragraph{$(1,2)$-block of $J_{{\hdm}}(\alpha^{\star}, \hat{z}^{\star}_i)$ and $J_{\os}(\alpha^{\star}, \hat{z}^{\star}_i)$.}
Using \eqref{eqn:pf-orbit-3}--\eqref{eqn:pf-orbit-5} and $\eta =
\tfrac{1}{L}$, we have
\[ v (\alpha^{\star}, \hat{z}^{\star}_i) = \tfrac{2 \alpha^{\star}
   \eta}{\langle \hat{z}^{\star}_i, \Lambda^2 \hat{z}^{\star}_i \rangle}
   \Lambda^2 ( \tfrac{\langle \hat{z}^{\star}_i, \Lambda^3
   \hat{z}^{\star}_i \rangle}{\langle \hat{z}^{\star}_i, \Lambda^2
   \hat{z}^{\star}_i \rangle} I - \Lambda ) \hat{z}^{\star}_i =
   \tfrac{1}{L^3 \gamma^2} \Lambda^2 (I - \alpha^{\star} \Lambda)
   \hat{z}^{\star}_i, \]
which, together with \eqref{eqn:pf-orbit-7} and \eqref{eqn:pf-orbit-8},
implies
\begin{align}
  v (\alpha^{\star}, \hat{z}^{\star}_1) & = \tfrac{\kappa - 1}{L^3 \gamma^2
  (\kappa + 1)} \Lambda^2 \hat{z}^{\star}_2 = \tfrac{\kappa - 1}{L^3 \gamma^2
  (\kappa + 1)} [- L^2 \gamma e_1 + \mu^2 \delta e_n] = - \tfrac{(\kappa - 1) 
  \sqrt{\kappa^2 + 1}}{L (\kappa + 1)} e_1 + \tfrac{\pm (\kappa - 1) 
  \sqrt{\kappa^2 + 1}}{L \kappa (\kappa + 1)} e_n ; \label{eqn:pf-jacobian-sim-0.1}\\
  v (\alpha^{\star}, \hat{z}^{\star}_2) & = \tfrac{\kappa - 1}{L^3 \gamma^2
  (\kappa + 1)} \Lambda^2 \hat{z}^{\star}_1 = \tfrac{\kappa - 1}{L^3 \gamma^2
  (\kappa + 1)} [L^2 \gamma e_1 + \mu^2 \delta e_n] = \tfrac{(\kappa - 1) 
  \sqrt{\kappa^2 + 1}}{L (\kappa + 1)} e_1 + \tfrac{\pm (\kappa - 1) 
  \sqrt{\kappa^2 + 1}}{L \kappa (\kappa + 1)} e_n . \label{eqn:pf-jacobian-sim-0.2}
\end{align}
\paragraph{$(2,1)$-block of $J_{\os}(\alpha^{\star}, \hat{z}^{\star}_i)$.}
Equations \eqref{eqn:pf-orbit-7}--\eqref{eqn:pf-orbit-8} imply
\begin{equation} \label{eqn:pf-jacobian-sim-1} 
u (\alpha^{\star}, \hat{z}^{\star}_1) = (I - \alpha^{\star} \Lambda)
\hat{z}^{\star}_1 = \tfrac{\kappa - 1}{\kappa + 1} \hat{z}^{\star}_2, \quad
u (\alpha^{\star}, \hat{z}^{\star}_2) = (I - \alpha^{\star} \Lambda)
\hat{z}^{\star}_2 = \tfrac{\kappa - 1}{\kappa + 1} \hat{z}^{\star}_1, \quad
\| u (\alpha^{\star}, \hat{z}^{\star}_1) \| = \| u (\alpha^{\star},
\hat{z}^{\star}_2) \| = \tfrac{\kappa - 1}{\kappa + 1},
\end{equation}
and hence
\begin{equation}\label{eqn:pf-jacobian-sim-2}
I - \tfrac{u (\alpha^{\star}, \hat{z}^{\star}_1) u (\alpha^{\star},
\hat{z}^{\star}_1)^{\top}}{\| u (\alpha^{\star}, \hat{z}^{\star}_1) \|^2} =
I - \hat{z}^{\star}_2 {(\hat{z}^{\star}_2)}^{\top} \quad \text{and} \quad I -
\tfrac{u (\alpha^{\star}, \hat{z}^{\star}_2) u (\alpha^{\star},
\hat{z}^{\star}_2)^{\top}}{\| u (\alpha^{\star}, \hat{z}^{\star}_2) \|^2} =
I - \hat{z}^{\star}_1 {(\hat{z}^{\star}_1)}^{\top} .
\end{equation}
Observe that ${(\hat{z}^{\star}_2)}^{\top} A \hat{z}^{\star}_1 =
{(\hat{z}^{\star}_1)}^{\top} A \hat{z}^{\star}_2 = - L \gamma^2 + \mu \delta^2$,
together with \eqref{eqn:pf-jacobian-sim-1}, we have
\begin{align}
  \tfrac{- 1}{\| u (\alpha^{\star}, \hat{z}^{\star}_1) \|} \Pi(\alpha^{\star}, \hat{z}^{\star}_1) \Lambda \hat{z}^{\star}_1 
  & = \tfrac{- 1}{\| u (\alpha^{\star}, \hat{z}^{\star}_1) \|} ( I - \tfrac{u (\alpha^{\star}, \hat{z}^{\star}_1) u (\alpha^{\star},
  \hat{z}^{\star}_1)^{\top}}{\| u (\alpha^{\star}, \hat{z}^{\star}_1) \|^2}
  ) \Lambda \hat{z}^{\star}_1 \nonumber\\
  & = - \tfrac{\kappa + 1}{\kappa - 1}
  ( \Lambda \hat{z}^{\star}_1 - \hat{z}^{\star}_2 
  {(\hat{z}^{\star}_2)}^{\top} \Lambda \hat{z}^{\star}_1 ) \nonumber\\
  & = - \tfrac{\kappa + 1}{\kappa - 1} [\Lambda \hat{z}^{\star}_1 - (- L
  \gamma^2 + \mu \delta^2) \hat{z}^{\star}_2] \nonumber\\
  & = - \tfrac{\kappa + 1}{\kappa - 1} [(L \gamma e_1 + \mu \delta e_n) - (-
  L \gamma^2 + \mu \delta^2) (- \gamma e_1 + \delta e_n)] \nonumber\\
  & = - \tfrac{\kappa + 1}{\kappa - 1} [(L \gamma - L \gamma^3 + \mu \gamma
  \delta^2) e_1 + (\mu \delta + L \gamma^2 \delta - \mu \delta^3) e_n]
  \nonumber\\
  & = \tfrac{- L \kappa (\kappa + 1)^2}{(\kappa - 1) (\kappa^2 + 1)^{3 / 2}}
  e_1 \mp \tfrac{L (\kappa + 1)^2}{(\kappa - 1) (\kappa^2 + 1)^{3 / 2}} e_n\label{eqn:pf-jacobian-sim-3}
\end{align}
\begin{align}
  \tfrac{- 1}{\| u (\alpha^{\star}, \hat{z}^{\star}_1) \|} \Pi(\alpha^{\star}, \hat{z}^{\star}_2) \Lambda \hat{z}^{\star}_2
  & = \tfrac{- 1}{\| u (\alpha^{\star}, \hat{z}^{\star}_2) \|} ( I - \tfrac{u (\alpha^{\star}, \hat{z}^{\star}_2) u (\alpha^{\star},
  \hat{z}^{\star}_2)^{\top}}{\| u (\alpha^{\star}, \hat{z}^{\star}_2) \|^2}
  ) \Lambda \hat{z}^{\star}_2 \nonumber\\
  & = - \tfrac{\kappa + 1}{\kappa - 1}
  ( \Lambda \hat{z}^{\star}_2 - \hat{z}^{\star}_1 
  {(\hat{z}^{\star}_1)}^{\top} \Lambda \hat{z}^{\star}_2 ) \nonumber\\
  & = - \tfrac{\kappa + 1}{\kappa - 1} [\Lambda \hat{z}^{\star}_2 - (- L
  \gamma^2 + \mu \delta^2) \hat{z}^{\star}_1] \nonumber\\
  & = - \tfrac{\kappa + 1}{\kappa - 1} [(- L \gamma e_1 + \mu \delta e_n) -
  (- L \gamma^2 + \mu \delta^2) (\gamma e_1 + \delta e_n)] \nonumber\\
  & = - \tfrac{\kappa + 1}{\kappa - 1} [(- L \gamma + L \gamma^3 - \mu \gamma
  \delta^2) e_1 + (\mu \delta + L \gamma^2 \delta - \mu \delta^3) e_n]
  \nonumber\\
  & = \tfrac{L \kappa (\kappa + 1)^2}{(\kappa - 1) (\kappa^2 + 1)^{3 / 2}}
  e_1 \mp \tfrac{L (\kappa + 1)^2}{(\kappa - 1) (\kappa^2 + 1)^{3 / 2}} e_n .
  \label{eqn:pf-jacobian-sim-4}
\end{align}

\paragraph{$(2,1)$-block of $J_{\hdm}(\alpha^{\star}, \hat{z}^{\star}_i)$.}
Note that the $(2,1)$-block of $J_{\hdm}(\alpha^{\star}, \hat{z}^{\star}_i)$ is a scalar multiple of the $(2,1)$-block of $J_{\os}(\alpha^{\star}, \hat{z}^{\star}_i)$ with the scalar $1 - \eta \tfrac{\langle \hat{z}^{\star}_i, \Lambda^3 \hat{z}^{\star}_i \rangle}{\langle \hat{z}^{\star}_i, \Lambda^2 \hat{z}^{\star}_i \rangle}$. Then using \eqref{eqn:pf-jacobian-sim-0} and \eqref{eqn:pf-jacobian-sim-3}--\eqref{eqn:pf-jacobian-sim-4}, we have
\begin{align*}
    \tfrac{- 1}{\| u (\alpha^{\star}, \hat{z}^{\star}_1) \|} \Pi(\alpha^{\star}, \hat{z}^{\star}_1) \big( 1 - \eta \tfrac{\langle \hat{z}^{\star}_1, \Lambda^3 \hat{z}^{\star}_1 \rangle}{\langle \hat{z}^{\star}_1, \Lambda^2 \hat{z}^{\star}_1 \rangle} \big) \Lambda \hat{z}^{\star}_1
    &= \tfrac{- L \kappa (\kappa + 1)^2}{2 \kappa (\kappa^2 + 1)^{3 / 2}}e_1 \mp \tfrac{L (\kappa + 1)^2}{2 \kappa (\kappa^2 + 1)^{3 / 2}} e_n;\\
    \tfrac{- 1}{\| u (\alpha^{\star}, \hat{z}^{\star}_2) \|} \Pi(\alpha^{\star}, \hat{z}^{\star}_2) \big( 1 - \eta \tfrac{\langle \hat{z}^{\star}_2, \Lambda^3 \hat{z}^{\star}_2 \rangle}{\langle \hat{z}^{\star}_2, \Lambda^2 \hat{z}^{\star}_2 \rangle} \big) \Lambda \hat{z}^{\star}_2
    &= \tfrac{L \kappa (\kappa + 1)^2}{2 \kappa (\kappa^2 + 1)^{3 / 2}}
  e_1 \mp \tfrac{L (\kappa + 1)^2}{2 \kappa (\kappa^2 + 1)^{3 / 2}} e_n .
\end{align*}

\paragraph{$(2,2)$-block of $J_{\os}(\alpha^{\star}, \hat{z}^{\star}_i)$.}
Let $I_{1, n} \assign [e_1, ~e_n] \in \mathbb{R}^{n \times 2}$ denote the first and last column of an $n \times n$ identity matrix.
Using \eqref{eqn:pf-jacobian-sim-1} and \eqref{eqn:pf-jacobian-sim-2}, we have
\begin{align}
  \tfrac{1}{\| u (\alpha^{\star}, \hat{z}^{\star}_1) \|} \Pi(\alpha^{\star}, \hat{z}^{\star}_1) (I - \alpha^{\star} \Lambda) 
  & = \tfrac{\kappa + 1}{\kappa - 1} \big( I - \hat{z}^{\star}_2 {(\hat{z}^{\star}_2)}^{\top} \big) (I - \alpha^{\star} \Lambda) 
  = \tfrac{\kappa + 1}{\kappa - 1} (I - \alpha^{\star} \Lambda) - \tfrac{\kappa + 1}{\kappa - 1} \hat{z}^{\star}_2 [(I - \alpha^{\star} \Lambda)
  \hat{z}^{\star}_2]^{\top} \nonumber\\
  & = \tfrac{\kappa + 1}{\kappa - 1} (I - \alpha^{\star} \Lambda) - \hat{z}^{\star}_2 {(\hat{z}^{\star}_1)}^{\top} \nonumber\\
  & = \tfrac{\kappa + 1}{\kappa - 1} (I - \alpha^{\star} \Lambda) - I_{1, n}
  \left(\begin{smallmatrix}
    - \gamma^2 & - \gamma \delta\\
    \gamma \delta & \delta^2
  \end{smallmatrix}\right) I^{\top}_{1, n}, \label{eqn:pf-jacobian-sim-5}\\
  \tfrac{1}{\| u (\alpha^{\star}, \hat{z}^{\star}_2) \|} \Pi(\alpha^{\star}, \hat{z}^{\star}_2) (I - \alpha^{\star} \Lambda) 
  & = \tfrac{\kappa + 1}{\kappa - 1} \big( I - \hat{z}^{\star}_1 {(\hat{z}^{\star}_1)}^{\top} \big) (I - \alpha^{\star} \Lambda) 
  = \tfrac{\kappa + 1}{\kappa - 1} (I - \alpha^{\star} \Lambda) - \tfrac{\kappa + 1}{\kappa - 1} \hat{z}^{\star}_1 [(I - \alpha^{\star} \Lambda) \hat{z}^{\star}_1]^{\top} \nonumber\\
  &= \tfrac{\kappa + 1}{\kappa - 1} (I - \alpha^{\star} \Lambda) - \hat{z}^{\star}_1 {(\hat{z}^{\star}_2)}^{\top} \nonumber\\
  & = \tfrac{\kappa + 1}{\kappa - 1} (I - \alpha^{\star} \Lambda) - I_{1, n}
  \left(\begin{smallmatrix}
    - \gamma^2 & \gamma \delta\\
    - \gamma \delta & \delta^2
  \end{smallmatrix}\right) I^{\top}_{1, n}. \label{eqn:pf-jacobian-sim-6}
\end{align}
Both \eqref{eqn:pf-jacobian-sim-5} and \eqref{eqn:pf-jacobian-sim-6} have $1 - \alpha^{\star} \lambda_i$ at the $(i, i)$-entry for $i
= 2, \ldots, n - 1$ and the only other non-zero entries lie on the four
corners indexed by $(1, 1)$, $(1, n)$, $(n, 1)$, and $(n, n)$.

\paragraph{$(2,2)$-block of $J_{\hdm}(\alpha^{\star}, \hat{z}^{\star}_i)$.}
Equation \eqref{eqn:pf-orbit-6} reads as $[\alpha^{\star}]^+
(\hat{z}^{\star}_i) = \alpha^{\star}$ for $i = 1, 2$, and thus the $(2, 2)$-block of
$J_{{\hdm}} (\alpha^{\star}, \hat{z}_i^{\star})$ becomes
\[ \tfrac{1}{\| u (\alpha^{\star}, \hat{z}^{\star}_i) \|} \Pi(\alpha^{\star}, \hat{z}^{\star}_i) [(I - [\alpha^{\star}]^+ (\hat{z}^{\star}_i) \Lambda) - \Lambda
   \hat{z}^{\star}_i v (\alpha^{\star}, \hat{z}^{\star}_i)^{\top}] =
   \tfrac{1}{\| u (\alpha^{\star}, \hat{z}^{\star}_i) \|} \Pi(\alpha^{\star}, \hat{z}^{\star}_i) [(I - \alpha^{\star} \Lambda) - \Lambda \hat{z}^{\star}_i v
   (\alpha^{\star}, \hat{z}^{\star}_i)^{\top}], \]
in which is first term has been evaluated in \eqref{eqn:pf-jacobian-sim-5} and \eqref{eqn:pf-jacobian-sim-6}. 
Using \eqref{eqn:pf-jacobian-sim-0.1}--\eqref{eqn:pf-jacobian-sim-0.2} and \eqref{eqn:pf-jacobian-sim-3}--\eqref{eqn:pf-jacobian-sim-4}, the second term is evaluated as
\begin{align}
  - \tfrac{1}{\| u (\alpha^{\star}, \hat{z}^{\star}_1) \|} \Pi(\alpha^{\star}, \hat{z}^{\star}_1) \Lambda \hat{z}^{\star}_1 v (\alpha^{\star},
  \hat{z}^{\star}_1)^{\top} & = - \tfrac{1}{L^3 \gamma^2} I_{1, n}
  \left(\begin{smallmatrix}
    L \gamma - L \gamma^3 + \mu \gamma \delta^2\\
    \mu \delta + L \gamma^2 \delta - \mu \delta^3
  \end{smallmatrix}\right) \left(\begin{smallmatrix}
    - L^2 \gamma\\
    \mu^2 \delta
  \end{smallmatrix}\right)^{\top} I_{1, n}^{\top} ; \label{eqn:pf-jacobian-sim-7}\\
  - \tfrac{1}{\| u (\alpha^{\star}, \hat{z}^{\star}_2) \|} \Pi(\alpha^{\star}, \hat{z}^{\star}_2) \Lambda \hat{z}^{\star}_2 v (\alpha^{\star},
  \hat{z}^{\star}_2)^{\top} & = - \tfrac{1}{L^3 \gamma^2} I_{1, n}
  \left(\begin{smallmatrix}
    - L \gamma + L \gamma^3 - \mu \gamma \delta^2\\
    \mu \delta + L \gamma^2 \delta - \mu \delta^3
  \end{smallmatrix}\right) \left(\begin{smallmatrix}
    L^2 \gamma\\
    \mu^2 \delta
  \end{smallmatrix}\right)^{\top} I^{\top}_{1, n}, \label{eqn:pf-jacobian-sim-8}
\end{align}
The non-zero entries of \eqref{eqn:pf-jacobian-sim-7} and \eqref{eqn:pf-jacobian-sim-8} again only lie on the four corners indexed by $(1, 1)$, $(1, n)$, $(n, 1)$, and $(n, n)$. 
Hence, the Jacobians $J_{{\hdm}} (\alpha^{\star}, \hat{z}^{\star}_i)$ and
$J_{{\os}} (\alpha^{\star}, \hat{z}^{\star}_i)$ for $i = 1, 2$ all have the structure
\[ \left(\begin{array}{ccccccc}
     \times & \times & 0 &  & \cdots & 0 & \times\\
     \times & \times & 0 &  & \cdots &  & \times\\
     0 & 0 & 1 - \alpha^{\star} \lambda_2 &  &  &  & 0\\
     &  &  & \ddots &  &  & \\
     \vdots & \vdots &  &  & \ddots &  & \vdots\\
     0 & 0 &  &  &  & 1 - \alpha^{\star} \lambda_{n - 1} & 0\\
     \times & \times & 0 &  & \cdots & 0 & \times
   \end{array}\right) . \]
We can simultaneously permute the rows and columns so that all the four
Jacobians are similar to a block diagonal matrix that consists of $n - 2$
individual diagonal entries $1 - \alpha^{\star} \lambda_2, \ldots, 1 -
\alpha^{\star} \lambda_{n - 1}$ and one $3 \times 3$ dense matrix, denoted by
$J_{{\hdm}}^i$ and $J_{{\os}}^i$ for $i = 1, 2$.
After careful but tedious calculations, the four $3 \times 3$ matrices match their definition in the lemma.

%% file: proofs/proof-jacob-specrad.tex
The result is obtained by direct computation using Mathematica.

%% file: app_nonconvex.tex
\section{Proof of results in Section \ref{sec:nonconvex}}

\subsection{Proof of Lemma \ref{lem:hx-nonconvexity}}
\input{proofs/proof-lem-nonconvexity.tex}

\subsection{Proof of Lemma \ref{lem:property-hxreg}}
\input{proofs/proof-property-hxreg.tex}

\subsection{Proof of Theorem \ref{thm:nonconvex-reduction}}
\input{proofs/proof-thm-nonconvex-reduction.tex}

\subsection{Proof of Lemma \ref{lem:nonconvex-feedback-progress}}
\input{proofs/proof-lem-nonconvex-feedback-progress.tex}

\subsection{Proof of Theorem \ref{thm:nonconvex-global-conv}}
\input{proofs/proof-nonconvex-global-conv.tex}

%% file: proofs/proof-lem-nonconvexity.tex
Since $h_x (P)$ is $L$-smooth, it is at most $L$-weakly convex. The second result needs several auxiliary lemmas.

\begin{lem}[Hessian of PL-functions]
  \label{lem:hess-PL}Suppose $f$ is $L$-smooth, has $H$-Lipschitz Hessian, and
  is $\mu$-PL, then $\lambda_{\min} (\nabla^2 f (x)) \geq - \mu^{- 1} H \|
  \nabla f (x) \|$.
\end{lem}
\begin{proof}
	Let $x^{\star} = \Pi_{\mathcal{X}^{\star}} [x]$. Using $H$-Lipschitz Hessian,
we have
\begin{align}
  \nabla^2 f (x) ={} & \nabla^2 f (x^{\star}) + \nabla^2 f (x) - \nabla^2 f
  (x^{\star}) \nonumber\\
  \succeq{} & \nabla^2 f (x^{\star}) - \| \nabla^2 f (x) - \nabla^2 f
  (x^{\star}) \| \cdot I \nonumber\\
  \succeq{} & \nabla^2 f (x^{\star}) - H \| x - x^{\star} \| \cdot I.
  \tag{$H$-Lipschitz Hessian}
\end{align}

Since PL condition implies quadratic growth \cite{karimi2016linear}, 
$\tfrac{1}{2 \mu} \| \nabla f (x) \|^2 \geq f (x) - f^{\star} \geq
\tfrac{\mu}{2} \| x - x^{\star} \|^2$ and $\| x - x^{\star} \| \leq
\tfrac{1}{\mu} \| \nabla f (x) \|$. Hence $\nabla^2 f (x) \succeq \nabla^2 f
(x^{\star}) - \tfrac{H}{\mu} \| \nabla f (x) \| \cdot I \succeq - \mu^{- 1}
H \| \nabla f (x) \| \cdot I$ and this completes the proof.
\end{proof}

\begin{lem}[Hessian of $h_x (P)$]
  \label{lem:hess-hxP} Suppose $f$ has $H$-Lipschitz Hessian and $x \nin
  \mathcal{S}^{\star}$, then
  \[ \| \nabla^2 h_x (P_1) - \nabla^2 h_x (P_2) \| \leq H \| \nabla f (x) \|
     \| P_1 - P_2 \| . \]
\end{lem}

\begin{proof}
Recall that $\nabla h_x (P) = \tfrac{\nabla f (x - P \nabla f (x)) \nabla f
(x)^{\top}}{\| \nabla f (x) \|^2}$ and $\nabla^2 h_x (P) = (g g^{\top})
\otimes \nabla^2 f (x - P \nabla f (x))$ where $g = \tfrac{\nabla f (x)}{\|
\nabla f (x) \|}$. Since $\| g \| = 1$, $H$-Lipschitz Hessian implies
\begin{equation}
\| \nabla^2 h_x (P_1) - \nabla^2 h_x (P_2) \| \leq \| \nabla^2 f (x - P_1
  \nabla f (x)) - \nabla^2 f (x - P_2 \nabla f (x)) \| \nonumber\\
  \leq H \| \nabla f (x) \|  \| P_1 - P_2 \|
\end{equation}
and this completes the proof.
\end{proof}

\paragraph{Proof of \Cref{lem:hx-nonconvexity}.} 
Take any fixed $\bar{P} \in \Pcal$. Using \Cref{lem:hess-hxP}, we deduce that
\begin{align}
  \nabla^2 h_x (P) \succeq{} & \nabla^2 h_x (\bar{P}) - H \| \nabla f (x) \| \|
  P - \bar{P} \| \cdot I \nonumber\\
  \succeq{} & [\lambda_{\min} (\nabla^2 f (x - \bar{P} \nabla f (x))) - H \|
  \nabla f (x) \| \| P - \bar{P} \|] \cdot I \nonumber\\
  \succeq{} & [\lambda_{\min} (\nabla^2 f (x - \bar{P} \nabla f (x))) - H \|
  \nabla f (x) \| D] \cdot I \nonumber
\end{align}
When $f$ is additionally $\mu$-PL and $\bar{P} = 0$, we can further apply $\lambda_{\min}
(\nabla^2 f (x - \bar{P} \nabla f (x))) \geq - \mu^{- 1} H \| \nabla f (x) \|$ from \Cref{lem:hess-PL} and this completes the proof.

%% file: app_extensions.tex
\section{Proof of results in Section \ref{sec:extensions}}

\subsection{Details of Subsection \ref{sec:bb}} \label{app:bb}
\input{sec_bb.tex}

\subsection{Proof of Theorem \ref{thm:agd-quadpot}}
\input{proofs/proof-agd-quadpot.tex}

\subsection{Details of in Subsection \ref{sec:proxgrad}} \label{sec:app-proxgrad}
\input{sec_proxgrad.tex}

%% file: sec_bb.tex
In this section, we replace the online gradient descent step in
{\osgmhandsonhx} with online proximal point. For the generality of the result, we still use a matrix stepsize $P$.
\vspace{10pt}

\begin{algorithm}[H]
{\textbf{input} initial point $x^1$, initial stepsize $P_1 \in \Rbb^{n \times n}$, online proximal regularizer $\eta > 0$}

\For{$k = 1, 2, \dots$}{
$x^{k + 1/2} = x^{k} - P_k \nabla f(x^k)$ \\
Choose $x^{k+1}$ that satisfies $f(x^{k+1}) \leq \min \{ f(x^{k+1/2} - \frac{1}{L} \nabla f(x^{k+1/2})
   ), f(x^k) \}$ \\
$P_{k+1} =  \argmin_P  \{ h_{x^k} (P) + \tfrac{1}{2 \eta} \| P -
  P_k \|_F^2 \}$ \label{line:proxbb}
}
\caption{{\osgmhandsonhx} with online proximal point update \label{alg:bb-osgm}}
\end{algorithm}
\vspace{10pt}
To show the convergence of \Cref{alg:bb-osgm}, we again proceed in three steps: 1)
bounding suboptimality with progress 2) bounding progress with feedback 3)
bounding feedback with online learning. Since we essentially use the same
hypergradient feedback, step 1) and 2) are the same as in Part I. It suffices
to prove the regret guarantee of online proximal point and combine it with the lookahead action. Recall the definition of progress $h_k = \tfrac{f (x^{k + 1}) - f (x^k)}{\|
\nabla f (x^k) \|^2}$ and feedback $h_{x^k} (P_k) = \tfrac{f (x^{k + 1 / 2}) -
f (x^k)}{\| \nabla f (x^k) \|^2}$. We have, from \cite[Lemma 4.1]{gao2025gradient}, that

\begin{lem}[Feedback and progress \cite{gao2025gradient}] \label{lem:bb-0}
  Let $f$ be $L$-smooth and convex. Then $h_k \leq \min \{ h_{x^k} (P_k) -
  \tfrac{1}{2 L} \| \nabla h_{x^k} (P_k) \|_F^2, 0 \}$ with monotone
  lookahead action.
\end{lem}

The proof requires an additional technical lemma that shows a relation between the value of the feedback function
$h_x (P)$ and its Moreau envelope function $h_x^{1 / \rho} (P) \assign \min_Q  \{ h_x (Q) +
\tfrac{\rho}{2} \| Q - P \|_F^2 \}$.

\begin{lem}
  \label{lem:bb-1}Let $h^{1 / \rho} (x) \assign \min_z  \{ h (z) +
  \tfrac{\rho}{2} \| z - x \|^2 \}$ be the Moreau envelope of a differentiable convex function $f$. Then $h (x) - h^{1 / \rho} (x) - \tfrac{1}{2
  \rho} \| \nabla h (x) \|^2 \leq 0$.
\end{lem}

\begin{proof}
  By convexity, we have, for all $z$, that
  \begin{align}
    h (z) + \tfrac{\rho}{2} \| z - x \|^2 \geq{} & h (x) + \langle \nabla h (x),
    z - x \rangle + \tfrac{\rho}{2} \| z - x \|^2 \nonumber\\
    \geq{}& \min_d  \{ h (x) + \langle \nabla h (x), d \rangle +
    \tfrac{\rho}{2} \| d \|^2 \} \nonumber\\
    ={} & h (x) - \tfrac{1}{2 \rho} \| \nabla h (x) \|^2 . \nonumber
  \end{align}
    Minimizing the LHS, $h^{1 / \rho} (x) = \min_z  \{ h (z) + \tfrac{\rho}{2} \| z - x \|^2
     \} \geq h (x) - \tfrac{1}{2 \rho} \| \nabla h (x) \|^2$. This completes the proof.
\end{proof}
Using the optimality condition on the proximal subproblem on \Cref{line:proxbb} of \Cref{alg:bb-osgm}, for any $\hat{P}$, we deduce that
\begin{align}
  h_{x^k} (P_{k + 1}) + \tfrac{1}{2 \eta} \| P_{k + 1} - P_k \|_F^2 \leq{} &
  h_{x^k} (\hat{P}) + \tfrac{1}{2 \eta} \|\hat{P}- P_k \|_F^2 - \tfrac{1}{2 \eta} \|
  P_{k + 1} -\hat{P}\|_F^2 . \nonumber
\end{align}
Re-arranging the terms, $h_{x^k} (P_{k + 1}) - h_{x^k} (\hat{P}) \leq \tfrac{1}{2
\eta} \|\hat{P}- P_k \|_F^2 - \tfrac{1}{2 \eta} \| P_{k + 1} - P_k \|_F^2 -
\tfrac{1}{2 \eta} \| P_{k + 1} -\hat{P}\|_F^2$ and
\begin{align}
h_{x^k} (P_k) - h_{x^k} (\hat{P})  ={} & h_{x^k} (P_k) - h_{x^k} (P_{k + 1}) + h_{x^k} (P_{k + 1}) - h_{x^k} (\hat{P})
  \nonumber\\
  \leq{} & h_{x^k} (P_k) - h_{x^k} (P_{k + 1}) + \tfrac{1}{2 \eta} \|\hat{P}- P_k
  \|_F^2 - \tfrac{1}{2 \eta} \| P_{k + 1} - P_k \|_F^2 - \tfrac{1}{2 \eta} \|
  P_{k + 1} -\hat{P}\|_F^2 \nonumber\\
  ={} & h_{x^k} (P_k) - [ h_{x^k} (P_{k + 1}) + \tfrac{1}{2 \eta} \| P_{k +
  1} - P_k \|_F^2 ] + \tfrac{1}{2 \eta} [\|\hat{P}- P_k \|_F^2 - \| P_{k +
  1} -\hat{P}\|_F^2] \nonumber\\
  ={} & h_{x^k} (P_k) - h_{x^k}^{\eta} (P_k) + \tfrac{1}{2 \eta} [\|\hat{P}- P_k
  \|_F^2 - \| P_{k + 1} -\hat{P}\|_F^2] \nonumber\\
  \leq{} & \tfrac{\eta}{2} \| \nabla h_{x^k} (P_k) \|^2 + \tfrac{1}{2 \eta} [\|
 \hat{P}- P_k \|_F^2 - \| P_{k + 1} -\hat{P}\|_F^2] \nonumber
\end{align}

Finally, we bound the progress using the monotone lookahead action using \Cref{lem:bb-0}:
\begin{align}
\textstyle \sum_{k = 1}^K h_k \leq & \textstyle \min \{ \sum_{k = 1}^K h_{x^k} (P_k) -
  \tfrac{1}{2 L} \| \nabla h_{x^k} (P_k) \|_F^2, 0 \} \nonumber\\
  \leq & \textstyle \min \{ \sum_{k = 1}^K h_{x^k} (\hat{P}) + ( \tfrac{\eta}{2} -
  \tfrac{1}{2 L} ) \| \nabla h_{x^k} (P_k) \|^2 + \tfrac{1}{2 \eta} \|
  P_1 -\hat{P}\|_F^2, 0 \} \nonumber
\end{align}

and for $\eta = \tfrac{1}{L}$, $\sum_{k = 1}^K h_k \leq \sum_{k = 1}^K h_{x^k}
(\hat{P}) + \tfrac{L}{2} \| P_1 -\hat{P}\|_F^2$. Plugging the progress relation back
into the reduction in Part I \cite[Theorem 4.2]{gao2025gradient} establishes the convergence.

%% file: proofs/proof-agd-quadpot.tex
By definition of the {\agd} update, we deduce that
\begin{align}
  y^k ={} & x^k + \tfrac{1}{\sqrt{\kappa} + 1} (z^k - x^k) \nonumber\\
  ={} & \tfrac{\sqrt{\kappa}}{\sqrt{\kappa} + 1} x^k + \tfrac{1}{\sqrt{\kappa} +
  1} z^k \nonumber\\
  x^{k + 1} - x^{\star} ={} & y^k - \tfrac{1}{L} \nabla f (y^k) - x^{\star}
  \nonumber\\
  ={} & y^k - \tfrac{1}{L} (A y^k - A x^{\star}) - x^{\star} \nonumber\\
  ={} & ( I - \tfrac{1}{L} A ) y^k + \tfrac{1}{L} A x^{\star} -
  x^{\star} \nonumber\\
  ={} & ( I - \tfrac{1}{L} A ) (
  \tfrac{\sqrt{\kappa}}{\sqrt{\kappa} + 1} x^k + \tfrac{1}{\sqrt{\kappa} + 1}
  z^k ) + \tfrac{1}{L} A x^{\star} - x^{\star} \nonumber\\
  ={} & ( I - \tfrac{1}{L} A ) \tfrac{\sqrt{\kappa}}{\sqrt{\kappa} +
  1} (x^k - x^{\star}) - ( I - \tfrac{1}{L} A )
  \tfrac{1}{\sqrt{\kappa} + 1} (z^k - x^{\star}) \nonumber\\
  z^{k + 1} - x^{\star} ={} & ( 1 - \tfrac{1}{\sqrt{\kappa}} ) (z^k -
  x^{\star}) + \tfrac{1}{\sqrt{\kappa}} ( y^k - \tfrac{1}{\mu} (A y^k - A
  x^{\star}) - x^{\star} ) \nonumber\\
  ={} & ( 1 - \tfrac{1}{\sqrt{\kappa}} ) (z^k - x^{\star}) +
  \tfrac{1}{\sqrt{\kappa}} ( I - \tfrac{1}{\mu} A ) (y^k -
  x^{\star}) \nonumber\\
  ={} & ( 1 - \tfrac{1}{\sqrt{\kappa}} ) (z^k - x^{\star}) +
  \tfrac{1}{\sqrt{\kappa}} ( I - \tfrac{1}{\mu} A ) (
  \tfrac{\sqrt{\kappa}}{\sqrt{\kappa} + 1} x^k + \tfrac{1}{\sqrt{\kappa} + 1}
  z^k - x^{\star} ) \nonumber\\
  ={} & \tfrac{1}{\sqrt{\kappa} + 1} ( I - \tfrac{1}{\mu} A ) (x^k -
  x^{\star}) + \tfrac{\sqrt{\kappa}}{\sqrt{\kappa} + 1} ( I -
  \tfrac{1}{L} A ) (z^k - x^{\star}) \nonumber
\end{align}

On the other hand, we can express the potential function as follows:
\begin{align}
  \varphi_{\mu} (x, z) ={} & \tfrac{1}{2} \langle (z - x^{\star}), A (z -
  x^{\star}) \rangle + \tfrac{1}{2 \mu} \| A (x - x^{\star}) \|^2 \nonumber\\
  ={} & \tfrac{1}{2} \| A^{1 / 2} (z - x^{\star}) \|^2 + \tfrac{1}{2} \|
  \tfrac{1}{\sqrt{\mu}} A (x - x^{\star}) \|^2 \nonumber\\
  ={} & \tfrac{1}{2} \Bigg\| \Bigg(\begin{array}{c}
    \tfrac{1}{\sqrt{\mu}} A (x - x^{\star})\\
    A^{1 / 2} (z - x^{\star})
  \end{array}\Bigg) \Bigg\|^2 . \nonumber
\end{align}
Using the {\agd} update, we can further write
\begin{align}
  \Bigg(\begin{array}{c}
    \tfrac{1}{\sqrt{\mu}} A (x^{k + 1} - x^{\star})\\
    A^{1 / 2} (z^{k + 1} - x^{\star})
  \end{array}\Bigg) ={} & \Bigg(\begin{array}{c}
    \tfrac{1}{\sqrt{\mu}} A ( I - \tfrac{1}{L} A ) (
    \tfrac{\sqrt{\kappa}}{\sqrt{\kappa} + 1} (x^k - x^{\star}) +
    \tfrac{1}{\sqrt{\kappa} + 1} (z^k - x^{\star}) )\\
    A^{1 / 2} [ \tfrac{1}{\sqrt{\kappa} + 1} ( I - \tfrac{1}{\mu} A
    ) (x^k - x^{\star}) + \tfrac{\sqrt{\kappa}}{\sqrt{\kappa} + 1}
    ( I - \tfrac{1}{L} A ) (z^k - x^{\star}) ]
  \end{array}\Bigg) = M \Bigg(\begin{array}{c}
    \tfrac{1}{\sqrt{\mu}} A (x^k - x^{\star})\\
    A^{1 / 2} (z^k - x^{\star})
  \end{array}\Bigg) \nonumber
\end{align}
where $M = \Bigg(\begin{array}{cc}
  ( I - \tfrac{1}{L} A ) \tfrac{\sqrt{\kappa}}{\sqrt{\kappa} + 1} &
  \tfrac{1}{\sqrt{\mu}} ( A^{1 / 2} - \tfrac{1}{L} A^{3 / 2} )
  \tfrac{1}{\sqrt{\kappa} + 1}\\
  \sqrt{\mu} ( A^{- 1 / 2} - \tfrac{1}{\mu} A^{1 / 2} )
  \tfrac{1}{\sqrt{\kappa} + 1} & ( I - \tfrac{1}{L} A )
  \tfrac{\sqrt{\kappa}}{\sqrt{\kappa} + 1}
\end{array}\Bigg)$. Therefore,
\[ \varphi (x^{k + 1}, z^{k + 1}) = \tfrac{1}{2} \Bigg\| M
   \Bigg(\begin{array}{c}
     \tfrac{1}{\sqrt{\mu}} A x^k\\
     A^{1 / 2} z^k
   \end{array}\Bigg) \Bigg\|^2 \leq \tfrac{1}{2} \| M^{\top} M \| \Bigg\|
   \Bigg(\begin{array}{c}
     \tfrac{1}{\sqrt{\mu}} A x^k\\
     A^{1 / 2} z^k
   \end{array}\Bigg) \Bigg\|^2 = \| M^{\top} M \| \cdot \varphi (x^k, z^k) \]
and it suffices to analyze $\| M^{\top} M \|$. Since $A \in \mathbb{S}^n_{+
+}$, without loss of generality, let $A = Q \Lambda^2 Q^{\top}$ and
\begin{align}
  M ={} & \Bigg(\begin{array}{cc}
    ( I - \tfrac{1}{L} A ) \tfrac{\sqrt{\kappa}}{\sqrt{\kappa} + 1}
    & \sqrt{\mu} ( A^{- 1 / 2} - \tfrac{1}{\mu} A^{1 / 2} )
    \tfrac{1}{\sqrt{\kappa} + 1}\\
    \tfrac{1}{\sqrt{\mu}} ( A^{1 / 2} - \tfrac{1}{L} A^{3 / 2} )
    \tfrac{1}{\sqrt{\kappa} + 1} & ( I - \tfrac{1}{L} A )
    \tfrac{\sqrt{\kappa}}{\sqrt{\kappa} + 1}
  \end{array}\Bigg) \nonumber\\
  ={} & \Bigg(\begin{array}{cc}
    Q ( I - \tfrac{1}{L} \Lambda^2 ) Q^{\top}
    \tfrac{\sqrt{\kappa}}{\sqrt{\kappa} + 1} & Q ( \sqrt{\mu} \Lambda^{-
    1} - \tfrac{1}{\sqrt{\mu}} \Lambda ) Q^{\top}
    \tfrac{1}{\sqrt{\kappa} + 1}\\
    \tfrac{1}{\sqrt{\mu}} Q ( \Lambda - \tfrac{1}{L} \Lambda^3 )
    Q^{\top} \tfrac{1}{\sqrt{\kappa} + 1} & Q ( I - \tfrac{1}{L}
    \Lambda^2 ) Q^{\top} \tfrac{\sqrt{\kappa}}{\sqrt{\kappa} + 1}
  \end{array}\Bigg) \nonumber\\
  ={} & \bigg(\begin{array}{cc}
    Q & \\
    & Q
  \end{array}\bigg) \Bigg(\begin{array}{cc}
    ( I - \tfrac{1}{L} \Lambda^2 )
    \tfrac{\sqrt{\kappa}}{\sqrt{\kappa} + 1} & ( \sqrt{\mu} \Lambda^{- 1}
    - \tfrac{1}{\sqrt{\mu}} \Lambda ) \tfrac{1}{\sqrt{\kappa} + 1}\\
    \tfrac{1}{\sqrt{\mu}} ( \Lambda - \tfrac{1}{L} \Lambda^3 )
    \tfrac{1}{\sqrt{\kappa} + 1} & ( I - \tfrac{1}{L} \Lambda^2 )
    \tfrac{\sqrt{\kappa}}{\sqrt{\kappa} + 1}
  \end{array}\Bigg) \bigg(\begin{array}{cc}
    Q & \\
    & Q
  \end{array}\bigg)^{\top} \nonumber\\
  ={} & \bigg(\begin{array}{cc}
    Q & \\
    & Q
  \end{array}\bigg) N \bigg(\begin{array}{cc}
    Q & \\
    & Q
  \end{array}\bigg)^{\top}, \nonumber
\end{align}

where $N = \Bigg(\begin{array}{cc}
  ( I - \tfrac{1}{L} \Lambda^2 )
  \tfrac{\sqrt{\kappa}}{\sqrt{\kappa} + 1} & ( \sqrt{\mu} \Lambda^{- 1} -
  \tfrac{1}{\sqrt{\mu}} \Lambda ) \tfrac{1}{\sqrt{\kappa} + 1}\\
  \tfrac{1}{\sqrt{\mu}} ( \Lambda - \tfrac{1}{L} \Lambda^3 )
  \tfrac{1}{\sqrt{\kappa} + 1} & ( I - \tfrac{1}{L} \Lambda^2 )
  \tfrac{\sqrt{\kappa}}{\sqrt{\kappa} + 1}
\end{array}\Bigg)$ and it suffices to analyze $\| N^{\top} N \| = \| M^{\top}
M \|$. By direct computation, $N^{\top} N$ is a 2 by 2 block matrix
\[ N^{\top} N = \frac{1}{( \sqrt{\kappa} + 1 )^2}
   \Bigg(\begin{array}{cc}
     \tfrac{L}{\mu} - \tfrac{1}{\mu} \Lambda^2 - \tfrac{1}{L \mu} \Lambda^4 +
     \tfrac{1}{\mu} \tfrac{1}{L^2} \Lambda^6 & \sqrt{L} ( \Lambda^{- 1} -
     \tfrac{1}{L} \Lambda - \tfrac{1}{L \mu} \Lambda^3 + \tfrac{1}{\mu L^2}
     \Lambda^5 )\\
     \sqrt{L} ( \Lambda^{- 1} - \tfrac{1}{L} \Lambda - \tfrac{1}{L \mu}
     \Lambda^3 + \tfrac{1}{\mu L^2} \Lambda^5 ) & \mu \Lambda^{- 2} +
     ( \tfrac{L}{\mu} - 2 ) I - \tfrac{1}{\mu} \Lambda^2 +
     \tfrac{1}{L \mu} \Lambda^4
   \end{array}\Bigg), \]
where each block is a diagonal matrix. Since permutation does not affect the
spectrum, we can permute $N^{\top} N$ to obtain multiple 2 by 2 matrices, each
taking form of
\[ \frac{1}{( \sqrt{\kappa} + 1 )^2} \Bigg(\begin{array}{cc}
     \tfrac{L}{\mu} - \tfrac{1}{\mu} \lambda^2 - \tfrac{1}{L \mu} \lambda^4 +
     \tfrac{1}{\mu L^2} \lambda^6 & \sqrt{L} \lambda^{- 1} -
     \tfrac{1}{\sqrt{L}} \lambda - \tfrac{1}{\sqrt{L} \mu} \lambda^3 +
     \tfrac{1}{\mu L^{3 / 2}} \lambda^5\\
     \sqrt{L} \lambda^{- 1} - \tfrac{1}{\sqrt{L}} \lambda - \tfrac{1}{\sqrt{L}
     \mu} \lambda^3 + \tfrac{1}{\mu L^{3 / 2}} \lambda^5 & \mu \lambda^{- 2} +
     \tfrac{L}{\mu} - 2 - \tfrac{1}{\mu} \lambda^2 + \tfrac{1}{L \mu}
     \lambda^4
   \end{array}\Bigg), \]
where $\lambda$ is one of the diagonal elements of $\Lambda$. Hence it
suffices to compute
\begin{align}
  & \max_{\sqrt{\mu} \leq \lambda \leq \sqrt{L}}  \Bigg \|
  \small{\Bigg(\begin{array}{cc}
    \tfrac{L}{\mu} - \tfrac{1}{\mu} \lambda^2 - \tfrac{1}{L \mu} \lambda^4 +
    \tfrac{1}{\mu L^2} \lambda^6 & \sqrt{L} \lambda^{- 1} -
    \tfrac{1}{\sqrt{L}} \lambda - \tfrac{1}{\sqrt{L} \mu} \lambda^3 +
    \tfrac{1}{\mu L^{3 / 2}} \lambda^5\\
    \sqrt{L} \lambda^{- 1} - \tfrac{1}{\sqrt{L}} \lambda - \tfrac{1}{\sqrt{L}
    \mu} \lambda^3 + \tfrac{1}{\mu L^{3 / 2}} \lambda^5 & \mu \lambda^{- 2} +
    \tfrac{L}{\mu} - 2 - \tfrac{1}{\mu} \lambda^2 + \tfrac{1}{L \mu} \lambda^4
  \end{array}\Bigg)} \Bigg \| \nonumber\\
  ={} & \max_{1 \leq \lambda \leq \sqrt{\kappa}}  \Bigg \|
  \small{\Bigg(\begin{array}{cc}
    \tfrac{L}{\mu} - \lambda^2 - \tfrac{\mu}{L} \lambda^4 + \tfrac{\mu^2}{L^2}
    \lambda^6 & \sqrt{\frac{L}{\mu}} \lambda^{- 1} - \sqrt{\frac{\mu}{L}}
    \lambda - \sqrt{\frac{\mu}{L}} \lambda^3 + ( \frac{\mu}{L} )^{3
    / 2} \lambda^5\\
    \sqrt{\frac{L}{\mu}} \lambda^{- 1} - \sqrt{\frac{\mu}{L}} \lambda -
    \sqrt{\frac{\mu}{L}} \lambda^3 + ( \frac{\mu}{L} )^{3 / 2}
    \lambda^5 & \lambda^{- 2} + \tfrac{L}{\mu} - 2 - \lambda^2 +
    \tfrac{\mu}{L} \lambda^4
  \end{array}\Bigg)} \Bigg \| \nonumber\\
  ={} & \max_{1 \leq \lambda \leq \sqrt{\kappa}}  \Bigg \|
  \small{\Bigg(\begin{array}{cc}
    \kappa - \lambda^2 - \tfrac{1}{\kappa} \lambda^4 + \tfrac{1}{\kappa^2}
    \lambda^6 & \sqrt{\kappa} \lambda^{- 1} - \kappa^{- 1 / 2} \lambda -
    \sqrt{\kappa^{- 1 / 2}} \lambda^3 + \kappa^{- 3 / 2} \lambda^5\\
    \sqrt{\kappa} \lambda^{- 1} - \kappa^{- 1 / 2} \lambda - \sqrt{\kappa^{- 1
    / 2}} \lambda^3 + \kappa^{- 3 / 2} \lambda^5 & \lambda^{- 2} + \kappa - 2
    - \lambda^2 + \kappa^{- 1} \lambda^4
  \end{array}\Bigg)} \Bigg \| \nonumber\\
  ={} & \max_{1 \leq \lambda \leq \sqrt{\kappa}}  \Bigg \|
  \small{\Bigg(\begin{array}{cc}
    \kappa - \lambda^2 - \tfrac{1}{\kappa} \lambda^4 + \tfrac{1}{\kappa^2}
    \lambda^6 & \sqrt{\kappa} \lambda^{- 1} - \kappa^{- 1 / 2} \lambda -
    \sqrt{\kappa^{- 1 / 2}} \lambda^3 + \kappa^{- 3 / 2} \lambda^5\\
    \sqrt{\kappa} \lambda^{- 1} - \kappa^{- 1 / 2} \lambda - \sqrt{\kappa^{- 1
    / 2}} \lambda^3 + \kappa^{- 3 / 2} \lambda^5 & \lambda^{- 2} + \kappa - 2
    - \lambda^2 + \kappa^{- 1} \lambda^4
  \end{array}\Bigg)} \Bigg \| \nonumber\\
  ={} & \max_{1 \leq \lambda \leq \sqrt{\kappa}}  \Bigg \|
  \small{\Bigg(\begin{array}{cc}
    \frac{(\kappa - \lambda^2)^2 (\kappa + \lambda^2)}{\kappa^2} &
    \frac{(\kappa - \lambda^2) (\kappa - \lambda^4)}{\kappa^{3 / 2} \lambda}\\
    \frac{(\kappa - \lambda^2) (\kappa - \lambda^4)}{\kappa^{3 / 2} \lambda} &
    \frac{\lambda^4}{\kappa} + \kappa - \lambda^2 + \frac{1}{\lambda^2} - 2
  \end{array}\Bigg)} \Bigg \| \nonumber
\end{align}
Two eigenvalues of the matrix are respectively
\[ \tfrac{1}{2} (\tfrac{\lambda^6}{\kappa^2} \pm (\lambda^4 - \kappa) \tfrac{
   \sqrt{4 (\kappa - 1) \kappa^2 \lambda^2 + \kappa^2 + 2 \kappa (1 - 2
   \kappa) \lambda^4 + \lambda^8}}{\kappa^2 \lambda^2} + 2 \kappa - 2
   \lambda^2 + \tfrac{1}{\lambda^2} - 2), \]
whose maximizers are respectively $\sqrt{\kappa}$ and 1, giving
\[ \max \{ - 2 + \tfrac{1}{\kappa} + \kappa, \tfrac{(\kappa -
   1)^2}{\kappa } ( 1 + \tfrac{1}{2
   \kappa} + \tfrac{\sqrt{4 \kappa + 1}}{2 \kappa} ) \} \leq
   \tfrac{(\kappa - 1)^2}{\kappa} ( 1
   + \tfrac{1}{2 \kappa} + \tfrac{\sqrt{4 \kappa + 1}}{2 \kappa} ) \]
and multiplying the upper bound by $\frac{1}{( \sqrt{\kappa} + 1 )^2}$ completes the proof.

\begin{figure}[h]
\centering
\includegraphics[scale=0.4]{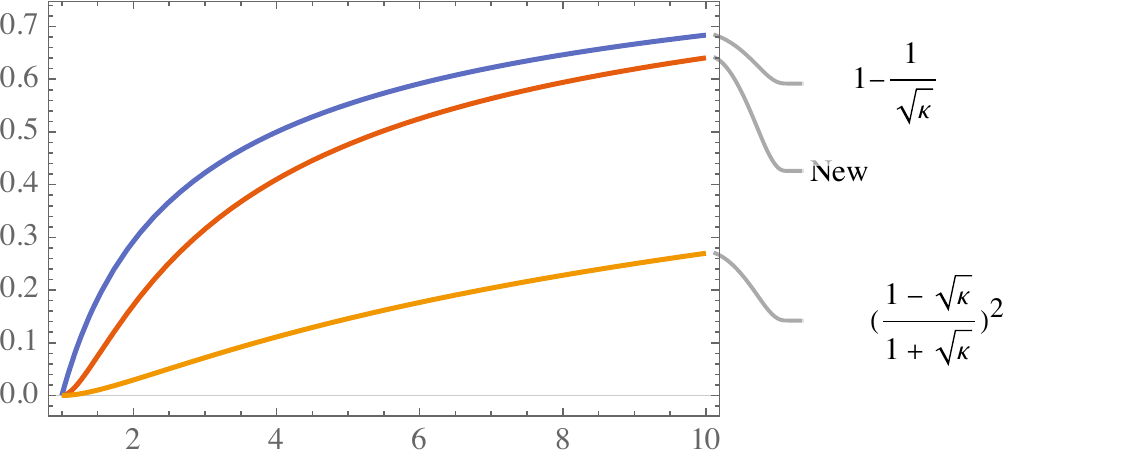}
\includegraphics[scale=0.23]{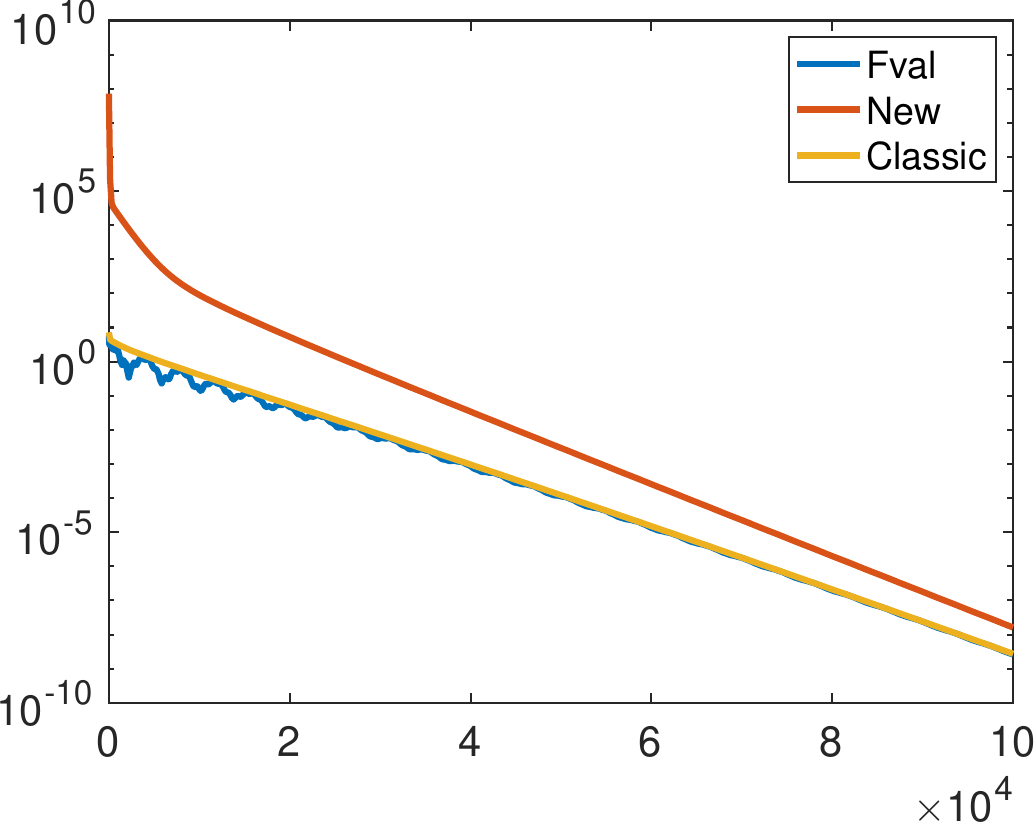}
\caption{Left: comparison between the convergence rate of the new potential function and classic ones \label{fig:pot}. Right: plots of the new potential function, the classic potential function, and the function value gap}
\end{figure}

\begin{rem}
As a final remark, we compare the linear convergence rate for the new potential function with two standard convergence rates in the literature: $1 - \frac{1}{\sqrt{\kappa}}$ for {\agd} and $(\frac{1 - \sqrt{\kappa}}{1 + \sqrt{\kappa}})^2$ for Chebeshev iteration \cite{d2021acceleration}. As \Cref{fig:pot} (left) shows, our new convergence rate is sandwiched between these two rates. \Cref{fig:pot} (right) verifies the validity of this new potential function on a toy quadratic minimization problem.
\end{rem}

%% file: sec_proxgrad.tex
In this section, we illustrate the convergence analysis of a variant of
{\oh}. We make the following assumptions throughout the section.

\begin{enumerate}[leftmargin=25pt,label=\textbf{A\arabic*:},ref=\rm{\textbf{A\arabic*}},start=1]
\item $f$ is $L$-smooth and convex with $\tmop{dom} f =\mathbb{R}^n$. \label{A3}
\item The proximal term $w$ is convex. \label{A4}
\end{enumerate}
Similarly, we denote $\varphi^\star \assign \varphi(x^\star)$. The challenging part of the analysis is that the proximal term $w$ will
introduce non-Lipschitzness of the feedback function if $\| \mathcal{G}_L (x) \| \rightarrow 0$. To
this end, we consider \Cref{alg:prox-osgm}.

\vspace{10pt}
\begin{algorithm}[H]
{\textbf{input} initial point $x^1$, initial stepsize $P_1 \in \Rbb^{n \times n}$, online proximal regularizer $\eta > 0$}

\For{$k = 1, 2, \dots$}{
$P_{k + 1} ={} \argmin_P \big\{ \tfrac{f (x^k - P_k \mathcal{G}_L (x^k)) +
  \langle \nabla f (x^k - P_k \mathcal{G}_L (x^k)) \mathcal{G}_L (x^k)^{\top}, P -
  P_k \rangle + \omega (x^k - P\mathcal{G}_{L} (x^k))}{\|
  \mathcal{G}_{L} (x) \|^2} + \tfrac{1}{2 \eta} \| P - P_k \|_F^2 \big\}$
  
Choose $x^{k+1}$ such that
  $\varphi (x^{k + 1}) \leq{} \min \{ \varphi (x^k - P_{k + 1} \mathcal{G}_L (x^k)), \varphi
  (x^k) \}$
}
\caption{{\oh} with proximal gradient  update \label{alg:prox-osgm}}
\end{algorithm}
\vspace{10pt}

This update corresponds to {\hdmclassic} we discussed in \Cref{sec:prescient}. As a first step, we define progress $h_k \assign \tfrac{\varphi (x^{k + 1}) - \varphi
(x^k)}{\| \mathcal{G}_L (x^k) \|^2}$. Note that $h_k \leq 0$ by definition,
which enables the following progress reduction.

\begin{lem}[Nonsmooth reduction] \label{lem:nonsmooth-reduction}
  Under \ref{A3} and \ref{A4}, $\min_{1 \leq k \leq K}  \|
  \mathcal{G}_L (x^k) \|^2 \leq \tfrac{\varphi (x^1) - \varphi^{\star}}{K}
  \tfrac{1}{\frac{1}{K} \sum_{k = 1}^K - h_k}$.
\end{lem}

\begin{proof}
The proof resembles that of the hypergradient reduction.
\begin{align}
  ( \textstyle \sum_{k = 1}^K - h_k ) \cdot \displaystyle \min_{1 \leq k \leq K} \|
   \mathcal{G}_L (x^k) \|^2 
  \leq{} & \textstyle \sum_{k = 1}^K - h_k \|  \mathcal{G}_L (x^k) \|^2 \tag{by $-h_k \geq 0$}\\
  ={} & - [ \textstyle \sum_{k = 1}^K \varphi (x^{k + 1}) - \varphi (x^k) ] \tag{by definition of $h_k$}\\
  ={} & \varphi (x^1) - \varphi (x^{K + 1}) \nonumber\\
  \leq{} & \varphi (x^1) - \varphi^\star  . \tag{by $f(x^{K+1}) \leq f^\star$}
\end{align}
\end{proof}
An important step is to establish online learning guarantees with
respect to $\sum_{k = 1}^K h_{x^k} (P_{k + 1})$. Below, we prove
the result for a more general setting.

\begin{lem}[Prescient online proximal gradient]
  \label{lem:proxgrad-regret}Suppose $\{ (\ell_k, \pi_k) \}$ are a sequence of
  convex functions and $\{ \ell_k \}$ are $L$-smooth. Consider the online
  proximal gradient update
  \[ P_{k + 1} = \argmin_{P \in \Pcal}  ~ \{ \ell_k (P_k) + \langle \nabla \ell_k
     (P_k), P - P_k \rangle + \pi_k (P) + \tfrac{1}{2 \eta} \| P - P_k \|^2_F
     \} . \]
  Then with $\eta \leq \frac{1}{L}$, $\textstyle \sum_{k = 1}^K \ell_k (P_{k + 1}) + \pi_k (P_{k + 1}) - [\ell_k (\hat{P}) +
     \pi_k (\hat{P})] \leq \frac{1}{2 \eta} \| P_1 - \hat{P} \|^2_F$
  for any $\hat{P} \in \Pcal$.
\end{lem}

\begin{proof}
By the optimality condition of the online update, for any $\hat{P} \in \Pcal$,
\begin{align}
  & \langle \nabla \ell_k (P_k), P_{k + 1} - P_k \rangle + \pi_k (P_{k + 1})
  + \tfrac{1}{2 \eta} \| P_{k + 1} - P_k \|^2_F \nonumber\\
  \leq{} & \langle \nabla \ell_k (P_k), \hat{P} - P_k \rangle + \pi_k (\hat{P}) +
  \tfrac{1}{2 \eta} \| \hat{P} - P_k \|^2_F - \tfrac{1}{2 \eta} \| P_{k + 1} - \hat{P}
  \|^2_F . \nonumber
\end{align}

Re-arranging the terms,
\begin{align}
\| P_{k + 1} - \hat{P} \|^2_F   \leq{} & \| P_k - \hat{P} \|^2_F - \| P_{k + 1} - P_k \|^2_F + 2 \eta \langle \nabla
  \ell_k (P_k), \hat{P} - P_{k + 1} \rangle + 2 \eta [\pi_k (\hat{P}) - \pi_k (P_{k+1})] \nonumber\\
  ={} & \| P_k - \hat{P} \|^2_F - \| P_{k + 1} - P_k \|^2_F + 2 \eta \langle \nabla
  \ell_k (P_k), \hat{P} - P_k \rangle \nonumber \\
  & \quad + 2 \eta \langle \nabla \ell_k (P_k), P_k -
  P_{k + 1} \rangle + 2 \eta [\pi_k (\hat{P}) - \pi_k (P_{k + 1})] 
   \nonumber\\
  \leq{} & \| P_k - \hat{P} \|^2_F - \| P_{k + 1} - P_k \|^2_F - 2 \eta [\ell_k (P_k) -
  \ell_k (\hat{P})] \nonumber \\
  & \quad + 2 \eta \langle \nabla \ell_k (P_k), P_k - P_{k + 1} \rangle
  + 2 \eta [\pi_k (\hat{P}) - \pi_k (P_{k + 1})] \tag{Convexity of $\ell_k$}\\
  ={} & \| P_k - \hat{P} \|^2_F - 2 \eta [\ell_k (P_{k + 1}) + \pi_k (P_{k + 1}) -
  \ell_k (\hat{P}) - \pi_k (\hat{P})] \nonumber\\
  & + 2 \eta [ \ell_k (P_{k + 1}) - \ell_k (P_k) + \langle \nabla
  \ell_k (P_k), P_k - P_{k + 1} \rangle - \tfrac{1}{2 \eta} \| P_{k + 1} -
  P_k \|^2_F ] . \nonumber
\end{align}

When $\eta \leq \tfrac{1}{L}$, we have, by smoothness of $\ell_k$, that
\[ \ell_k (P_{k + 1}) - [\ell_k (P_k) + \langle \nabla \ell_k (P_k), P_{k + 1}
   - P_k \rangle] \leq \tfrac{L}{2} \| P_{k + 1} - P_k \|^2_F \leq \tfrac{1}{2
   \eta} \| P_{k + 1} - P_k \|^2_F \]
and $\| P_{k + 1} - \hat{P} \|^2_F \leq \| P_k - \hat{P} \|^2_F - 2 \eta [\ell_k (P_{k + 1}) + \pi_k(P_{k+1})
- \ell_k (\hat{P}) - \pi_k(\hat{P})]$. Telescoping completes the proof.
\end{proof}

Chaining \Cref{lem:nonsmooth-reduction}, \Cref{lem:proxgrad-regret} and the fact $h_x(\frac{1}{L}I) \leq -\frac{1}{2L}$ \cite[Corollary 10.18]{beck2017first} gives the desired result.

\begin{thm}[\Cref{thm:proxgrad}] \label{thm:proxgrad-formal}
Suppose $f$ is $L$-smooth convex and $w$ is convex. Then for any benchmark stepsize $\hat{P} \in \Pcal$, \Cref{alg:prox-osgm} achieves 
\[ \min_{1 \leq k \leq K} \| \mathcal{G}_L (x^k) \|^2\leq \tfrac{\varphi(x^1) - \varphi^\star}{K} \tfrac{1}{\frac{1}{K} \max \{\sum_{k=1}^K h_{x^k}(\hat{P}) - \frac{L}{2}\|P_1 - \hat{P}\|_F^2, 0\}}.\]
In particular, taking $P_1 = \hat{P} = \frac{1}{L} I$ yields $\min_{1 \leq k \leq K} \| \mathcal{G}_L (x^k) \|^2 \leq 
 \tfrac{2L [\varphi(x^1) - \varphi(x^\star)] }{K}$.
\end{thm}

\begin{proof}
Using \Cref{lem:proxgrad-regret}, we can bound
\[ \textstyle \sum_{k = 1}^K h_k \leq \sum_{k = 1}^K h_{x^k} ({P_{k + 1}})
   \leq \sum_{k = 1}^K h_{x^k} (\hat{P}) + \tfrac{1}{2 \eta} \| P_1 - \hat{P} \|_F^2. \] Plugging the relation back into \Cref{lem:nonsmooth-reduction}  and using $h_{x^k} ( \tfrac{1}{L} I ) \leq - \tfrac{1}{2 L}$ completes the proof.
\end{proof}